\documentclass[12pt,a4paper]{article}
\usepackage{amsmath,amssymb,amsthm,amsfonts,amscd,euscript,verbatim, t1enc, newlfont}
\usepackage{hyperref}
 \theoremstyle{plain}
\newtheorem{theo}{Theorem}[subsection]

\newtheorem{pr}[theo]{Proposition}

 \newtheorem{lem}[theo]{Lemma}
 \newtheorem{coro}[theo]{Corollary}
  \newtheorem{conj}[theo]{Conjecture}
\theoremstyle{remark}
\newtheorem{rema}[theo]{Remark}

\theoremstyle{definition}
\newtheorem{defi}[theo]{Definition}
\newtheorem*{notat}{Notation}

 \newcommand\lan{\langle}
\newcommand\ra{\rangle}
\newcommand\bl{\bigl(} \newcommand\br{\bigl)}

\newcommand\ob{^{-1}}
\newcommand\smc{{SmCor}}
\newcommand\ssc{{Shv(SmCor)}}

\newcommand\dmk{{D^-(Shv(SmCor))}}

\newcommand\hk{\mathfrak{H}}

\newcommand\dmge{DM^{eff}_{gm}{}}
\newcommand\dmgm{DM_{gm}}
\newcommand\dme{DM_-^{eff}}
\newcommand\dmeq{DM_\q^{eff}}

\newcommand\dms{DM^{s}{}}
\newcommand\hsing{H^{sing}}
\newcommand\hsin{H_{sing}}
\newcommand\gr{gr_F}

\newcommand\mg{M_{gm}}
\newcommand\obj{Obj}
\newcommand\mo{Mor}
\newcommand\id{id}
\newcommand\cu{\underline{C}}
\newcommand\du{\underline{D}}
\newcommand\eu{\underline{E}}

\newcommand\hrt{{\underline{Ht}}}
\newcommand\hw{{\underline{Hw}}}

\newcommand\n{\mathbb{N}}
\newcommand\z{{\mathbb{Z}}}
\newcommand\q{{\mathbb{Q}}}
\newcommand\af{\mathbb{A}}
 \newcommand\h{\mathbb{H}{}}
\newcommand\p{\mathbb{P}}

\newcommand\pt{pt}

\newcommand\lamb{\Lambda}
\newcommand\al{\alpha}
\newcommand\om{\omega}
\newcommand\si{\sigma}

\newcommand\ns{\{0\}}
\DeclareMathOperator\homm{\operatorname{Hom}}
\DeclareMathOperator\prli{\varprojlim}
\DeclareMathOperator\inli{\varinjlim}

\newcommand\ihom{{\underline{Hom}}}
\newcommand\chow{Chow}
\newcommand\chowe{Chow^{eff}}

\newcommand\tcho{t_{Chow}}

\newcommand\cho{Corr_{rat}}
\newcommand\ab{Ab}
\newcommand\var{Var}
\newcommand\sv{SmVar}
\newcommand\spv{SmPrVar}

\newcommand\card{Card}
\newcommand\nde{\triangleleft}

 \DeclareMathOperator\ke{\operatorname{Ker}}
 \DeclareMathOperator\cok{\operatorname{Coker}}
\DeclareMathOperator\imm{\operatorname{Im}}
\DeclareMathOperator\co{\operatorname{Cone}}
\DeclareMathOperator\prt{Pre-Tr}
\DeclareMathOperator\prtp{Pre-Tr^+}

\DeclareMathOperator\en{\operatorname{End}}
\DeclareMathOperator\cha{\operatorname{char}}
\DeclareMathOperator\adfu{\operatorname{AddFun}}

\newcommand\trp{Tr^+}
\newcommand\tr{Tr}
\newcommand\w{{\mathfrak{w}}}

\newcommand\hml{H}
\newcommand\hel{H_{et}}
\newcommand\het{{\mathcal{H}}^{et}_{l^n}{}}

\newcommand\gds{\mathfrak{D}_s}

\begin{document}

 \title{
Weight structures vs.  $t$-structures; weight filtrations,
 spectral sequences, and complexes (for motives and in general)\\
 {\small\it{It is a great pleasure for me to dedicate this paper to Andrei Alexandrovich Suslin on the occasion of his sixtieth birthday.}}
 }
 
 \author{M.V. Bondarko
   \thanks{ 
The author acknowledges the support by  RFBR
(grants no.
 08-01-00777a and 10-01-00287a), by a Saint-Petersburg State University research grant no. 6.38.75.2011, and by the Federal Targeted Programme "Scientific and Scientific-Pedagogical Personnel of the Innovative Russia in 2009-2013" (Contract No. 2010-1.1-111-128-033).} }
 \maketitle
\begin{abstract}
In this paper we introduce a new notion of {\it weight structure}
($w$) for a triangulated category $\cu$; this notion is an important
natural counterpart of the notion of  $t$-structure.
 It allows 
extending several results of the preceding paper \cite{mymot} to a
large class of triangulated categories and functors.

The {\it heart} of $w$ is an additive category $\hw\subset \cu$.
 We prove that a weight
structure yields Postnikov towers for any $X\in \obj\cu$ (whose
'factors'  $X^i\in \obj \hw$).  For any (co)homological
functor $H:\cu\to A$  ($A$ is abelian) such a tower yields a  {\it weight
spectral sequence} $T:H(X^i[j])\implies H(X[i+j])$;  $T$ is canonical
and functorial in $X$ starting from $E_2$. $T$
specializes to the usual (Deligne) weight spectral sequences for 'classical'
realizations of Voevodsky's motives $\dmge$ (if we consider
 $w=w_{Chow}$ with $\hw=\chowe$) and to  Atiyah-Hirzebruch
spectral sequences in topology.

We prove that there often exists an exact conservative {weight complex}
functor
$\cu\to K(\hw)$. This is a generalization of the functor
$t:\dmge\to K^b(\chowe)$ constructed in \cite{mymot} (which is an
extension of the weight complex of Gillet and Soul\'e).
We prove that $K_0(\cu)\cong K_0(\hw)$ under certain restrictions.

We also introduce the concept of {\it adjacent} structures: a $t$-structure
  is adjacent 
  to $w$ if their negative parts coincide.
  This
is the case for the Postnikov $t$-structure
for the stable homotopy category $SH$ (of topological spectra) and a certain weight structure for
it that corresponds to the cellular filtration.
 We also define a new ({\it Chow}) $t$-structure $\tcho$ for
$\dme\supset \dmge$ which is adjacent to the Chow weight structure.
 We have $\hrt_{Chow}\cong \adfu (\chowe{}^{op},\ab)$;
$\tcho$ is
related to unramified cohomology. Functors left adjoint to those that
  are $t$-exact with respect to some $t$-structures are weight-exact
 with respect to the corresponding adjacent weight structures,
and vice versa. Adjacent structures
identify the following spectral sequences converging to $\cu(X,Y)$: the
one that comes from weight truncations of $X$ and the one
coming from $t$-truncations of $Y$ (for $X,Y\in \obj \cu$). Moreover,
the philosophy of adjacent
structures allows  expressing torsion motivic cohomology of
certain motives in terms of the \'etale cohomology of their
'submotives'. This is an extension of the calculation of $E_2$ of
the corresponding coniveau spectral sequences (by Bloch and Ogus).

\end{abstract}

\tableofcontents

 \section*{Introduction}

The purpose of this paper is twofold: (i) the introduction of a
 new formalism of {\it weight structures}
 for triangulated categories; (ii) the application of the
general theory of these weight structures to Voevodsky's motives (we also
consider the stable homotopy category $SH$ of topological spectra).

The notion of weight structure $w$ for a triangulated $\cu$
 is an important natural
counterpart to the notion of  $t$-structure; also, these
two types of structures are connected by interesting relations.
Similarly to $t$-structures, weight
structures are described in terms of $\cu^{w\le 0}$ (the
non-positive part of $\cu$) and
$\cu^{w\ge 0}$ (the non-negative part). We assume that 
any object of $\cu$ can be 'decomposed' (non-uniquely) into a
non-positive and a non-negative part; we call the corresponding
distinguished triangles {\it weight decompositions}. The
heart $\hw$ of  $w$ is defined as the intersection of these two parts.

A very simple (yet quite interesting) example  of a
weight structure is as follows:  $\cu=  K^?(B)$ (the homotopy
category of complexes over an additive $B$; $?$ is
any boundedness condition). We define  $\cu^{w\le 0}$
(resp. $\cu^{w\ge 0}$) as the class of complexes that are
homotopy equivalent to those concentrated in non-positive
(resp. non-negative) degrees. Then $\hw$ contains $B$
(it is equivalent to $B$ if the latter is pseudo-abelian).
Weight decompositions are given by stupid truncations of
complexes (so any object of $\cu$ has a large class of
non-isomorphic weight decompositions). In this case our
theory shows what functorial information may be obtained
using stupid truncation of complexes (this is a significant
amount of information that was never considered in the literature).
 Note that this example of a weight structure is 'almost universal'
 if one fixes the heart: for a 'reasonable' $\cu$ endowed with
a weight structure
one  has an exact {\it weight complex} functor $\cu\to K(\hw)$
(this comparison functor exists 'much more often' than its
$t$-structure counterpart, i.e.,  a functor $D^b(\hrt)\to\cu$ for
$\hrt$ being the heart of a $t$-structure for $\cu$); the weight
 complex is 'usually'
conservative (though it is very rarely an isomorphism).

In this example (as well as in the general case) there
are no $\cu$-morphisms of positive degree (one may also call them
positive Ext's with respect to $\cu$) between objects
 of $\hw$ (note in contrast that there are no morphisms of
{\bf negative} degree between objects of the heart of a $t$-structure);
 we will call
 such subcategories of $\cu$ {\it negative}. In fact, if $w$
 is {\it bounded} (as in the case of $K^b(B)$)
 then it can be easily recovered from $\hw$.
Moreover, any negative additive $B\subset \cu$ 
uniquely extends to a weight structure $w$ such that $B\subset \hw$ (at least) in the
case when $\cu$ is generated by $B$ as a triangulated category.

This observation allows  constructing
weight structures in the  following two important (and more complicated)
cases. The first case is $\cu=\dmge$ (Voevodsky's  category
of effective geometric motives), $\hw=\chowe$ (effective Chow motives);
 we call the corresponding weight structure $w_{Chow}$ the
{\it Chow weight structure}.  In this case our theory yields  canonical $\dmge$-functorial {\it (Chow)-weight spectral sequences}
$T(H,X)$ for any (co)homological functor $H:\dmge\to A$ (and $X\in \obj \dmge$). $T$
relates the cohomology of Voevodsky motives to that of Chow
motives; it is a vast (and 'motivic') generalization of the
usual (i.e., Deligne) weight spectral sequences (for \'etale
and singular cohomology of varieties). In particular, we obtain
certain (Chow)-weight spectral sequences for motivic cohomology.  We also have
(as was said above) an exact conservative weight complex functor
 $t:\dmge\to K^b(\chowe)$ (it extends to a functor $\dmgm\to K^b(\chow)$;
 the restriction of $t$ to $\chowe\subset \dmge$ is the obvious
embedding $\chowe\to K^b(\chowe)$).

Note that it is traditionally expected that the weight filtration
for singular and \'etale cohomology of motives is induced by a certain
 {\it weight filtration} of $\dmge$ (or $\dmge\q$). We briefly
explain the relation of this (conjectural) picture to our
results. It is conjectured that singular and \'etale cohomology
 for $\dmge$ can be factored through the category $MM$
(of mixed motives) which is the heart of a certain   motivic
 $t$-structure for $\dmge$. It is easily seen that the restriction
 of $w_{Chow}$ to $MM$  'should yield' the restriction of the weight
 filtration to $MM$. This situation, along with its  non-conjectural
(and compatible) analogues  for $1$-motives (over a smooth base), graded polarizable mixed Hodge
structures (lying inside their bounded derived category), and (Saito's) mixed Hodge modules 
was treated in detail in \cite{btrans} (see also \cite{bmsh} and \S\ref{stfiltr} below).  Furthemore, when the base field is a number field,
J. Wildeshaus has justified this picture for the restriction of
$w_{Chow}$ to {\it Artin-Tate motives} in the case when the base field is a number field (in \cite{wildat}).

The second case is $\cu=SH_{fin}$ with $\hw$ consisting of finite coproducts
of  sphere spectra. In this case weight decompositions are
induced by  cellular filtrations; weight spectral sequences yield
Atiyah-Hirzebruch spectral sequences!

Now we describe some other constructions and results (in the
literature) that are closely related to weight structures.

The simplest and oldest of them are projective and injective
(hyper)resolutions of
objects of (or complexes over) an abelian category $A$ and
the corresponding descriptions
of derived functors; note here that $D^?(A)$ is 'often'
isomorphic to $K^?(Proj\,A)$ (and to $K^?(Inj\,A)$). This example is related to the
 stupid weight structure of the latter category. Note  that
 stupid truncation of complexes does not yield a functor on
$K(Proj\,A)$;
 still one can describe derived functors  in terms of it.

 Another example is given by Deligne's weight filtrations and spectral
sequences (for \'etale and singular cohomology of varieties), and the definition of weights for mixed Hodge complexes 
(cf. \cite{btrans}). Note
that weight spectral sequences in \cite{de3} depend on the choice
 of nice compactifications for smooth quasi-projective varieties and
 of proper smooth hypercoverings for non-smooth projective varieties;
 Deligne's theory of weights yields the functoriality of the result
only after tensoring by $\q$. In \cite{gs} it was proved that the
$E_2$-terms of spectral sequences mentioned factor through
the weight complex functor (this term was also introduced by
 Gillet and Soul\'e);
 this allowed defining functorial weight spectral sequences with
integral coefficients. Finally, in the preceding paper \cite{mymot}
certain weight spectral sequences were defined (and described completely)
for a wide class of cohomological functors defined on $\dmge$.
Note here that the results of \cite{mymot} were deduced from the
fact that $\dmge$ has a certain (negative) differential graded
description (in terms of {\it twisted complexes} introduced in
\cite{bk}); they can be easily extended to any category that
has such a description. The description of $\dmge$ mentioned is
  somewhat similar to the definition of Hanamura's motives
(see \cite{ha3}; in \cite{mymot} it was also proved that Hanamura's
 motives are anti-isomorphic to Voevodsky ones); so it is no wonder
that some of the constructions of \cite{mymot} (including certain stupid
truncations for motives) are similar to that of \cite{ha3}.
  In the current paper we prove that weights can be introduced
for any (co)homological functor $\dmge\to A$ ($A$ is an abelian category).

 (More obvious)
predecessors of the definition of a weight structure were the
classical notions of
 'filtration bete' (see \S3.1.7 of \cite{BBD}) and of connective
 spectra (see \S7 of
 \cite{axstab}). Possibly, the so-called cofiltrations considered in
\cite{bar} are also related to the subject.
 Still,  our axiomatics and most of our main results are completely new.
 The only exception known to the author is that a part of Theorem
 \ref{madts} (the one that concerns  $t$-structures) is a slight
 generalization of Theorem 1.3 of \cite{hoshino}.
Recently weight structures were also (independently)
introduced by D. Pauksztello (he called them co-$t$-structures; see
\cite{konk}); some of the (easier) results of the
current paper were also proved there.
We would also like to mention the paper \cite{wild}, where our
results are applied to the study of
so-called boundary motives.

Half of "general" results of the current paper are extensions of the
 results of \cite{mymot}
to the case of arbitrary triangulated categories endowed with weight
 structures (that do not necessarily have a differential graded
description). We define weight spectral sequences and complexes.
In the bounded case we prove that $\cu$ is pseudo-abelian if and only if
$\hw$ is; in this case  $K_0(\cu)\cong K_0(\hw)$. We also define
 and study certain $K_0(\en \cu)$; this allows us to calculate
$K_0(\en SH_{fin})$ very explicitly.

The other half of the results of the  current paper is related to
 the yoga of
{\it adjacent structures}. This concept seems
to be completely new (though some examples for it are well-known and very
important, it seems that nobody studied this aspect of them).  We
will say that $w$ is {\it left adjacent} to a $t$-structure $t$
if $\cu^{w\le 0}=\cu^{t\le 0}$.
The
most simple (yet non-trivial) example here is the canonical $t$-structure
for $\cu=D^?(A)$  and  the (stupid) weight structure coming from the
 isomorphism $\cu\cong K^?(Proj\,A)$ (if there is one).
 Another
example is the {\it spherical} weight structure for $SH$ (mentioned
above) and the Postnikov  $t$-structure for it; in this case our
formalism implies that both of these structures yield the same
Atiyah-Hirzebruch spectral sequence converging to $SH(X,Y)$
(for $X,Y\in \obj SH$). Also, we extend $w_{Chow}$ from $\dmge$
to $\dme$ and construct  $t_{Chow}$ that is adjacent to it. $\tcho$
is related to unramified cohomology;
 $\hrt_{Chow}\cong \adfu (\chowe{}^{op},\ab)$.
One more illustration of the relevance of this notion is that functors left adjoint
 to those that
  are $t$-exact with respect to some $t$-structures are
{\it weight-exact} with respect to the corresponding
 adjacent weight structures, and vice versa.

We relate $t$-truncations
and weight truncations for adjacent
structures; this gives a collection of new interesting formulae.
 As in the partial case  $\cu=SH$ we also have the following: for $X,Y\in \obj \cu$
the spectral sequence whose $E_2$-terms  are given by the morphism groups 
from $X$ to the $t$-cohomology of $Y$, is isomorphic to the weight spectral
sequence $T(\cu(-,Y))$ (applied to $X$); see Theorem 2.6.1 of \cite{bger}. 
 The
corresponding exact couples are isomorphic also.
These spectral sequence
calculations are closely related 
 to the well-known calculations of
the $E_2$-terms of coniveau spectral sequences (by Bloch and Ogus
in \cite{blog}; see also \cite{suger}); we extend the latter and
generalize them to certain motives. This allows us to  express
 torsion motivic cohomology of these motives in terms of the \'etale
cohomology of their 'submotives' (in a certain sense).

In the next paper \cite{bger} the author constructs
a certain {\it Gersten weight structure} for a certain category of
{\it comotives} $\gds\supset \dmge$
(under a technical restriction that $k$ is countable that was lifted in \cite{bgern}; see  Remark \ref{rws}(4) and
Remark \ref{gersten} below). This result allows me
to extend the coniveau results mentioned to arbitrary motives
and makes them
more functorial; yet for general motives
it seems difficult to calculate the corresponding ('Gersten')
weight decompositions
explicitly. 

The author also develops the (general) theory of weight structures
further in \S2 of \cite{bger}. There we study
{\it orthogonal} weight and $t$-structures with respect
 to a {\it (nice) duality} $\cu^{op}\times \du\to A$;
this is a generalization of the notion of adjacent structures.
I relate weight spectral sequences and
virtual $t$-truncations of functors (introduced in the
current paper) to $t$-structures in this more general
situation (see also
Remark \ref{rsdt} below).

Now we describe the contents of the paper. More information of this sort can be
found at the beginnings  of sections.

In section \ref{wstrbdef} we give the definition of a weight
structure. We also give some
other basic definitions and prove  their (relatively) simple
properties. Our central objects of study  are {\it weight
decompositions} of objects and morphisms. We also describe certain ({\it weight})
Postnikov
towers for objects of $\cu$ that come from weight structures.

In section \ref{wspe} we define the {\it weight spectral
sequence} $T(H,X)$ (for $X\in\obj \cu$ and a (co)homological
functor $H:\cu\to A$); it comes from any choice of a weight Postnikov tower
mentioned. $T$ is canonical and functorial (in $X$) starting from $E_2$. It
specializes to the 'usual' weight spectral sequences for
'classical realizations' of Voevodsky's motives (at least  with
rational coefficients). Moreover, in this case $T$
degenerates at $E_2$, and  the $E_2^{*,*}$-terms are
exactly the graded pieces of the weight filtration of $H^*(X)$.

We also study the $D$-terms of the derived exact couple for
 $T(F,X)$ (for
a (co)homological $F$). For a homological $F$ they equal
$F_2(X)=\imm (F(w_{\ge k+1}X)
\to F(w_{\ge k}X))$
(or $F_1(X)=\imm (F(w_{\le k+1}X)\to F(w_{\le k}X))$
 for the other
possible version of the exact couple; see Remark \ref{tdss})
in our notation. $F_1$ and $F_2$
 are both (co)homological; they behave as if they were given by
truncations of $F$ in some triangulated 'category of functors' $\du$
with respect
to some $t$-structure. Composing these {\it virtual $t$-truncations}
'from different
sides' one obtains  $E_2^{*,*}(T)$.

In section \ref{swcomp} we define  the weight complex functor $t$.
Its target is a certain {\it weak category of complexes} $K_\w(\hw)$.
$K_\w(\hw)$ is a factor of $K(\hw)$ which is no longer triangulated.
 Yet the kernel of the projection $K(\hw)\to K_\w(\hw)$ is an ideal
of morphisms
whose square is zero; so our (weak) weight complex functor is not much
worse than the 'strong' one (as constructed in \cite{mymot} in the
differential graded case). In particular, $t$ is conservative, {\it
weakly exact}, and preserves the filtration given by the weight
structure. 
We conjecture that the strong
weight complex functor exists also; see Remark \ref{rctst} and
\S\ref{sftriang}. Besides, in some cases (for example, for
$SH_{fin}\subset SH$) we have $K_\w(\hw)=
K(\hw)$. Our main tool of study is the {\it weight decomposition
functor} $WD:\cu\to K_\w^{[0,1]}(\cu)$; see Theorem \ref{fwc}.

A reader only interested in motives could skip  section \ref{swcomp} since we
construct the strong version of the weight complex functor for
motivic categories
in \S\ref{dgmot}.

In  section \ref{dualh}  we relate weight structures to
$t$-structures (via the notion of {\it adjacent structures}). 
   In the case when weight and
$t$-structures are {\it adjacent}, we have a certain duality of their
hearts, and spectral sequences coming from these structures are
closely related (they differ only by a certain shift of indices). A functor left adjoint to a
$t$-exact one is {\it weight-exact} (if we consider left adjacent
weight and $t$-structures). 

We prove several results on existence of weight structures  (and also their adjacent $t$-structures).  In particular, we prove that a weight
structure can often be described in terms of some {\it negative}
additive subcategory of $\cu$.

In \S\ref{sh} we apply our results to the study of $SH$ (the stable
homotopy category). The corresponding {\it spherical} weight
structure constructed is generated by the sphere spectrum; it is
{\it left adjacent} to the (usual) Postnikov $t$-structure on $SH$.
Postnikov towers corresponding to this weight structure  are called
{\it cellular towers} (by topologists). It turns out that the corresponding 
weight complex functor calculates the singular (co)homology of spectra in
this case.

In section \ref{idemp} we prove that  a bounded $\cu$  is
idempotent complete if and only if $\hw$ is; the idempotent completion of a
general bounded $\cu$ has a weight structure whose heart is the
idempotent completion of $\hw$. If $\cu$  is bounded and idempotent
complete then $K_0(\cu)\cong K_0(\hw)$. In \S\ref{skzen} we study a
certain Grothendieck group of endomorphisms in $\cu$. Though  it
is not always isomorphic to $K_0(\en \hw)$, it is so if $\hw$ is {\it
regular} in a certain sense. Besides, we can still say something
about $K_0(\en \cu)$ in the general case also. In particular, this
allows us to generalize Theorem 3.3 of \cite{bloesn} (on
independence of $l$ for traces of certain {\it open correspondences}). 
As an application of our results, we
also calculate explicitly the groups  $K_0(SH_{fin})$ and $K_0(\en
SH_{fin})$ (along with their ring structure). We also extend these
results to the calculation of certain $K_0(\en^n SH_{fin})$ for
$n\in\n$.

In section \ref{dgmot} we translate (some of) the results of \cite{mymot}
into the language of weight structures. In particular, we show
that Voevodsky's $\dmge$ ($\subset \dmgm$) admits a {\it Chow
weight structure} whose heart it $\chowe$ (resp. $\chow$). This
allows us to prove that (Chow)-weight spectral sequences for
realizations (almost the same as those constructed in  \S7 of
\cite{mymot}, see \S\ref{twmot}) exist for all realizations and do
not depend on any choices.

In section \ref{othermot} we show that the Chow weight structure of
$\dmge$ extends to $\dme$ and admits an adjacent
$t$-structure $\tcho$ (whose heart is the category $\chowe_*=\adfu
(\chowe{}^{op},\ab)\supset\chowe$). $\tcho$ is closely related to
unramified cohomology!
 We also  prove that any possible (conjectural)  motivic
$t$-structure for $\dmge\q$ would automatically induce a canonical
weight filtration on its
heart (i.e., on mixed motives). We also prove
that (a certain version of) the weight complex functor can be
defined on $\dmge\subset \dme$ without using  resolution of singularities (so
one can define it for motives over any perfect field).

Next, we apply the philosophy of adjacent structures to the study of
coniveau spectral sequences. We express the cohomology of a motif
$X$ with coefficients in the homotopy ($t$-structure) truncations of an
arbitrary $H\in \obj\dme$ in terms of the limit of $H$-cohomology of
certain 'submotives' of $X$. 
In particular, one can express torsion motivic
cohomology of certain motives in terms of the \'etale cohomology of
their 'submotives' (this requires the Beilinson-Lichtenbaum conjecture).
As a partial case, we obtain a formula for  (torsion) motivic
cohomology with compact support of a smooth quasi-projective
variety.

In section \ref{ssupl} we show that a weight structure $w$ on $\cu$
which induces a weight structure on a triangulated $\du\subset\cu$
yields also a weight structure on the Verdier quotient $\cu/\du$. We
also prove that weight structures can be glued in a manner  similar to that for $t$-structures. This fact was used in (one of the versions of) 
the construction of the Chow  weight structures for  relative
motives (Voevodsky's motives over a base scheme $S$ that is not a field) in \S2.3 of \cite{brelmot} (and was central for the corresponding construction in \cite{bonivan}).

Next we make the following funny observation: functors represented by compositions
of $t$-truncations with respect to distinct $t$-structures can be
expressed in terms of the corresponding adjacent weight structures
(as certain images). We prove (by an argument due to A. Beilinson)
that any $f$-{\it category enhancement} of $\cu$ yields a 'strong' weight
complex functor $\cu\to K(\hw)$. We also describe other possible
sources of  'weight complex-like' functors (they are
usually conservative)  and related spectral sequences. We conclude
with the discussion of relevant types of filtrations for triangulated
categories, and of the conjectural picture  for these in the case of $\dmge\q$.

The author is deeply grateful to prof. A. Beilinson, prof. J.
Wildeshaus, prof. B. Kahn, prof.
F Deglise, prof. F. Morel, prof. S. Schwede,
prof. S. Podkorytov, and to the referees for their interesting remarks.

\begin{notat}

For categories $C,D$ we write 
$D\subset C$ if $D$ is a full 
subcategory of $C$.

 For a category $C,\ X,Y\in\obj C$, we denote by
$C(X,Y)$ the set of  $C$-morphisms from  $X$ to $Y$.
We will say that $X$ is  a {\it
retract} of $Y$ if $\id_X$ can be factored through $Y$. Note that if $C$ is triangulated or abelian 
then $X$ is a  retract of
$Y$ if and only if $X$ is its direct summand.

 For an additive $D\subset C$ the subcategory $D$ is called
{\it Karoubi-closed}
  in $C$ if it
contains all retracts  of its objects in $C$.
We will say that $D'\subset C$ is the {\it Karoubi-closure} of
$D$ if the objects of $D'$ are exactly all retracts of objects of
$D$ (in $C$).

$X\in \obj C$ will be called compact if the functor $X^*= C(X,-)$
commutes with all small coproducts that exist in $C$ (contrary
to  tradition, we do not assume that arbitrary coproducts exist).

For a category $C$ we denote by $C ^{op}$ the opposite
category.

$\cu$ will usually denote a triangulated category; usually it will
be endowed with a weight structure $w$ (see Definition \ref{dwstr}
below). We will use the term 'exact functor' for a functor of
triangulated categories (i.e.,  for a functor that preserves the
structures of triangulated categories). We will call a covariant
additive functor $\cu\to A$ for an abelian $A$ {\it homological} if
it converts distinguished triangles into long exact sequences;
homological functors $\cu^{op}\to A$ will be called {\it
cohomological} when considered as contravariant functors $\cu\to A$.

For an additive category $A$ we denote by $C(A)$ the unbounded
category of cohomological complexes over $A$; $K(A)$  is the homotopy category
of $C(A)$, i.e.,
  morphisms of complexes are considered up to homotopy
equivalence;  $C^-(A)$  denotes the category of complexes over $A$
bounded above; $C^b(A)\subset C^-(A)$ is the subcategory of
bounded complexes;  $K^b(A)$ denotes the homotopy category  of
bounded complexes. We will denote by $C(A)^{\le i}$ (resp.
$C(A)^{\ge i}$) the 
category of complexes concentrated
in degrees $\le i$ (resp. $\ge i$).

Below for a complex denoted by $\dots\to X^{-1}\stackrel{d^{-1}_X}\to X^0\stackrel{d^{0}_X} \to
X^1\to\dots$ (or similarly) we will assume that $X^i$ is in degree
$i$. For other complexes we will assume (by default) that the last
term specified is in degree $0$. Since (if $A$ is an abelian category) the functor $X\mapsto \ke(d^{i}_X)/\imm (d_X^{i-1})$ is homological, we will call it homology and denote it by $H_i(X)$ (since this convention is used in some other papers of the author; cf. Definition \ref{numbers} below). Below we will also introduce a similar convention for $t$-structures.

For an abelian $A$ we will denote by $D^b(A)\subset D^-(A) \subset D(A)$  the
corresponding versions of the derived category  of $A$.

When dealing with triangulated categories we will often  use
conventions and auxiliary statements of \cite{gelman}.
For $f\in\cu (X,Y)$, $X,Y\in\obj\cu$, we will call the third vertex
of (any) distinguished triangle $X\stackrel{f}{\to}Y\to Z$ a cone of
$f$. Recall that different choices of cones are connected by
non-unique isomorphisms (easy, see \S IV.1.7 of ibid.).
Besides, in $C(A)$ (see below) we have canonical cones of morphisms
(see section
\S III.3 of ibid.).

We will often specify a distinguished triangle by two of its
morphisms.

 For a set of
objects $C_i\in\obj\cu$, $i\in I$, we will denote by $\lan C_i\ra$
the smallest strictly full triangulated 
subcategory containing all $C_i$; for
$D\subset \cu$ we will write $\lan D\ra$ instead of $\lan C:\ C\in\obj
D\ra$. 

For $X,Y\in \obj \cu$ we will write $X\perp Y$ if $\cu(X,Y)=\ns$.
For $D,E\subset \obj \cu$ we will write $D\perp E$ if $X\perp Y$
 for all $X\in D,\ Y\in E$.
For $D\subset \cu$ we will denote by $D^\perp$ the class
$$\{Y\in \obj \cu:\ X\perp Y\ \forall X\in D\}.$$
Sometimes we will denote by $D^\perp$ the corresponding
 full subcategory of $\cu$. Dually, ${}^\perp{}D$ is the class
$\{Y\in \obj \cu:\ Y\perp X\ \forall X\in D\}$.

We will say that  $C_i$ generate $\cu$ if $\cu$ equals $\lan
C_i\ra$. We will say that $C_i$ {\it weakly generate} $\cu$ if for
$X\in\obj\cu$ we have $\cu(C_i[j],X)=\ns\ \forall i\in I,\
j\in\z\implies X=0$ (i.e., if $\{C_i[j]\}^\perp$ contains only
 zero objects).
Dually, $C_i$ {\it weakly cogenerate} $\cu$ if ${}^\perp\{C_i[j]\}=\ns$.

In this paper all complexes will be cohomological, i.e., the degree
of all differentials is $+1$; respectively, we will use
cohomological notation for their terms.

 $\ab$ is the category of abelian groups; $\ab_{fr}$ is its
 subcategory of free
abelian groups; $\ab_{fin.fr}$ is the category of finitely
generated free abelian groups.

 For additive $C,D$ we denote by $\adfu(C,D)$ the
category of additive functors from $C$ to $D$ (we will always
be able to assume that $C,D$ are small). 
For an additive $A$ we will
denote by $A^*$ the category $\adfu(A,\ab)$ and by $A_*$ the
category $\adfu(A^{op},\ab)$. Note that both of these are abelian.
Moreover, Yoneda's lemma gives full embeddings of $A$ into $A_*$ and
of $A^{op}$ into $A^*$ (these send $X\in\obj A$ to $X_*=A(-,X)$ and
to $X^*=A(X,-)$, respectively). $A_*'$ will denote the full abelian
subcategory of $A_*$  generated by $A$.

It is easily seen that any object of $A$ is projective in
$A_*$. Besides, any object of $A_*$ has a resolution by (infinite)
coproducts of objects of $A$. These statements are quite easy; the
proofs can be found at the beginning of \S8 of \cite{vbook}.

The definition of a cocompact object is dual to that of a compact
one, i.e., $X\in \obj\cu$ is cocompact if $\cu(\prod_{i\in I}
Y_i,X)=\bigoplus_{i\in I} \cu(Y_i,X)$ for any set $I$ and any
$Y_i\in\obj \cu$ such that the product exists (we do not assume
that it always exists).

We list the main definitions of this paper. Weight structures,
$\cu^{w\le 0},\cu^{w\ge 0}$, and weight decompositions of objects
are defined in \ref{dwstr};  $\hw$ (the {\it heart} of $w$),
$\cu^{w=0}$, $\cu^{w\le l}$, $\cu^{w\ge l}$, $\cu^{[j,i]}$,
non-degenerate, and bounded (above, below or both) weight structures
are defined in Definition \ref{d2}; $X^{w\le i}$, $X^{w\ge i+1}$,
 $w_{\le i}X$, and $w_{\ge i+1}X$
are defined in Remark \ref{rwd}; extension-stable subcategories are
defined in Definition \ref{exstab}; $\cu^-$, $\cu^+$, and $\cu^b$
are defined in Definition \ref{cuplus}; notation and
several definitions related to weight decomposition of morphisms, (infinite) weight
decomposition of objects, and weight Postnikov towers for objects are 
introduced in \S\ref{wedecompm};  weight filtrations on functors
are introduced in Definition \ref{dwfilf}; weight complexes of object are defined  in Definition \ref{dwcompl}; weight spectral
sequences (denoted by $T(H,X)$) are introduced in \S\S\ref{swss}-\ref{swssc}; virtual $t$-truncations of (co)homological functors are
defined in Remark \ref{rtrfun};
 $T\nde \mo A$ ($T$ is an ideal of morphisms for $A$)
is defined in Definition \ref{didmo}; $A/T$ is defined in Remark
\ref{ridmo}; the weak category of complexes $K_\w(A)$, distinguished
triangles in it, and weakly exact functors are defined in Definition
\ref{dkw}; the weight decomposition functor $WD$ and the weight
complex functor $t$ are introduced in
Theorem \ref{fwc}; $t$-structures are recalled
 in \S\ref{dtst}; countable homotopy colimits in $\cu$ and their
properties are
 described in \S\ref{colim};
 negative subcategories, 
  and
small envelopes are introduced in Definition \ref{negth}; the
categories $SH$ and $SH_{fin}$ of spectra are mentioned in
Corollary \ref{wsh}; adjacent (weight and $t$)-structures are
defined in Definition \ref{deadj}; weight-exact (and also $t$-exact)
functors are defined in Definition \ref{dadad};
negatively well-generating sets
of objects are defined in Definition
\ref{wg};
more categories of spectra, singular cohomology, and singular
homology 
of spectra are considered in \S\ref{sh}; we
discuss idempotent completions in \S \ref{ridcomp}; $K_0$-groups of
$\hw$, $\cu$, $\en \hw$,  $\en\cu$, $\en^n \hw$, and $\en^n\cu$ are
defined in \S\S\ref{sksws}--\ref{skzen}; regular additive
categories are defined in Definition \ref{dreg}; differential graded
categories and twisted complexes over them are defined in
\S\ref{sbdedg}; truncation functors $t_N$ are constructed in
\S\ref{cwecom}; the spectral sequence $S(H,X)$ is considered in
\S\ref{twmot}; we recall $\smc$, $J$, $\hk$, $\dme$, $\dms$,
$\dmge$, and $\dmgm$  in  \S\ref{remmymot}; we recall $\cho$, $\chowe$
and $\chow$ in \S\ref{hchow}; $\hml^{i}(X,\z/l^n\z(s))$ and
$\hel^{i}(X,\z/l^n\z(s))$ are considered in \S\ref{motet}; gluing data
is defined in \S\ref{sglu}; weight filtration (for motives) is mentioned
in \S\ref{stfiltr}. 

\end{notat}

\section{Weight structures for triangulated categories: basic
definitions and properties; auxiliary statements}\label{wstrbdef}

In this section we give the definition of a weight structure $w$
in a triangulated category $\cu$ (in \S\ref{sdwstr}) (this
includes the notion of  {\it weight decomposition} of an object).
We give other basic definitions and prove  certain simple
properties  of them in \S\S\ref{sotherd}--\ref{sbasicpr}. We recall
certain auxiliary statements that will help us to prove that 
weight decompositions are functorial (in a certain sense) in
\S\ref{auxss}. We study weight decompositions of morphisms 
and weight
Postnikov towers for objects 
in \S\ref{wedecompm}.

\subsection{Weight structures: definition}
\label{sdwstr}

\begin{defi}[Definition of a weight structure]\label{dwstr}

A pair of subclasses $\cu^{w\le 0},\cu^{w\ge 0}\subset\obj \cu$ for
a triangulated category $\cu$ will be said to define a weight
structure $w$ for $\cu$ if
they  satisfy the following conditions:

(i) $\cu^{w\le 0},\cu^{w\ge 0}$ are additive 
and Karoubi-closed in $\cu$
(i.e., contain all retracts of their objects that belong to
$\obj\cu$).

(ii) {\bf 'Semi-invariance' with respect to translations}.

$\cu^{w\ge 0}\subset \cu^{w\ge 0}[1]$; $\cu^{w\le 0}[1]\subset
\cu^{w\le 0}$.

(iii) {\bf Orthogonality}.

$\cu^{w\ge 0}\perp \cu^{w\le 0}[1]$ (i.e., for any
$X\in\cu^{w\ge 0}$, $Y\in \cu^{w\le 0}[1]$ we have
$\cu(X,Y)=\ns$).

(iv) {\bf Weight decompositions}.

 For any $X\in\obj \cu$ there
exists a distinguished triangle

\begin{equation}\label{wd}
B[-1]\to X\to A\stackrel{f}{\to} B
\end{equation} such that $A\in \cu^{w\le 0}, B\in \cu^{w\ge 0}$.

\end{defi}

The triangle (\ref{wd}) will be called a {\it weight
decomposition} of $X$.

The  basic example of a weight structure is given by the stupid
filtration on
$K(A)$ (for an arbitrary additive  $A$; we will call it the
{\it stupid} weight structure). We will omit $w$ in this case and denote by
$K(A)^{\le 0}$ (resp. $K(A)^{\ge 0}$) the class of complexes that are
(homotopy equivalent to complexes)
 concentrated in degrees $\le 0$ (resp. $\ge 0$). Its {\it heart} (see
Definition \ref{d2} below) is the  the Karoubi-closure  of $A$
 in $K(A)$ (so, it is equivalent to $A$ if the latter is pseudo-abelian).
  Moreover, we will see below (cf. Theorems
\ref{fwc}, \ref{wecomp}, and Remark \ref{rctst}) that this example
is 'almost universal' if one fixes the heart.

Below we will construct several more examples (for example, see Theorem \ref{recw} for  strong results on the existence of weight structures).

Note here that for $A$ being an abelian category with enough
projectives and injectives the appropriate version (i.e.,
we impose some boundedness conditions) of $D(A)$ is "often" equivalent to
$K^?(Inj\,A )$ and to $K^?(Proj\,A)$ (here $Proj\,A$ and $Inj\,A$
denote the categories of
projective and injective objects of $A$). We obtain that some
triangulated categories can support
at least two distinct weight structures with non-isomorphic hearts.

\begin{rema}\label{rdd}

1. Obviously, the axioms of weight structures are self-dual (recall
that the same is true for axioms of triangulated categories). This
means that  $(C_1,C_2)$ define a weight structure for $\cu$ if and only if
$(C_2^{op}, C_1^{op})$ define a weight structure for $\cu^{op}$.
Recall also that the same is true for $t$-structures (see Definition
\ref{dtstr}).

2. Besides, if $A$ is an abelian category, $F:\cu\to A$  a
(co)homological functor,  then the same is true for the functor
$F^{op}:\cu^{op}\to A^{op}$ obtained from $F$ in the natural way.
Hence one can interchange $\cu^{w\ge 0}$ with $\cu^{w\le 0}$ without
changing the variance of $F$.

 We will apply these observations several times.

3. A major distinction between the axioms of weight structures and those of $t$-structures (see Definition \ref{dtstr} below) is that their orthogonality conditions are opposite; also, the arrows in $t$-decompositions 'go in the opposite direction'. We will see below that this results in a drastic difference between the properties of these two types of structures. 

\end{rema}

\subsection{Other definitions}\label{sotherd}

We will also need the following definitions.

\begin{defi}\label{d2} [Other basic definitions]

1. The full category $\hw\subset \cu$ whose object class is
$\cu^{w=0}=\cu^{w\ge 0}\cap \cu^{w\le 0}$ 
 will be called the {\it heart} of the weight structure
$w$. Obviously, $\hw$ is additive.

2. $\cu^{w\ge l}$ (resp. $\cu^{w\le l}$) will denote $\cu^{w\ge
0}[-l]$ (resp. $\cu^{w\le 0}[-l]$).

3. For all $i,j\in \z$, $i\ge j$ we define $\cu^{[j,i]}=\cu^{w\ge
j}\cap \cu^{w\le i}$. By abuse of notation, we will sometimes
identify $\cu^{[j,i]}$ with the corresponding full additive
subcategory of $\cu$.

4. $w$ will be called {\it non-degenerate} if $$\cap_l \cu^{w\ge
l}=\cap_l \cu^{w\le l}=\ns.$$

5. $w$ will be called {\it bounded above} (resp. bounded below) if
$\cup_{l\in\z} \cu^{w\le l}=\obj \cu$
(resp. $\cup_{l\in\z} \cu^{w\ge l}=\obj
\cu$).

6. $w$ will be called {\it bounded} if it is bounded both above
and  below.

\end{defi}

Now we observe  an important difference between decompositions of
objects with respect to $t$-structures and 
with respect to
weight structures.

\begin{rema}\label{rwd}

1. In contrast to the $t$-structure situation, the presentation of $X$
in the form (\ref{wd}) is (almost) never canonical. The only exception
is the  {\it totally degenerate} situation when $\cu^{w=0}=\ns$.

Note that in this case the classes $\cu^{w\le i}$
coincide for all $i\in\z$;
all $\cu^{w\ge i}$ coincide also. It easily follows that $\cu^{w\le 0}$
and
$\cu^{w\ge 0}$ yield (full) triangulated subcategories of $\cu$
(cf. Proposition
\ref{bw}(3) below). In this case the weight
decomposition axiom yields that the inclusion $\cu^{w\le 0}\to \cu$
possesses a left adjoint.
Whereas such situations are certainly important, it does not
seem to make much sense to study
them using the (general) formalism of weight structures. In
particular, below we will mostly be
interested in the non-degenerate situation (since we will
'cut' objects of $\cu$ into pieces that belong to
 $\cu^{w=0}$).

2. Yet we will need to choose some $(A,B,f)$ several times. We will
write that $A=X^{w\le 0}$, $B=X^{w\ge 1}$ if there exists a
distinguished triangle (\ref{wd}). In Theorem \ref{fwc}  below we
will verify that $X\to (A,B,f)$ is a 'functor up to morphisms that are zero on
(co)homology'.

We will also often denote $(X[i])^{w\le 0}$ by $X^{w\le i}$ and
$(X[i])^{w\ge 1}$ by $X^{w\ge i+1}$ for all $i\in\z$. Note that we
have $X^{w\le i}\in  \cu^{w\le 0}$ and $X^{w\ge i+1}\in  \cu^{w\ge
0}$.

Below we will introduce a similar convention for the {\it weight
complex} of $X$.

 Besides, we will sometimes denote $X^{w\le i}[-i]$ by $w_{\le i}X$
and $X^{w\ge i}[-i]$ by $w_{\ge i}X$. So, for any $i\in\z$ we
have a distinguished triangle
 $$w_{\ge i+1}X\to X\to w_{\le i}X.$$

Yet   if $X$ admits a weight decomposition that {\it avoids weight
$0$} (a term proposed by J. Wildeshaus, see Definition 1.6 of
\cite{wild}) then the choice of such a weight decomposition for $X$
is unique; see  Remark \ref{bmorwc}(2) below.

\end{rema}

\subsection{Basic properties of weight structures} \label{sbasicpr}

We will need the following  definition several times.

\begin{defi}\label{exstab}

 $D\subset \obj \cu$ will be called extension-stable if for any
distinguished triangle $A\to B\to C$  in $\cu$ we have the following: $A,C\in
D\implies B\in D$.

We will also say that the corresponding full subcategory
is extension-stable.

\end{defi}

\begin{rema}\label{rexstab}
Certainly, any extension-stable subclass of $\obj \cu$ is additive
(i.e., closed with respect to finite direct sums) since a
triangle of the form $A\to A\oplus C\to C$ is always distinguished.
\end{rema}

For any $\cu,w$  the following basic properties are
fulfilled. Most of these properties are parallel to those of
$t$-structures; assertion 7 illustrates the distinction between these
notions.

\begin{pr}\label{bw}

1.  $\cu^{w\le 0}=(\cu^{w\ge 1})^\perp$ (see Notation).

2. Vice versa, $\cu^{w\ge 0}={}^\perp \cu^{w\le -1}$.

3. $\cu^{w\ge 0}$, $\cu^{w\le 0}$, and $\cu^{w=0}$ are
extension-stable in the sense of Definition \ref{exstab}.

4. All $\cu^{w\le i}$ are closed with respect to arbitrary (small)
products (that exist in $\cu$).

 5.  All $\cu^{w\ge i}$ are closed with respect to arbitrary
 (small) coproducts (that exist in $\cu$).

6. For any weight decomposition of $X\in \cu^{w\ge 0}$ (see
(\ref{wd})) we have $A\in \cu^{w=0}$. 

7. If $A\to B\to C\to A[1]$ is a distinguished triangle and
$A,C\in \cu^{w= 0}$ then $B\cong A\oplus C$.

8. If $X\in \cu^{w=0}$, $X[-1]\to A{\to} B$ is a weight
decomposition (of $X[-1]$) then $B\in \cu^{w=0}$; $B\cong A\oplus
X$.

\end{pr}
\begin{proof}

1. We should prove the following: if $\cu(Y,X)=\ns$ for some $X\in\obj \cu$
and all $Y\in \cu^{w\ge
1}$ then $X\in \cu^{w\le 0}$.

Let $B[-1]\to X\to A\to B$  be a weight decomposition of $X$.
Since $B[-1]\perp X$ we obtain that $X$ is a retract of $A$;
hence $X\in \cu^{w\le 0}$.

2. The proof is similar to that of assertion 1 and can be obtained by
dualization (see Remark \ref{rdd}). If $ B[-1]\to X[-1]\to A\to B$
is a weight decomposition of $X[-1]$ then $X[-1]\perp A$. Hence
$X$ is a retract of $B$.

3. Let $A,C\in \cu^{w\ge 0}$. For any $Y\in \obj \cu$ we have a
(long) exact sequence $\dots\to \cu(C,Y)\to \cu(B,Y)\to \cu(A,Y)\to
\dots$; hence by  Definition \ref{dwstr}(ii) we obtain that
$B\perp \cu^{w\le -1}$. Now
assertion 2 
implies that $B\in \cu^{w\ge 0}$.

The proof for the case $A,C\in \cu^{w\le 0}$ can be obtained by
dualization.

The statement for the case $A,C\in \cu^{w= 0}$ now follows
immediately from the definition of $\cu^{w= 0}$.

4. Obviously, assertion 1  implies that
$\cu^{w\le i}=(\cu^{w\ge i+1})^\perp$. This yields
the result immediately.

5. Similarly, by assertion 2 we have
$\cu^{w\ge i}={}^\perp\cu^{w\le i-1}$; this yields the result.

6. $A\in \cu^{w\le 0}$ by definition. Since we have a
distinguished triangle $X\to A\to B\to X[1]$,  assertion 3 implies
 that $A\in \cu^{w\ge 0}$.

7. Since $C\in \cu^{w\ge 0}$ and $A[1]\in \cu^{w\le -1}$, the
morphism $C\to A[1]$ in  the distinguished triangle is zero; so
the triangle splits.

8. We have a distinguished triangle $A\to B\to X$. By assertion 3 we
obtain that $B\in \cu^{w=0}$. Then assertion 7 yields the result.

\end{proof}

\begin{rema}\label{rthick}

1. We try to answer the questions when a morphism $b[-1]\in
\cu(B[-1],X)$ for $B\in \cu^{w\ge 0}$ extends to a weight
decomposition of $X$ and  $a\in \cu(X, A)$ for $A\in \cu^{w\le 0}$
extends to a weight decomposition of $X$ (i.e., $\co(f)\in
\cu^{w\ge 0}$) using Proposition \ref{bw}(1,2).

We apply the long exact sequence corresponding to the functor
$C^*$ for $C\in \cu^{w\ge 0}$ (resp. to $C_*$ for $C\in \cu^{w\le
0}$).  In the first case we obtain that $b[-1]$ extends to a
weight decomposition if and only if  the map $\cu(C[i],B[-1])\to \cu(C[i],X)$
induced by $b$ is bijective for $i=-2$ and is surjective for
$i=-1$ for all $C\in \cu^{w\ge 0}$. Dually, $a$ extends to a
weight decomposition if and only if
  for any
$C\in \cu^{w\le 0}$ the map $\cu(A,C)\to \cu(X,C)$ induced by $a$
is bijective for $i=1$ and is injective for $i=0$.

Moreover, in many important cases (cf. section \ref{dualh} below)
it suffices to check the conditions of part 1  (resp. part 2) of
Proposition \ref{bw} only for $Y=C[i]$ for $C\in \cu^{w=0}$, $i<0$
(resp. for $i>0$). Then these conditions are equivalent to the
bijectivity of all maps $\cu(C[i],B[-1])\to \cu(C[i],X)$ induced
by $b$ for $i<-1$ and their surjectivity for $i=-1$ for all $C\in
\cu^{w=0}$ (resp. to the bijectivity of all maps $\cu(A,C)\to
\cu(X,C)$ induced by $a$ for $i>0$ and their injectivity for
$i=0$).

We will use this observation below.

2. Certainly, all
$\cu^{w\ge i},\ \cu^{w\le i}$, and $\cu^{w=i}$ are additive.

3. Since all (co)representable functors are additive, for any class
 $C\subset\obj \cu$ the classes $C^\perp$ and ${}^\perp C$
are Karoubi-closed (in $\cu$). We will use this fact below.
 \end{rema}

\begin{defi}\label{cuplus}

We consider $\cu^-=\cup \cu^{w\le i}$ and $\cu^+=\cup \cu^{w\ge
i}$.

We call $\cu^b=\cu^+\cap \cu^-$ the class of {\it bounded} objects
of $\cu$.
\end{defi}

\begin{pr}\label{cub}
1. $\cu^-,\ \cu^+,\ \cu^b$ are
 Karoubi-closed triangulated subcategories of $\cu$.

2. $w$ induces weight structures for $\cu^-,\ \cu^+$, and $\cu^b$, whose
hearts equal $\hw$.

3. $w$ is non-degenerate when restricted to $\cu^b$.
\end{pr}
\begin{proof}
1.  From   Proposition \ref{bw}(3) we easily deduce that
$\cu^-,\ \cu^+,\ \cu^b$ are closed with respect to finite coproducts,
cones of morphisms, and retracts.

2. It suffices to verify that for any object  $X$ of $\cu^-,\
\cu^+$, or  $\cu^b$, the components of  all of its possible
weight decompositions belong to the corresponding category.

Let a distinguished triangle $B[-1]\to X\to A\to B\to X[1]$ be a
weight decomposition of $X$, i.e $A\in \cu^{w\le 0}, B\in
\cu^{w\ge 0}$.

If $X\in \cu^{w\le i}$ for some $i>0$ then  Proposition
\ref{bw}(3) implies that $B\in \cu^{w\le i-1}$. Similarly, if $X\in
\cu^{w\ge i}$ for some $i\le 0$ then $A\in \cu^{w\ge i}$. We
obtain the claim.

3. Let $X\in\obj \cu^b\bigcap (\cap \cu^{w\ge i})$; in particular,
$X\in \cu^{w\le j}$ for some $j\in\z$. Then by the orthogonality
property for $w$ we have $X\perp X$; hence $X=0$.

A similar argument proves that $\obj \cu^b\bigcap (\cap \cu^{w\le
i})=\ns$.

\end{proof}

 $\cu^b$ is especially important; note that it equals $\cu$
if $(\cu,w)$ is bounded.

Now we prove a simple lemma; it will help us several times (below) to
verify for a pair of subcategories that they satisfy axioms of weight
structures.

\begin{lem}\label{lsimple}
Let $C\subset \obj\cu$. Consider the classes $C_1= (C^\perp)[-1]$
 and $C_2=(^\perp C)[1]$.
 
1. $C_1,C_2$   are
Karoubi-closed and extension-stable.

2. Suppose that $C$ is additive, Karoubi-closed, and satisfies $C\subset C[1]$.
Suppose also that for any $X\in\obj\cu$ there exist $A\in C_1,\ B\in
C$, and a distinguished triangle $B[-1]\to X\to A{\to} B$. Then the
pair $(C_1,C)$ defines a weight structure for $\cu$.

3. Suppose  now that $C$ is additive, Karoubi-closed, and satisfie $C[1]\subset C$.
Suppose also that for any $X\in\obj\cu$ there exist $A\in C,\ B\in
C_2$ and a distinguished triangle $B[-1]\to X\to A{\to} B$. Then the
pair $(C,C_2)$ defines a weight structure for $\cu$.

4. If for $D_1,D_2\subset \obj \cu$ we have $D_2\perp D_1[1]$,
then the same is true for the
Karoubi-closures of $D_1,D_2$.

\end{lem}
\begin{proof}
1. The assertion is immediate from the fact that
(co)representable functors are additive and cohomological (resp.
homological); cf. the proof of  Proposition \ref{bw}(3).

2. By assertion 1,  it suffices to check that
$C_1[1]\subset C_1$.  Now, for any $X\in C_1$, $Y\in C$ we have
$\cu(Y,X[2])=\cu (Y[-1],X[1])=\ns$ (by the definition of $C_1$ and the equivalence of
$C[-1]\subset C$ and $ C\subset C[1]$).

3. This is exactly the dual of assertion 2 (see Remark \ref{rdd}).
\end{proof}

4. Immediate from the biadditivity of $\cu(-,-)$.

 Lastly we prove a simple statement on comparison of weight
 structures.

 \begin{lem}\label{lodn}
Suppose that $v,w$ are weight structures for $\cu$; let $\cu^{v\le
0}\subset \cu^{w\le 0}$ and $\cu^{v\ge 0}\subset \cu^{w\ge 0}$.
Then $v=w$ (i.e., the inclusions are equalities).\end{lem}
\begin{proof}

Let $X$ belong to $ \cu^{w\le 0}$; let $B[-1]\stackrel{h}{\to} X\to
A\to B$ be a weight decomposition of $X$ with respect to $v$.
Since $B[-1]\in \cu^{w\ge 1}$, the orthogonality property for $w$
implies $h=0$. Hence $X$ is a retract of $A$. Since $\cu^{v\le 0}$
is Karoubi-closed, we have  $X\in \cu^{v\le 0}$.

We obtain that $\cu^{v\le 0}= \cu^{w\le 0}$. The equality
$\cu^{v\ge 0}= \cu^{w\ge 0}$ is proved similarly.
\end{proof}

\subsection{Some auxiliary statements: 'almost functoriality'
 of distinguished
triangles}\label{auxss}

We will prove below that  weight decompositions are  functorial in
a certain sense ('up to morphisms that are zero on cohomology'). We will need some
(general) statements on 'almost functoriality' of distinguished
triangles for this. This means that a morphism  between single
vertices of two distinguished  triangles can often be completed to a
morphism of these triangles.

\begin{lem}\label{lbbd}

Let  $T:X\stackrel{a}{\to} A\stackrel{f}{\to} B\stackrel{b}{\to} X[1]$
and $T':X'\stackrel{a'}{\to} A'\stackrel{f'}{\to}
 B'\stackrel{b'}{\to} X'[1]$ be
distinguished triangles.

1.  Suppose that $B\perp A'[1]$.
 Then for any morphism $g:X\to X'$ there exist $h:A\to A'$ and
$i:B\to B'$ completing $g$ to a morphism of triangles $T\to T'$.

2. Assume moreover $B\perp A'$. Then  $g$ and $h$ are unique.

\end{lem}
\begin{proof}

This fact can be easily deduced from Proposition 1.1.9 of \cite{BBD}
(or Corollary IV.1.4 of \cite{gelman}); we use the same argument
here.

1. Since the sequence $\cu(B,A')\to \cu(B,B')\to \cu(B,X'[1])\to
\cu(B,A'[1])$ is exact, there exists $i:B\to B'$ such that
$b'\circ i=g[1]\circ b$. By axiom TR3 (see \S IV.1 of
\cite{gelman}) there also exist a morphism $h:A\to A'$ that
completes $(g,i)$ to a morphism of triangles.

2. Now we also have $\cu(B,A')=\ns$. Hence the exact sequence
mentioned in the proof of assertion 1 now also yields the uniqueness of
$i$.

The condition on $h$ is that $h\circ a =a'\circ g$.  We  have an
exact sequence $\cu(B,A') \to \cu(A,A')\to \cu(X,A')$. Since
$B\perp A'$, we obtain that $h$ is unique also.
\end{proof}

\begin{pr}\label{3na3}[$3\times 3$-Lemma]

Any commutative square
$$\begin{CD}
X@>{a}>>A\\
@VV{g}V@VV{h}V \\
X'@>{a'}>>A'
\end{CD}$$
can be completed to a  $4\times 4$ diagram (we will mainly need
its upper left $3\times 3$ part) of the following sort:
\begin{equation}\label{dia3na3}
\begin{CD}
X@>{a}>>A@>f>> B@>{}>>X[1]\\
 @VV{g}V@VV{h}V @VV{i}V@VV{g[1]}V\\
X'@>{a'}>>A'@>{f'}>>B'@>{}>>X'[1]\\
 @VV{}V@VV{}V @VV{}V@VV{}V\\
X''@>{a''}>>A''@>{f''}>>B''@>{}>>X''[1]\\
 @VV{}V@VV{}V @VV{}V@VV{}V\\
X[1]@>{a[1]}>>A[1]@>f[1]>> B[1]@>{}>>X[2]\\
\end{CD}
\end{equation}
such that all rows and columns are distinguished triangles and all
squares are commutative, except the right lowest square which
anticommutes.
\end{pr}
\begin{proof}
The proof is mostly a repetitive use of the octahedral axiom.
However it requires certain unpleasant diagrams. It is written in
\cite{BBD}, Proposition 1.1.11.
\end{proof}

We will also apply the octahedral axiom  (see  \S IV.1.1 of
\cite{gelman}) directly.  We recall that  it states the following: any diagram
$X\stackrel{f}{\to}Y\stackrel{g}{\to}Z$ can be completed to an
octahedral diagram. In particular, there exists a distinguished
triangle $\co(g\circ f)\to\co(g)\to \co(f)[1]  $, whereas  the morphism
$\co(g)\to \co(f)[1]$ is obtained by composing of two of the morphisms in
the distinguished triangles that define $\co(f)$ and $\co(g)$ (see
\S IV.1.8 of \cite{gelman}).

In particular, we will  need the following very easy application of
the octahedral axiom (we will apply it for the study of Postnikov
towers below).

\begin{lem}\label{lpost}
Let   $X$ be an object of $\cu$; consider a  (bounded above, below or both)
sequence of $\cu$-morphisms $\dots \to Y_{m-1}\to Y_{m}\to \dots$
equipped with morphisms $Y_i\to X$ such that all the corresponding
triangles commute.  Consider distinguished triangles
$Y_{ i} \to Y_{i+1}\to X_{i}$
for some $X_i\in \obj \cu$  (when the corresponding $X_l$ are
defined). Then for $Z_i=\co(Y_i\to X)$ there also exist
distinguished triangles
\begin{equation}\label{wdeck5} Z_i\to Z_{i+1}\to X_i[1]. \end{equation}
\end{lem}
\begin{proof}
The proof is immediate if one completes the commutative triangle
$Y_{i}\to Y_{i+1}\to X$ to an octahedral diagram.

\end{proof}

\subsection{Weight decompositions of morphisms;
 multiple weight decompositions of
objects}\label{wedecompm}

Starting from this moment the triangle
\begin{equation}\label{wdeck}
T_k[k]: X[k]\stackrel{a^k}{\to} X^{w\le k} \stackrel{f^k}{\to}
X^{w\ge k+1}\stackrel{b^k}{\to}X[k+1]\end{equation} will be
(an arbitrary choice of) a
weight decomposition of $X[k]$  for some $X\in\obj\cu$, $k\in\z$;
$T_k'[k]: X'[k]\stackrel{a'^k}{\to} X'{}^{w\le k}
\stackrel{f'^k}{\to} X'{}^{w\ge k+1}\stackrel{b'^k}{\to}X'[k+1]$
will be a weight decomposition of $X'[k]$. Sometimes we will drop
the index $k$ in the case $k=0$.

\begin{lem}\label{fwd}

1.  Let $l\le m$. Then for any morphism $g:X\to X'$ there exist
$h:X^{w\le m}[-m]\to X'{} ^{w\le l}[-l]$ and $i:X^{w\ge m+1}[-m]\to
X'{}^{w\ge l+1}[-l]$ completing $g$ to a morphism of triangles
$T_m\to T_l'$.

2. Let $l<m$. Then  $h$ and $i$ are unique.

3. For $l=m$, and  any two choices of $(h,i)$ and $(h',i')$ as above, we have
$h-h'=(s\circ f^m)[-m]$ and $i-i'=(f'^m\circ s')[-m]$ for some
$s,s'\in \cu(X^{w\ge m+1},X'{}^{w\le m})$.

\end{lem}
\begin{proof}
1,2: Immediate from Proposition \ref{lbbd}.

3. If suffices to consider the case $g,h,i=0$. Since $a^k[-k]\circ
h=0$, $T_k$ is a distinguished triangle, we obtain that $h'$ can
be presented as $(s\circ f^m)[-m]$. Dually, $i'$ can be presented
as $(f'^m\circ s')[-m]$.

\end{proof}

\begin{rema}\label{bmorwc}

1. For $l<m$ we will denote $i,h$ constructed by $g_{X^{w\le m},
X'{}^{w\le l}}$ and $g_{X^{w\ge m+1}, X'{}^{w\ge l+1}}$,
respectively.

For $l=m=0$ we will call any possible pair $(h,i)$ a {\it  weight
decomposition} of $g$.

2.  Suppose that $X$ admits a weight decomposition that {\it avoids
weight $0$} (this notion was introduced in Definition 1.6 of \cite{wild}), 
i.e.,  a weight decomposition such that $X^{w\le 0}\in \cu^{w\le
-1}$. Then such a decomposition is unique up to a unique
isomorphism. Indeed, in this case we can take $(X[-1])^{w\le
0}=X^{w\le 0}[-1]$. Therefore, we can apply part 2 of the previous
lemma for $l=-1$, $m=0$, $X'=X$, $g=\id_X$.

Moreover, the same method yields functoriality of such weight
decompositions; see
 Proposition 1.7. of loc. cit.

In fact, for these statements to be true it suffices to demand
$\hw(X^0,X^{-1})=\ns$. In particular, this argument yields that weight
decompositions of mixed motives are unique; cf. \S\ref{stfiltr} below and \cite{btrans}.

3. The statement of  Lemma \ref{fwd}(3) is the best possible
 in a certain sense.
It is not possible (in general) to choose $s=s'$. In particular, one
can take $$X=X'=\z/4\z\stackrel{\times 2}{\to} \z/4\z \in \obj
K^{[0,1]}((\z/4\z)-\operatorname{mod})\subset K(\ab) .$$ Then for
$g=0$ there
exists a pair $(h,i)= (\times 2,0)$ that is not homotopic to $0$.
Certainly, this example can be generalized to $X=X'=(R/r^2R\stackrel{\times r}{\to} R/r^2R)$ for
any commutative ring $R$, where $r\in R$ and $ r^2\nmid r$. In
particular, this problem is not a torsion phenomenon.

Note that the example of the weight decomposition described is
obviously not a 'nice' one. In particular, it cannot be extended
to a $3\times 3$ diagram. 
Yet adding this example to the obvious weight decomposition of
$\id_X$ one obtains another weight decomposition of $\id_X$ that
is not homotopy equivalent to the first one; yet it does not seem
to be 'bad' in any sense.

Still one can check that extension of morphisms ($X\to X'$ to
$(X^{w\le
i},X^{w\ge i+1})\to (X'{}^{w\le i},X'{}^{w\ge i+1})$) via part 1
of the Lemma
is
sufficient to prove the functoriality of the homology of the
weight complex of $X$ as defined in \S\ref{dwc} below (here the
homology objects belong to $\hw_*'$, see the Notation and 
Remark \ref{rkw}(2) below).

\end{rema}

We check that $g_{X^{w\le m}, X'{}^{w\le l}}$ and $g_{X^{w\ge m+1},
X'{}^{w\ge l+1}}$ are functorial in $g$ and the corresponding weight decompositions.

\begin{lem}\label{cwd}

Let $T''[j]: X''[j]\stackrel{a''^j}{\to} X''{}^{w\le j}
\stackrel{f''^j}{\to} X''{}^{w\ge
j+1}\stackrel{b''^j}{\to}X''[j+1]$ be a weight decomposition of
$X''[j]$ for $X''\in\obj\cu$ for some $j\le i\le 0$; let $p\in
\cu(X,X')$ and $q\in\cu(X',X'')$.

 If $j<0$ then for any choice of $(h',h'')$ satisfying
 \begin{equation}\label{wmor}
 h'\circ a^0 =a'^i[-i]\circ p\ \text{and } h''\circ a'^i[-i]
=a"^j[-j]\circ
 q\end{equation}
 we have $(q\circ p)_{X^{w\le 0},X''{}^{w\le j}}=h''\circ h'$, whereas
  for any choice of $(i',i'')$
 satisfying  $b'^i[-i]\circ i'=p [1]\circ b^0$ and
$b"^j[-j]\circ i''=q [1]\circ b'^i[-i]$ we have
 $(q\circ p)_{X^{w\ge 1},X''{}^{w\ge j+1}}=i''\circ i'$.

\end{lem}
\begin{proof}

We apply the uniqueness proved in the previous lemma.

 Both sides of the first equality calculate the only morphism $h$
that satisfies $h\circ a^0 =a''^j\circ (q\circ p)$, while both
sides of the second equality calculate the only morphism $i$ that
satisfies $b''^j\circ i=(q\circ p) [1]\circ b^0$.

\end{proof}

Now we prove that weight decompositions are 'exact' (in a certain sense).
We relate weight decompositions of $\co(X\to X')$ to those of
$X$ and $X'$;
this statement is easily seen to be related to the definition
of a
cone of a morphism
in $C(A)$.

\begin{lem}\label{twd}

Let $DT: X\stackrel{g}{\to} X'\to C[1]$ be a distinguished triangle.
Then $DT$ can be completed
to a diagram
\begin{equation}\label{dia3na3t}
\begin{CD}
X@>{a^i[-i]}>>w_{\le i}X @>{f^i[-i]}>> X^{w\ge i+1}[-i]\\
@VV{g}V@VV{g_{X^{w\le i},X'{}^{w\le i-1}} }V @VV{g_{X^{w\ge
i+1},X'{}^{w\ge i}} }V\\ X'@>{a'^{i-1}[1-i]}>>w_{\le i-1} X'
@>{{f'^{i-1}[1-i]}}>>X'{}^{w\ge i}[1-i]\\ @VV{}V@VV{}V @VV{}V\\
C[1]@>{}>>C_i[1-i]@>{}>>C'_{i}[1-i]\\
\end{CD}
\end{equation}
whose rows and columns are distinguished triangles, all squares
commute, $C_i,C_i'\in\obj \cu$. Moreover, the last row (shifted by
$[i-1]$) gives a weight decomposition of $C[i]$.

Besides, the choice of the   part of (\ref{dia3na3t}) consisting
of six upper objects and arrows connecting them is unique (even if
we don't demand that this part can be completed to the whole
(\ref{dia3na3t})).

\end{lem}
\begin{proof}
By  Lemma \ref{fwd}(2), $g$  can be uniquely completed to
a morphism of triangles that are the first two rows of (\ref
{dia3na3t}). Since the left upper square of (\ref {dia3na3t}) is
commutative, it can be completed to a $3\times 3$-diagram (see
Proposition \ref{3na3}). Hence the first two rows of this diagram
will be as in (\ref {dia3na3t}).  It remains to study the third
row.

 By  Proposition
\ref{bw}(3), the second column yields $C_i\in \cu^{w\le 0}$, whereas
the third column yields  $C_i'\in \cu^{w\ge 0}$. Hence $C[i]\to
C_i\to C_i'$ is a weight decomposition of $C[i]$.
\end{proof}

\begin{rema}\label{rtwd}

1. In fact, the lemma is valid in a more general situation. Suppose
that we have a pair of full subcategories $D,E\subset \cu$
that satisfy the orthogonality condition (for  weight structures, i.e.,
$E\perp D[1]$) and are
extension-stable (see Definition \ref{exstab}).
 Then any 'weight decompositions' of $X[i]$ and $X'[i-1]$
(defined similarly to the case when $D,E$ form a weight structure)
can be completed to a diagram (\ref{dia3na3t}) with $C_i\in \obj D$,
$C_i'\in\obj E$.

Indeed, it suffices to use the orthogonality to construct the
diagram required (for some $C_i,C_i'\in\obj \cu$). Next, the
second column yields $C_i\in \obj D$, whereas  the third
column yields $C_i'\in\obj E$.

We will use this statement below  for constructing weight
decompositions for certain 'candidate weight structures'.

2. Lemma \ref{twd} and its expansion described above show that
 it suffices to know
weight decompositions for some basic objects of $\cu$ in order to
obtain weight decompositions for all objects; see Theorems
\ref{recw} and \ref{madts}, and \S\ref{ddme} below. The situation is
quite different for $t$-structures; for this reason weight
structures are 'more likely to exist' (than $t$-structures),
especially in 'small' triangulated categories; cf.  Remark
\ref{exoticcond}(4).

\end{rema}

Now we study what happens if one combines  more than one weight
decomposition $T_k$.

\begin{pr}\label{mwdec}[Multiple weight decompositions]

1. [Double weight decomposition]

Let $T_k$ be (arbitrary and) fixed for some $X\in\obj \cu$
for  $k$ being equal to
some $i,j\in\z$, $i>j$.

Then there exist unique morphisms $s^{ij}: X^{w\le i}[j-i]\to
X^{w\le j}$, $q^{ij}:X^{w\ge i+1}[j-i]\to X^{w\ge j+1}$ making the
corresponding squares commutative. There also exists
$X^{[i,j]}\in \cu^{[0,i-j-1]}$, and distinguished triangles
\begin{equation}\label{wdeck1}
 X^{w\le i}[j-i]\stackrel{s^{ij}}\to X^{w\le j}\stackrel{c^{ij}}\to
X^{[i,j]}\stackrel{d^{ij}}\to X^{w\le i}[j-i+1]
\end{equation}
 and
\begin{equation}\label{wdeck2}
X^{[i,j]}[-1]\stackrel{x^{ij}}\to X^{w\ge
i+1}[j-i]\stackrel{q^{ij}}\to X^{w\ge j+1}\stackrel{y^{ij}}\to
X^{[i,j]} \end{equation}
 for some $\cu$-morphisms $c^{ij}, d^{ij},x^{ij},y^{ij}$.

2. [Infinite weight decomposition]

Let $T_k$ be (arbitrary and) fixed for all $k\in\z$.
Then for all $k\in\z$ there
exist unique morphisms $s^{k}: X^{w\le k}[-1]\to X^{w\le k-1}$,
$q^{k}:X^{w\ge k+1}[-1]\to X^{w\ge k}$ making the corresponding
squares commutative. There also exist $X^{k}\in \cu^{w=0}$, and
distinguished triangles
\begin{equation}\label{wdeck3}
 X^{w\le k}[-1]\stackrel{s^{k}}\to X^{w\le k-1}\stackrel{c^{k}}\to
X^{k}\stackrel{d^{k}}\to X^{w\le k}
\end{equation}
 and
\begin{equation}\label{wdeck4}
X^{k}[-1]\stackrel{x^{k}}\to X^{w\ge k+1}[-1]\stackrel{q^{k}}{\to}
X^{w\ge k}\stackrel{y^{k}}\to X^{k} \end{equation}
 for some $\cu$-morphisms $c^{k}, d^{k},x^{k},y^{k}$.

Moreover, $c^k$ and $x^k$ can be chosen equal to $y^k\circ
f^{k-1}$ and $(f^k\circ d^k)[-1]$, respectively.

\end{pr}
\begin{proof}

1. Applying Lemma \ref{fwd} for $X=X'$ and $g=\id_X$ we obtain the
existence and uniqueness of $s^{ij}, q^{ij}$. It remains to study
cones of these morphisms.

 The $3\times 3$-Lemma (i.e., Proposition \ref{3na3}) implies that
the morphism of distinguished triangles $T_i[i]\to T_j[j] $ can be completed to a
$3\times 3$ diagram whose rows and columns are distinguished
triangles. Hence there exists a distinguished triangle
$\co(\id_X)\to \co (s^{ij})\to \co(q^{ij})$; therefore $\co
(s^{ij})\cong \co({q^{ij}})$. Since $\cu^{w\le i-j-1}$ is
extension-stable (see  Proposition \ref{bw}(3))
 the distinguished triangle (\ref{wdeck1}) yields
$X^{[i,j]}\in\cu^{w\le i-j-1}$; the same argument applied to the
distinguished triangle (\ref{wdeck2}) yields
$X^{[i,j]}\in\cu^{w\ge 0}$.

2. The first part of the assertion is immediate from assertion 1
applied for $(i,j)=(k,k-1)$ for all $k\in\z$.

To prove the second part it suffices to complete the commutative
triangle $X[k]\stackrel{a^{k}}{\to} X^{w\le k}
\stackrel{s^{k}[1]}{\to} X^{w\le k-1}[1]$ to an octahedral
diagram.

\end{proof}

\begin{coro}\label{gcu}

  $\cu^{w=0}$ generates $\cu^b$ (as a triangulated category).

\end{coro}
\begin{proof}

 Since $\cu^b$ is a triangulated category that contains $\hw$, it
suffices to prove that any object of $\cu^b$ can be obtained
from objects of $\hw$ by a finite number of taking cones of
morphisms.

Let $X$ belong to $ \cu^{w\ge j}\cap \cu^{w\le i}$. Then we can take
$X^{w\le k}=0$ for $k<j$ and $X^{w\ge k}=0$ for $k>i$ in
(\ref{wdeck}). Then $X=w_{\le i}X$ and the formula (\ref{wdeck3})
gives a sequence of distinguished triangles implying that $X\in
\lan \cu^{w=0} \ra$.

\end{proof}

We will need the following definition several times.

\begin{defi}\label{dpoto}
We will denote by $Po(X)$ (a {\it weight Postnikov tower} for $X$)
all the 
 data of the (distinguished) triangles   (\ref{wdeck}), (\ref{wdeck3}),  and
(\ref{wdeck4}).\end{defi}

\begin{rema}\label{rpoto} 

1. For any $X,X'$, and arbitrary
weight Postnikov 
towers for them, any $g\in \cu(X,X')$ can be extended to a
morphism of Postnikov towers (i.e., there exist morphisms  $X^{w\le
k}\to X'{}^{w\le k},\ X^{w\ge
k}\to X'{}^{w\ge k},\ X^k\to X'{}^k$ for all $k$, such that  the corresponding
squares commute). Indeed,  Lemma \ref{fwd}  implies that we can
construct morphisms $T_k\to T_k'$ desired. Next we can
complete this data to morphisms of distinguished triangles of the type (\ref{wdeck3}) and
(\ref{wdeck4}). 

The main difficulty here is to prove that we can choose morphisms $X^k\to X'^k$ that will be compatible both with (\ref{wdeck3}) and
(\ref{wdeck4}). We will not really need this fact below (since it will always be sufficient for our purposes to consider either only (\ref{wdeck3}) or only
(\ref{wdeck4})). So, we will only sketch the proof of this statement. 

We should verify that a morphism of commutative triangles $(X[k]{\to} X^{w\le k}
{\to} X^{w\le k-1}[1])\to (X'[k]{\to} X'^{w\le k}
{\to} X'^{w\le k-1}[1])$ may be completed to a morphism of the corresponding octahedral diagrams (see the proof of  Proposition \ref{mwdec}(2)). This is certainly true if $\cu$ has a 'reasonable model' and our octahedra are compatible with this model; cf. the Corollary in \S1.7 of \cite{malt}. In the general case one can check the following: it suffices to verify that the morphisms of (\ref{wdeck3}) and of
(\ref{wdeck4}) yield the same morphisms on the homology of weight complexes (of $X$ and $X'$; the target of this homology is $\hw_*'$; see  Remark \ref{rkw}(2) below). The latter fact could be proved using the isomorphism (\ref{virtttr}) below; see  Theorem 2.3.1(II3) of \cite{bger}. 

We also note here that morphisms of weight Postnikov towers obtained this way are usually
not unique;
the choice of a 
 weight Postnikov tower (for $X$) is  not unique also.
 Still we will
prove below that 
weight Postnikov towers are 'unique and functorial up to
morphisms that are zero on homology'
(in a certain sense).

The corresponding functoriality facts for cellular Postnikov
 towers of
spectra are described in (Lemma 14 of) \S6.3 of \cite{marg}.

2.  Now, suppose that
we have an unbounded  system of morphisms $\dots Y_{m-1}\to Y_m\to
\dots$ equipped with compatible morphisms $Y_i\to X$ (as in Lemma
\ref{lpost}). Suppose also that $Y_l=0$ for $l\ll 0$, and that for
the corresponding $X_i$ we have $X_i\in \cu^{w=-i}$ for all
$i\in\z$.
 Then one can easily
check that  $Y_i\in \cu^{w\ge 1-i}$. Moreover, if we also have
$Y_l=X$ for $l\gg 0$, then $Y_i$  can be chosen for $w_{\ge
1-i}X$, whereas $Z_i=\co (Y_i\to X)$ can be taken for the
corresponding $w_{\le -i}X$.

Indeed, if  $Y_l=0$ for $l<N$ then $Y_l\in \cu^{w\ge 1-l}$ for
$l<N$. Now, for $l\ge N$ the distinguished triangles relating $Y_i$ with $X_j$
easily yield the same statement for all $l$ (by induction on $l$; 
here we apply  Proposition \ref{bw}(3)).

Moreover, the distinguished triangles $Y_i\to X\to Z_i$
 corresponding to 'weight decompositions' exist
by the definition of $Z_i$. Hence it suffices to check that $Y_l=X$
for $l>N$ implies $Z_i\in \cu^{w\le -i}$. We have $Z_i=0$ for $i>N$.
Hence our (last) claim can be easily deduced from the (distinguished) triangles
(\ref{wdeck5}).

3. More generally,  the 
argument 
above 
can be naturally extended to the case of extension-stable
'candidate weight
structures' (see Remark \ref{rtwd}).

 \end{rema}

\section{Weight filtrations and spectral sequences}\label{wspe}

The goal of this section is to study the {\it weight spectral
sequence} $T(H,X)\implies H(X)$ for $X\in\obj \cu$ and a
(co)homological functor $H:\cu\to A$.
$T$ is a  generalization of  Deligne weight spectral sequences
(see 
Remark \ref{rwfiltcoh}(2)), Atiyah-Hirzebruch ones (see \S\ref{sh}),
 and (essentially)
 motivic descent spectral sequences of \S7 of \cite{mymot} (cf.
Remark 7.4.4 of ibid. and \S\ref{twmot} below). 
So, 
if $H$ is a 'classical' realization of motives then
$T$ degenerates at $E_2$
(rationally) and its $E_2$-terms are exactly the graded pieces of
the weight
filtration. 

In \S\ref{swefil} we define the weight filtration for any functor
from $\cu$ to an abelian category.

In \S\ref{dwc} we define the {\it weight complex} of $X$ (in the
terms of   $Po(X)$). Its 'homology' computes the $E_2$-terms of
weight spectral sequences.

For simplicity we only consider  in detail only  spectral
sequences for homological functors (in \S\ref{swss}); dualization
immediately extends the result  to the cohomological functor case
(see \S\ref{swssc}). We also note that our weight spectral sequences induce
the standard weight filtration for the rational \'etale and Hodge
realizations of varieties (and motives); we also obtain a new
spectral sequence for motivic cohomology; see Remark
\ref{rwfiltcoh}.

Lastly, in \S\ref{strfun} we study the $D$-terms of the derived
exact couple, i.e., $D_2$ for a (co)homological weight spectral
sequence $T(F,X)$. There are two methods for constructing
these spectral sequences (see Remark \ref{tdss});
 the $D^{*,*}_2$-terms (for
a homological $F$) will be either $F'(X)=\imm (F(w_{\le
k+1}(X[l]))\to F(w_{\le k}(X[l])))$ or $F''(X)=\imm (F(w_{\ge
k+1}(X[l]))\to F(w_{\ge k}(X[l]))$. We
prove that $F',F''$ are both (co)homological. For $l=0$ the
sequence $F''(X)\to F(X)\to F'(X)$ extends to a long exact sequence
of functors. So, the properties of $F',F''$ are very similar to the
properties that would be satisfied by $t$-truncations of
$F$ if it is considered as an object
of some triangulated
'category of functors' $\du$ (so we call $F',F''$
 {\it virtual $t$-truncations} of $F$).
A certain explanation for this will be
given in \S\ref{sadst}  below; see also Remark
\ref{tmotco}(2) and (\S2.5 of) \cite{bger}. 
 Besides, we observe that $E^{*,*}_2(T(F,X))$ can be
obtained by composing two of our virtual $t$-truncations
('from different
sides').

\subsection{Weight filtration for (co)homological
functors}\label{swefil}

Let $A$ be an abelian
category.

\begin{defi}\label{dwfilf}

1. If $H:\cu\to A$  is any covariant functor then for any $i\in\z$
we define $W_i(H)(X)=\imm(H(w_{\ge i}(X))\to H(X))$.

2. If $H:\cu\to A$ is contravariant then we define
$W^i(H)(X)=\imm(H(w_{\le
i}X)\to H(X))$.

In both cases we will call the filtration obtained the
 {\it weight filtration} of $H(X)$.
\end{defi}

\begin{pr} \label{fwfil}

1. Let $H$ be covariant. Then the correspondence $X\to W_i(H)(X)$
gives a canonical subfunctor of $H(X)$. This means that
$W_i(H)(X)$ does not depend on the choice of a weight
decomposition of $X[i]$ and for any $f:X\to Y$ we have
$H(f)(W_i(H)(X))\subset W_i(H)(Y)$ (for $X,Y\in\obj\cu$).

2. The same is true for contravariant $H$ and $W^i(H)(X)$.
\end{pr}
\begin{proof}

1.  Part 1 of Lemma \ref{fwd} implies that for any choice of
weight decompositions of any $X[i],Y[i]$ we have
$H(f)(W_i(H)(X))\subset W_i(H)(Y)$. In particular, taking $Y=X$,
$f=\id_X$ we obtain that $W_i(H)(X)$ does not depend on the choice
of the weight decomposition of $X[i]$

2. The statement is exactly the dual of assertion 1 (see Remark
\ref{rdd}).

\end{proof}

\begin{rema}\label{semim}

1. In the case when $H$ is  (co)homological, we can replace the
images in Definition \ref{dwfilf} by certain kernels (since they
coincide).

Moreover,  Proposition \ref{fwfil} is also true for $A$ being {\bf any}
category with well-defined images of morphisms.

2. A partial case of this method for defining weight filtration
for cohomology  was (essentially) considered in Proposition 3.5 of
\cite{ha3}.

3. 
Recall now that we have a natural embedding $i:\cu\to \cu_*$, that sends
$X\in \obj \cu$ to $X_*=\cu(-,X)$. Then for any $j\in\z$ and
$X\in\obj\cu$ we can define a functor $W_i(X)$ that sends $Y\to
W_j(Y^*(X)):\cu\to\ab$. This yields an object of $\cu_*$; it equals
$W_j(i)(X)$.  We obtain a sequence of functors $W_j:\cu\to \cu_*$.
The usual Yoneda's isomorphism $F(X)\cong
\mo_{\adfu(\cu^{op},\ab)}(X_*,F)$ for a contravariant functor $F:\cu\to
\ab$ can be easily generalized to
\begin{equation}\label{covar} W^j(F)(X)\cong
\mo_{\adfu(\cu^{op},\ab)}((X_*/W_j(X)),F)\end{equation}
 Hence the sequence of functors
$W_j$ yields a description of weight filtrations for all
contravariant functors $\cu\to \ab$. Moreover, for any %
homological
$F:\cu\to
\ab$ 
one can similarly define an isomorphism
\begin{equation}\label{contrvar} W_j(F)(X)\cong
\mo_{\adfu(\cu,\ab)}(W^i(X),F)\end{equation}

In particular, one can apply this construction to the Chow weight
filtration of Voevodsky's motives (see \S\ref{remmymot} below). For
any motif $X$ (an so, for any variety) one obtains a sequence of
objects of $\dmge{}_*$ which could be called {\it semi-motives}.
These objects contain important  information on $X$. In particular,
(\ref{covar}) and (\ref{contrvar}) show that they have
both homological and cohomological realizations!

\end{rema}

\subsection{Weight complexes of objects (definition)}\label{dwc}

Now we describe the weight complex of $X\in\obj \cu$. We will
prove that it is canonical and functorial (in a certain sense) in
\S\ref{sfwc} below.

We adopt the notation of subsection \ref{wedecompm}.

\begin{defi}\label{dwcompl}
We define  the morphisms $h^i:X^i\to X^{i+1}$ as $c^{i+1}\circ
d^i$. We will call $t(X)=(X^i, h^i)$ the {\it weight complex} of $X$.

\end{defi}

Note  that all information on $t(X)$ is contained in $Po(X)$
(including the relation of $t(X)$ to $X$).

\begin{pr}\label{fprwc}

1. Weight complex (of any $X\in\obj\cu$) is a complex
 indeed, i.e., we have $d^2=0$.

2. If for some $i\in\z$ and $X\in \cu^{w\le i}$ (resp. $X\in
\cu^{w\ge i}$) then there exists a choice of the weight complex of
$X$ belonging to $C(\hw)^{\le i}$ (resp. to $C(\hw)^{\ge i}$).

\end{pr}
\begin{proof}

1. We have
$$h^{i+1}\circ h^i=c^{i+2}\circ (d^{i+1}\circ c^{i+1})\circ d^i
=c^{i+2}\circ 0\circ d^i=0$$ for all $i$.

2. Similarly to the proof of Corollary \ref{cub}, we can take
$X^{w\ge k}=0$ for $k>i$ (resp. $X^{w\le k}=0$ for $k<i$). Then we
would have $w_{\le k}X=X$ for $k\ge i$ (resp. $w_{\ge k}X=X$ for
$k\le i$). Therefore the corresponding choice of the weight
complex of $X$ belongs to $C(\hw)^{\le i}$ (resp. to $C(\hw)^{\ge
i}$) by definition.

\end{proof}


\subsection{Weight spectral sequences for
homological functors}\label{swss}

Let $A$ be an abelian category; let $H:\cu\to A$ be a homological
functor (i.e., a covariant additive functor that transfers
distinguished triangles into long exact sequences). The
cohomological functor case will be obtained from the homological
one by dualization.

Let $X$ be an object of $\cu$, $(X^i,h^i)=t(X)$. We construct
a spectral sequence 
whose $E_1$-terms are $H(X^i[j])$, which converges to $H(X[i+j])$
in many important cases.

\begin{defi}\label{numbers}
 We denote $H(Y[p])$ by $H_p(Y)$ for any $Y\in\obj\cu$.

 For a cohomological $H$ we will denote   by $H^p(-)$
the functor $Y\mapsto H(Y[-p])$. 
\end{defi}

First we describe the exact couple. It is obtained by applying $H$
to the data contained in a 
weight Postnikov tower for $X$
(see Definition \ref{dpoto}). In the
first three parts of Theorem \ref{hsps} we will fix the choice 
 of this tower.

Our exact couple is almost the same as the  couple in \S IV2,
Exercise 2, of \cite{gelman}. We take $E_1^{pq}=H_q(X^p)$,
$D_1^{pq}=H_q(X^{w\ge p})$.
 Then the
distinguished triangles (\ref{wdeck4}) endow $(E_1,D_1)$ with the
structure of an exact couple.


\begin{theo} \label{hsps}[The homological weight spectral sequence]

I There exists a spectral sequence $T=T(H,X)$  coming from our $(E_1,D_1)$ with $E_1^{pq}=
H_q(X^p)$  such that the map
$E_1^{pq}\to E_1^{p+1q}$ equals $H_{q}(h^p)$.

II $T(H,X) $ converges to $ H_{p+q}(X)$ in either of the following cases:

(i) $X\in \cu^b$.

(ii) $H$ vanishes on $\cu^{w\ge q}$ for $q$ large enough and on
$\cu^{w\le q}$ for $q$ small enough.

(iii) $X\in \cu^-$ (resp $\cu^+$) and $H$ vanishes on $\cu^{w\le
q}$ for $q$ small enough (resp. on $\cu^{w\ge q}$ for $q$ large
enough).

In all these cases  the corresponding
filtration on $H_*(X)$ coincides with the (weight)
 filtration described in  Definition
\ref{dwfilf}.

III  $T$ is  functorial with respect to $H$, i.e., for any
transformation of functors $H\to H'$ we have a canonical morphism
of spectral sequences $T(H,X)\to T(H',X)$; these morphisms respect
sums and compositions of transformations.

 IV $T$ is canonical and functorial with respect to  $X$
starting from $E_2$.

\end{theo}
\begin{proof}

I  These are just the standard properties of a spectral sequence coming from a Postnikov tower; see the Exercises after \S IV.2 of \cite{gelman}.

II In case (ii) $E_1(T)$ is obviously bounded.

In case (i) this will be also true if we choose the weight complex
of $X$ to be bounded (we can do this by the definition of
$\cu^b$). Now, for an arbitrary choice of the weight complex
 we  also
obtain that $T$ will be bounded starting from
$E_2$  by assertion IV. 

The proof of boundedness in case (iii) is similar.

The connecting maps $w_{\ge p} X\to X$
 also yield the connection desired of $E^{pq}_\infty$
with $H_{p+q}(X)$ 
(since they are compatible with (\ref{wdeck4})).
Moreover, the induced filtration $F_*$ on $H_*(X)$ is the
weight filtration
          (of Definition \ref{dwfilf}) indeed, since
 \begin{equation} F_pH_{p+q}(X)=\imm (D_1^{pq}\to H_{p+q}(X))=
W_p H_{p+q}(X).\end{equation}

III This is obvious since all the components of the exact couple are
functorial with respect to $H$.

IV It suffices to check  that the correspondence  sending $X $ to
the derived exact couple (i.e., $(E_2,D_2)$ + the connecting
morphisms) defines a functor. 
Now, since any
$g\in\cu(X,X')$ can be extended to a morphism
 $Po(X)\to Po(X')$ (see  Remark \ref{rpoto}(1); here one can take any possible
$Po(X), Po(X')$ and does not have to consider the triangles (\ref{wdeck3}); the extension is not unique), we obtain that $g$
is compatible with at least one morphism of the original
couples (i.e., of $C_1=(D_1\to D_1\to E_1\to D_1)$).
It remains to prove that the induced morphism of the
derived couples $C_2(X)\to C_2(X')$ coming from this
construction is uniquely determined by $g$; 
hence it suffices to prove that the correspondences
$X\mapsto D_2^{pq}(T)$ and $X\mapsto E_2^{pq}(T)$ define
(canonical) functors (so we don't have to mind the connecting
morphisms of $C_2$). 

Now, $D_2^{pq}(T)=\imm (H_{q-1}(X^{w\ge p+1})\to
H_{q}(X^{w\ge p}))$ is canonical and functorial by 
Proposition \ref{trfun}(I) below. $E_2$ is also functorial
since it can be factored through the weight complex
$t$ (whose functoriality is checked \S\ref{sfwc}
 below);  see   Remark \ref{rkw}(3) and  also 
Remark \ref{rtrfun}(4).

\end{proof}

\begin{rema}\label{tdss}

One can easily 'dualize' the exact couple above (see Remark \ref{rdd}). This means the following:
 there exists an exact couple that yields the same spectral sequence
(so all $E_*^{**}$ do not change) but with $D_1^{pq}=H_q(X^{w\le p-1})$.

Certainly, the same observation can be applied to cohomological
weight spectral sequences; see below.

\end{rema}

\subsection{Weight spectral sequences for
cohomological functors; examples}\label{swssc}

Inverting the arrows in $\cu$,  we
obtain the following cohomological analogue of the previous theorem.

\begin{rema}\label{rtdss}
If we dualize the exact couple from Theorem
\ref{hsps} directly, then we will obtain
 $D_1^{pq}=
H^{q}(X^{w\le -p})$.
For an 'alternative' exact couple (see Remark \ref{tdss})
we have $D_1^{pq}=
H^{q}(X^{w\ge 1-p})$. Yet this does not affect the spectral
 sequence.\end{rema}

\begin{theo} \label{csps} [The cohomological weight spectral sequence]

I There exists a spectral sequence $T=T(H,X)$ with $E_1^{pq}=
H^{q}(X^{-p})$   such that
the map $E_1^{pq}\to E_1^{p+1q}$ equals $H^{q}(h^{-1-p})$. 

II $T(H,X)$ converges to $H^{p+q}(X)$ in either of the following cases:

(i) $X\in \cu^b$.

(ii) $H$ vanishes on $\cu^{w\ge q}$ for $q$ large enough and on
$\cu^{w\le q}$ for $q$ small enough.

(iii) $X\in \cu^-$ (resp $\cu^+$) and $H$ vanishes on $\cu^{w\le
q}$ for $q$ small enough (resp. on $\cu^{w\ge q}$ for $q$ large
enough).

The
corresponding filtration on $H^*(X)$ coincides with the weight filtration
 of
Definition \ref{dwfilf}.

III  $T$ is  functorial with respect to $H$, i.e., for any
transformation of functors $H\to H'$ we have a morphism of
spectral sequences $T(H,X)\to T(H',X)$.

 IV $T$ is canonical and (contravariantly) functorial
with respect to  $X$
starting from $E_2$.

\end{theo}

\begin{proof}
It suffices to apply Theorem \ref{hsps} to the (homological) functor
$H':\cu^{op}\to A$.
\end{proof}

\begin{rema}\label{rwfiltcoh} \label{rmotco}[Examples:
'classical' realizations and motivic cohomology]

1. Suppose that $w$ be bounded and there are no maps between
distinct weights for $H$, i.e.,  there exists a family of full
abelian subcategories $A_i\subset A$ that contain all $A$-subquotients  of their objects,
 such that $H^i(P)\in \obj A_i$
for all $i\in\z,\ P\in \cu^{w=0}$, and there are no non-zero
$A$-morphisms between distinct $A_i$. Then we easily obtain that
$T(H,X)$
 degenerates at $E_2$. Besides,  for any $a\in \obj A$
there cannot exist more than one
finite filtration $V^j$ on $a$ such that $V^j(a)/V^{j-1}(a)\in \obj
A_j$ for all $j$; whereas the Theorem shows that for any $X\in \obj \cu$ 
our filtration   $W^i(H)(-)$ is such a filtration (for $H(X)$).

2. We will see in \S\ref{dgmot} below that Voevodsky's
$\dmgm(\supset\dmge)$ admits a {\it Chow weight structure} whose heart
is $\chow(\supset\chowe)$. Hence we obtain certain weight spectral sequences (that we will call {\it Chow-weight} ones) 
and weight
filtrations for all realizations of motives. In particular, we have them
for the \'etale and singular realizations of motives, and for
motivic cohomology.

Now, it is well known that for the rational \'etale and singular 
realization of motives there are no non-zero morphisms between  distinct weights
(in the corresponding categories
 of mixed structures) in the sense described above.
Therefore for the rational \'etale and singular realization of motives
our weight filtration coincides with the usual one (up to a shift of
indices; see also \S7.4 of \cite{mymot}).

 Recall also that
'classically' the weight filtration is well-defined only for
cohomology with rational coefficients. Yet
our method allows us to
define canonical weight filtrations integrally; this generalizes the
construction of Theorem 3 of \cite{gs}.

Now we consider the case of   motivic cohomology. A simple example
of the spectral sequence obtained comes from 
Bloch's long exact
localization sequence for higher Chow groups of varieties. 
Recall that it relates the motivic cohomology of
$X\setminus Z$ to  that of $X$ and $Z$, where $Z,X$ are smooth,
$Z$ is closed in $X$.
In the motivic setting it comes from the Gysin distinguished triangle
 $$\mg(X\setminus Z)\to \mg(X)\to\mg(Z)(c)[2c]$$ ($c$ is the
codimension of $Z$; see \S\ref{remmymot}
below
and Proposition
    5.21 of \cite{deg2}).
Now, suppose additionally that $X$ is  projective (so $Z$ also is).
Then  $\mg(X)$
and $\mg(Z)(c)[2c]$ are Chow motives; so they belong to
$\dmge^{w_{Chow}=0}$ (see \S\ref{hchow}). Since motivic cohomology is a
(representable) cohomological functor on $\dmge$, we obtain the following:
the Chow-weight spectral sequence converging to
the motivic cohomology of $X\setminus Z$
(corresponding to $w_{Chow}$) reduces to the corresponding Bloch's long exact
                                              sequence.
Since the latter is non-trivial in general, the Chow-weight
spectral sequences obtained is non-trivial either; 
 it appears  not to be mentioned in the literature.
This filtration is compatible with  regulator maps (whose
targets are classical cohomology theories). Unfortunately,
 morphisms of motivic cohomology of motives (induced by
$\dmge$-morphisms) 
are not necessarily strictly compatible with the weight filtration for this theory
(in contrast
with the properties of the weight filtration for rational singular and
\'etale cohomology).

\end{rema}

\subsection{ Higher $D$-terms of exact couples;
virtual $t$-truncations for (co)homological functors}\label{strfun}

Now we study the higher $D$-terms of exact couples
 (i.e., the $D$-terms for derived exact couples)
and especially $D_2$ for  (co)homological weight
spectral sequences, in more detail.

Let $A$ be an abelian category;
let $j>0$, $k\in\z$, be fixed.

\begin{pr}\label{trfun}

I Let $F:\cu\to A$ be a covariant functor. Then the assignments
$F_1=F_1^{kj}:X\mapsto \imm (F(w_{\le k+j}X)\to F(w_{\le k}X))$ and
 $F_2=F_2^{kj}:X\mapsto
\imm (F(w_{\ge k+j}X)\to F(w_{\ge k}X))$ define functors $\cu\to
A$ that do not depend (up to a canonical isomorphism) from the
choice of weight decompositions. 
Besides, there exist natural transformations  $F_2\to F\to F_1$.

II Let $F:\cu\to A$  be a homological (covariant) functor; let $j=1$. Then the
following statements are valid.

1. $F_l$ ($l=1,2$) are also homological.

2. The natural transformations $F_2\to F\to F_1$ extend canonically to a
  complex of functors
$\dots\to F_1\circ [-1]\to F_2\to F\to F_1\to F_2\circ [1]\to \dots$ 
  that is long exact when applied to any $X\in\obj\cu$.

III. Let $F:\cu\to A$ be a contravariant functor.

1. The assignments
$F_1=F^{kj}_1:X\mapsto \imm (F(w_{\le k}X)\to F(w_{\le
k+j}X))$ and $F_2=F_2^{kj}:X\mapsto \imm (F(w_{\ge k}X)\to F(w_{\ge
k+j}X))$  define contravariant functors $\cu\to\ab$ that do not
depend (up to a canonical isomorphism) from the choice of weight
decompositions. 
There exist natural transformations $F_1\to F\to F_2$.

2. If $F$ is cohomological, $j=1$, then $F_l$ ($l=1,2$) also are;
the transformations $F_1\to F\to F_2$ extend canonically to a long
exact sequence of functors
$\dots\to F_2\circ [1]\to F_1\to F\to F_2\to F_1\circ [-1]\to \dots$ 
(i.e., the sequence is exact when applied
 to any $X\in\obj\cu$).

\end{pr}
\begin{proof} 
I We use the following very simple observation: for any commutative square in $A$
$$\begin{CD}
X@>{f}>>Y \\
@VV{h}V@VV{}V \\
Z@>{g}>>T
\end{CD}$$
if we fix the rows then the morphism $g\circ h:X\to T$ completely determines the morphism $\imm f\to \imm g$ induced by $h$.

Hence Lemma \ref{cwd} easily implies that both $F_l$ ($l=1,2$) are
  well-defined and
functorial. Indeed, by this lemma the
morphisms $w_{\le k+j}X\to w_{\le k}X'$ and
 $w_{\ge k+j}X\to w_{\ge k}X'$ 'induced' by 
$g\in \cu(X,X')$ are compatible with
arbitrary morphisms of weight decompositions of $X$ and $X'$
 that come from morphisms of objects.

It follows that fixing  arbitrary weight decompositions
for all $X\in\obj\cu$ one obtains $F_1(X)$, $F_2(X)$,
and functorial connecting
morphisms $F_1(X)\to F_1(X')$, $F_2(X)\to F_2(X')$ for $g\in \cu(X,X')$,
that depend  on the choices made only up to a canonical isomorphism.

It remains to note that the connecting maps
$w_{\le k+j}X\to w_{\le k}X'$ and
 $w_{\ge k+j}X\to w_{\ge k}X'$ were chosen to be compatible with
  $w_{\ge k}X\to X\to w_{\le k+j}X$; this yields the existence of
the transformations in question.

II 1. Since $F$ is additive, both $F_l$ also are.

Now we check that $F_1$ is homological. It suffices to check that
for any distinguished triangle $U=(C\to X\to X')$ the sequence
$F_1(X'[-1])\to F_1(C)\to F_1(X)$ is half-exact (i.e., exact in the
middle term). Note that this sequence
 is obviously a complex (since the composition of morphisms
is zero, $F_1$ of them also is).

Note also that we can use arbitrary choices of weight
decompositions (by assertion I). We fix arbitrary choices of weight decompositions for (shifted)
$X,X'$  and look for a 'nice' weight decomposition
of $C$.

We use the notation of \S\ref{wedecompm}. In particular, we easily
obtain that $F_1(X)=\ke (F(w_{\le k}X)\to F(X^{k+1}[-k]))$ and
$F_1(X'[-1])=\ke (F(w_{\le k-1}(X')[-1])\to F(X'{}^{k}[-k]))$. Our
goal is to find a weight decomposition of $C$ such that there will
be a distinguished triangle $w_{\le k-1}(X')[-1]\to w_{\le k}C\to
w_{\le k}X$ compatible with $U$ and $F_1(C)=\ke (F(w_{\le k}C)\to
F((X^{k+1}\oplus X'^{k})[-k]))$.

We apply Lemma \ref{twd} for
$i=k,k+1$. By the lemma, the triangles $C[i]\to C_i\to C_i'$
obtained from (\ref{dia3na3t}) by shifting the last row are weight
decompositions of $C[i]$ (for all $i\in\z$).

Next, we apply Lemma \ref{twd} to the morphism $g_{ X^{w\le k+1},
X'{}^{w\le k}}[k+1]$ and the weight decompositions
$X^{k+1}{\to} X^{w\le k+1}\to X^{w\le k}[1]$ and
$X'{}^{k}[1]{\to} X'{}^{w\le k}[1]\to
X'{}^{w\le k-1}[2]$ of the corresponding objects. We obtain a
diagram
\begin{equation}\label{dia3na30}
\begin{CD}
X^{k+1}@>{}>>X^{w\le k+1} @>{}>> X^{w\le k}[1]\\
@VV{}V@VV{}V @VV{}V\\
X'{}^{k}[1]@>{}>>X'{}^{w\le k}[1] @>{}>>X'{}^{w\le
k-1}[2]\\ @VV{}V@VV{}V @VV{}V\\
D_{k+1}[1]@>{}>>C_{k+1}[1]@>{t[2]}>>C_{k}[2]\\
\end{CD}
\end{equation}
for some $D_{k+1}\in\obj \cu$ and some $t$. The first column gives
$D_{k+1}\cong X^{k+1}\oplus X'^{k}$. Indeed, $X^{k+1},X'^{k}\in
\cu^{w=0}$, whereas $\cu^{w=0}$ is extension-stable by 
Proposition \ref{bw}(3), and any extension in $\hw$ splits by part 7 of
loc. cit.  Besides,  $t$ equals the
 morphism connecting the weight
 decompositions of $C[i]$ for $i=k,k+1$, since this morphism is unique
 by Lemma \ref{fwd}(2). So, we have achieved our goal (as declared above).

 Now, we have a half-exact sequence
 $F(w_{\le k-1}(X')[-1])\to F(w_{\le k}C)\to F(w_{\le k}X)$.
 Therefore, if $x\in F_1(C)$ vanishes in $F_1(X)$, we obtain that
 it comes from some
  $y\in F(w_{\le k-1}(X')[-1])$.
 Hence it suffices to check that the image of $y$ vanishes in
$F(X'{}^{k}[-k])$.
 Since the image of $x$ in
$F((X'{}^{k}\oplus X^{k+1})[-k])$ vanishes, it suffices to note that
 $F(X'{}^{k}[-k])\hookrightarrow F((X'{}^{k}\oplus X^{k+1})[-k])$.
  Certainly, this reasoning could have been written down without
the use of
'elements' (of objects of $A$).

$F_2$ is homological for similar reasons; this fact also immediately
 follows
 from the previous one by  Remark \ref{rdd}(2).

2.  Let weight decompositions of all $X[i]$ be fixed.

The exactness of
$$ \begin{aligned} F_2(X)=\imm (F(w_{\ge k+1}(X))\to
F(w_{\ge k}X)) \to F(X)\\
\to F_1(X)= \imm (F(w_{\le k+1}X)\to F(w_{\le k}X))\end{aligned} $$
in $F(X)$ is immediate from the exactness of $F(w_{\ge k+1}X)\to
F(X)\to F(w_{\le k}X)$ (in the middle).

Next, by   Proposition \ref{fwd}(2) we obtain that $\id_X$
yields canonically a diagram
\begin{equation}\label{dia3na4}
\begin{CD}
F(w_{\ge k+2}X) @>{}>> F(X)@ >{}>>F(w_{\le k+1}X)
@>{}>>F((w_{\ge k+2}X)[1])\\
@VV{}V@VV{F(\id_{X})}V @VV{}V @VV{}V\\ F(w_{\ge k+1}X) @>{}>> F(X)@ 
>{}>>F(w_{\le k}X) @>{}>>F((w_{\ge k+1}X)[1])
\\
@VV{}V@VV{}V @VV{}V @VV{}V\\ F(X^{k+1}[-1-k]) @>{}>> 0 @>{}>>
F(X^{k+1}[-k]) @>{}>> F(X^{k+1}[-k])
\end{CD}
\end{equation}
the rows and columns are exact in  non-edge terms.

Hence we obtain a well-defined functorial morphism
$$\begin{gathered} F_1(X) 
 \to \imm(F((w_{\ge k+2}X)[1])\to F((w_{\ge k+1}X)[1]))
=F_2(X[1]).\end{gathered}$$
This map will be our boundary morphism (of functors).

Now, to check the exactness of the complex $F(X)\to F_1(X)\to F_2(X[1])$
 in the middle it suffices to check the following inclusion:
$$\begin{gathered}\cok (F(w_{\le k+1}X)\to F(w_{\le k}X))=
\ke (F(w_{\le k}X)\to F(X^{k+1}[-k])) \\ \hookrightarrow \cok
(F(w_{\ge k+2}X[1])\to F(w_{\ge k+1}X[1]))= \ke (F(w_{\ge
k+1}X[1])\to F(X^{k+1}[-k]));\end{gathered}$$ cf. the proof of
assertion II1 (above). Hence the last row of (\ref{dia3na4}) yields
the result.

Lastly, the exactness of $F_1(X[-1])\to F_2(X)\to F(X)$ in the
middle follows easily from the previous exactness by 
Remark \ref{rdd}(2).

III These statements are exactly the duals  of assertions I,II;
 see (part 1 of) Remark \ref{rdd}.
\end{proof}

\begin{rema}\label{rtrfun}[On virtual $t$-truncations of $F$]

 1. For a (co)homological $F$ we will call $F^{k1}_l$, $l=1,2,\ k\in\z$,
{\it virtual $t$-truncations} of $F$.
Note that $F$ often lies in a certain triangulated 'category of
(co)homological functors' $\du$ (of functors $\cu\to A$). The
simplest case is when we consider  $\du\supset \cu$ (or $=\cu$) and
$F$ is the restriction to $\cu$ of the functor (co)represented by
some $Y\in\obj\du$; see also   Remark \ref{tmotco}(2). If such a
$\du$ exists, then
  the virtual $t$-truncations defined
are often actual $t$-truncations of $F$ (corresponding to  a certain
 $t$-structure $t$ on $\du$; see  Theorem
\ref{sdt}(7,8), and
  Remark \ref{rsdt}(5) for the general case).  In the case when $\cu=\du$ this
 $t$-structure is called  
 {\it adjacent} to
 $w$; see Definition \ref{deadj}.

Still, it is very amusing that these $t$-truncated functors as well
as their transformations corresponding to $t$-decompositions (see
Definition \ref{dtstr}) can be defined without specifying any
particular $\du$!

2. In particular, for any $Y\in \obj SH$ the functors
(co)represented by $t$-truncations of $Y$ with respect to the
Postnikov $t$-structure are exactly our $F^{k1}_l$ (for $l=1,2,\
k\in\z$; here $F$ is either $SH(Y,-)$ or $SH(-,Y)$; see \S\ref{sh}
for the definition of $t_{Post}$). Hence one can express the 
restrictions of $F^{k1}_l$ to $SH_{fin}$ in terms of $F$ restricted
to $SH_{fin}$ (i.e., without considering infinite spectra). Note that
one can obtain an Eilenberg-Maclane spectrum for $\z$ by considering
the zeroth Postnikov '$t$-cohomology' of $S^0$.

 A similar observation
can be applied to the pair $\dmge\subset \dme$ (see \S\ref{ddme}
below). So, though the {\it Chow $t$-truncations}
 of geometric motives are
(usually) non-geometric, we can express functors represented by them
in terms of morphisms of $\dmge$.

3.  $F^{k1}_2(X[p])$ defined in part I of the proposition,
yield the $D$-terms of the
derived couple for the exact couple used for the proof
of Theorem \ref{hsps} (for a homological $F$).
By Remark \ref{rtdss}, for a cohomological $F$
the corresponding
$D$-terms 
are given by
$F^{k1}_1(X[p])$ (defined in part III of the proposition).

Recall also that making the alternative choice of exact couples interchanges
the roles of $F^{k1}_l$ ($l=1,2$) here, without changing the spectral
sequences (see Remark \ref{tdss}). Moreover, for any $j\ge 1$ the
functors $F^{jk}_l(X[p])$
($l=1,2$) yield $D_{j+1}(T)$ (i.e., they calculate the $D$-terms of
both possible choices of higher derived exact couples).

4. The definition  implies that any two of our virtual
$t$-truncation operations (for different $k,l$) commute. Besides for
a homological $F$ one can easily check that
\begin{equation}\label{virtttr} (F^{1,0}_1)^{1,-1}_2
=\imm (F(X^{0}\to
X^1)\to F(X^{-1}\to X^0))\cong E^{0,0}_2(T(F,X));\end{equation} a
similar fact is also valid in the case when $F$ is cohomological.
 Thus the terms $E^{*,*}_2(T(F,X))$ of the weight spectral
sequence can be
expressed in terms of our virtual $t$-truncations; in particular,
 they are
given by well-defined (co)homological functors. 
We prove that this approach yields a full description of the derived exact
couple for  $T$ 
in (part II of) Theorem 2.4.2 of \cite{bger}.

 5. 
The remarks above can be vastly extended (see \S\S2.3--2.6 of \cite{bger}).
Let   $\Phi:\cu^{op}\times \du\to A$ be a {\it (nice) duality}
(see   Remark \ref{rsdt}(5) below and  \S2.5 of \cite{bger});
Suppose that a $t$-structure $t$ for $\du$ is {\it orthogonal} to $w$
(for $\cu$) with respect to $\Phi$. In this case the   virtual
 $t$-truncations 
  of functors
 of the type $\Phi(-,Y),\ Y\in \obj\du$ are exactly the
 functors 'represented via $\Phi$' by the actual
 $t$-truncations of $Y$ (corresponding to $t$; see Proposition
2.5.4 of ibid.). This allows us to establish a natural isomorphism
for the weight spectral sequence for $\Phi(-,Y)$ to the one coming
 from $t$-truncations of $Y$ (in Theorem 2.6.1 of ibid.).

 Moreover, one can give a certain 'axiomatic' description of
virtual $t$-truncations; see Theorems 2.3.1 and 2.3.5 of ibid.

6. In general, for $X\in \obj \cu$ and $i\in \z$ there  does not exist an $X'$ such that for any cohomological $H$ we have $\tau_{\le i}H(X)\cong H(X')$. Yet if there exists a choice of $w_{\le i+1}X$ and $w_{\le i}X$ such that $w_{\le i}X$ is a retract of $w_{\le i+1}X$, then one can take $w_{\le i}X$ for such an $X'$. In particular, this is the case if $X$ {\it avoids weight $i+1$}; see Remark \ref{bmorwc}(2).

\end{rema}

\section{The weight complex functor}\label{swcomp}

In \S5 of \cite{mymot} for a triangulated category $\cu$ with a {\it
negative differential graded enhancement}
(this includes Voevodsky's motives)
 an exact conservative {\it
weight complex} functor $t_0:\cu\to K(\hw)$ (in our notation) was
constructed.  The goal of this section is to
extend this result to the case of arbitrary $(\cu,w)$. A reader only
interested in motives could skip this section (since a 'stronger'
version of the weight complex functor for all 'enhanceable'
categories will be constructed in  \S\ref{cwecom} below).

 As shown in Remark \ref{rpoto}, any $g\in\cu(X, X')$ (for $X,X'\in\obj\cu$) can be
extended to a morphism of (any possible) 
weight Postnikov towers for $X,X'$. Moreover, for compositions of $g$'s the
corresponding morphisms of weight Postnikov towers can be composed. Yet,
as the example of   Remark \ref{bmorwc}(3) shows, this
construction cannot give a canonical morphism of weight complexes in
$K(\hw)$. We have to consider a certain factor $K_\w(\hw)$ of this
category. This factor is no longer triangulated (in the general
case; yet cf. Remark \ref{rctst}). Still the kernel of the
projection $K(\hw)\to K_\w(\hw)$ is an ideal (of morphisms) whose
square is zero; so our (weak) weight complex functor is not much
worse than the 'strong' one of \cite{mymot}.

We define and study $K_\w(\hw)$ in \S\ref{wcatc}. We construct the
weight complex functor $t$ in \S\ref{sfwc} and prove its main
properties in \S\ref{smpwc}. One of our main tools is the {\it
weight decomposition functor} $WD:\cu\to K_\w^{[0,1]}(\cu)$; see
Theorem \ref{fwc}.

One of the main properties of the functor $t$ is that it
calculates the $E_2$-terms of the weight spectral sequence $T$,
see  Remark \ref{rkw}(3). In fact, this is why $t$ it called
the weight complex; this term was used for the first time in
\cite{gs} (see \S2 and \S3.1 of ibid.).

\subsection{The weak category of complexes}\label{wcatc}

Let $A$ be an additive category. We will need the following, very
natural definition. The author would like to note that this definition,
 as well as several related (and interesting) definitions
 and results were also independently introduced in \cite{pir}. 

\begin{defi}\label{didmo}
 A class $T$ of morphisms in $A$
will be called a (two-sided) ideal if it is closed with respect to
sums and differences (of two morphisms of $T$ lying in the same
morphism group), finite directs sums, and compositions with any
morphisms of
$A$.

We will abbreviate these properties as $T\nde \mo A$.
\end{defi}

\begin{rema}\label{ridmo}
For any $T\nde \mo A$ we can consider an additive
category $A/T$ whose object are the same as for $A$, and
$A/T(X,Y)=A(X,Y)/T(X,Y)$ for all $X,Y\in \obj A$.

Besides, it is easily seen that one can naturally 'multiply' ideals
of $\mo A$ via the composition operation.
\end{rema}

Now, we will denote by $Z(X,Y)$ for $X,Y\in \obj K(A)$ the subgroup
of $K(A)(X,Y)$ consisting of morphisms that can be presented as
$(s^{i+1}\circ d^i_X+d^{i-1}_Y\circ t^i)$ for some set of
$s^i,t^i\in A(X^i,Y^{i-1})$ (here $X=(X^i),\ Y=(Y^i)$).

\begin{rema}\label{shdef}

 We will often use the fact that $sd+dt=(s-t)d+(dt+td)$ is homotopy
equivalent to $(s-t)d$; hence we may assume that $t=0$ in the
definition of $Z$.
\end{rema}

Now we check that for $Z=\cup_{X,Y\in \obj K(A)} Z(X,Y)$  we have
$Z\nde \mo K(A)$ and $Z^2=0$. A easy standard argument also shows
that
 for any $\cu$ all ideals $Z\nde\mo\cu$  satisfying $Z^2=0$ also possess a
 collection of nice properties.

\begin{lem}\label{lzm}

I1.  $Z\nde\mo K(A)$.

2. Let $L,M,N$ be objects of $K(A)$; let $g\in Z(L,M)\subset K(A)(L,M),\
h\in Z(M,N)\subset K(A)(M,N)$. Then $h\circ g=0$ (in $K(A)$).

II Let $T\nde\mo\cu$  for some additive category $\cu$, suppose also
that $T^2=0$; let $D$ be an additive category. Let $p:\cu\to D$ be
an additive functor such that for any $X,Y\in \obj \cu$ we have $\ke
(\cu(X,Y)\to D(p(X),p(Y)))=T(X,Y)$.

Then the following statements are valid.

1.  Let $p$ be a full functor. Then it is conservative, i.e., $p(g)$
is an isomorphism if and only if $g$ is (for any morphism $g$ in $\cu$).

2. For any $X\in\obj \cu$ and $r\in \cu(X,X)$, if $p(r)$ is an
idempotent then it can
be 
lifted to an idempotent $r'\in \cu(X,X)$ (i.e., $p(r')=p(r)$).

3. If $\cu$ is idempotent complete then its categorical image in $D$
also is. Here we consider a not necessarily full subcategory of $D$
such that all of its objects and morphisms are exactly those that come
from $\cu$.

\end{lem}
\begin{proof}

I1. Obviously, $Z$ is closed with respect to sums and direct sums. 

 Lastly, let $d$ denote the differential, let $f,g$, and $h$ be
 composable morphisms;
 let $g=s\circ d$, for
 $s$ being a collection of arrows shifting the degree by $-1$. Then we
 have $f\circ g=(f\circ s)\circ d$ and
$g\circ h=-(s\circ h)\circ d$; note that $h$ 'anticommutes with the
differential'.

2. Let $L=(L^i),\ M=(M^i),\ N=(N^i)$. Suppose that   for all
$i\in\z$ we have $g^i=s^{i+1}\circ d^i_L$ for some set of $s^i\in
A(L^i,M^{i-1})$, whereas $h^i=u^{i+1}\circ d^i_M$ for some set of
$u^i\in A(M^i,N^{i-1})$

  Then $h^i\circ g^i=u^{i+1}\circ
d^i_M\circ s^{i+1}\circ d^i_L$. Recall now that $g$ is a morphism
of complexes; hence for all $i\in\z$ we have $d^i_M\circ
s^{i+1}\circ d^i_L=d^i_M\circ d^{i-1}_M\circ s^i=0$.  We obtain
that $h\circ g$ is homotopic to $0$.

II1. Since $p$ is a functor, it sends isomorphisms to
isomorphisms.

Now we prove the converse statement. Let $g$ belong to $\cu(X,X')$ for
$X,X'\in \obj \cu$; let $p(h)$ for some  $h\in \cu(X',X)$ be the
inverse to $p(g)$. We have $h\circ g-\id_X\in T(X,X)$ and $g\circ
h-\id_{X'}\in T(X',X')$. It suffices to check that $h\circ g$ and
$g\circ h$ are invertible in $\cu$. The last assertion follows from
equalities $(h\circ g-\id_X)^2=0$ and $(g\circ h-id_{X'})^2=0$ in
$\cu$, that yield $(h\circ g)(2\id_X- h\circ g)=\id_X$ and $(g\circ
h)(2\id_{X'}- g\circ h)=\id_{X'}$.

2. This is just the 
 well-known fact that idempotents can be
lifted (in rings).

We consider  $r'=-2r^3+3r^2$. Since $p(r)^2=p(r)$ in $D$ and
$r'=r+(r^2-r)\circ (id_X-2r)$, we have $p(r')=p(r)$.  Since
$r'^2-r'=(r^2-r)^2\circ (4r^2-4r-3\id_X)$,  we obtain that
$r'$ is an idempotent.

3. The assertion follows immediately from  II2. Indeed, any
idempotent $d$ in the image can be lifted to an idempotent $c$ in
$\cu$. Since $c$ splits in $\cu$, $p(c)=d$ splits in the image.


\end{proof}

\begin{rema}\label{rlzm}
The assertions of part II remain valid for any nilpotent $T$. For
 $l$ that satisfies $l^n=0$, $n>0$, the inverse to $\id_X-l$ is given by
$\id_X+l+l^2+\dots+l^{n-1}$. If $l^n=0$, $r^2-r=l$, then the
equality $(x-(x-1))^{2n-1}=0$ allows one to construct explicitly a
polynomial $P(x)$ such that $P\equiv 0\mod x^n\z[x]$ and $P\equiv
1 \mod (x-1)^n\z[x]$. Then $P(r)^2=P(r)$; $P(r)-r$ can be
factored through $l$.

\end{rema}

\begin{defi}\label{dkw} [The definition of $K_\w(A)$]

We define $K_\w(A) $ as $ K(A)/Z$ (in the sense of Remark
\ref{ridmo}) 
with isomorphic  objects  (i.e., homotopy equivalent complexes)
identified.

We have the obvious shift functor $[1]:K_\w(A)\to K_\w(A)$.

A triangle $X\to Y\to Z\to X[1]$ in $K_\w(A)$ will be called
distinguished if any of its two sides can be  lifted to two sides
of some distinguished triangle   in $K(A)$.

An additive functor $F:\cu\to K_\w(A)$ for a triangulated $\cu$ will
be called {\it weakly exact} if it  commutes with shifts and sends
distinguished triangles to distinguished triangles.

The bounded subcategories  of $K_\w(A)$ are defined in the obvious
way.

\end{defi}

\begin{rema}\label{rkw} [Why $K_\w(A)$ is a category; homology ]

1.  $K_\w(A)$ is a category since we just  factorize the class of
objects of $K(A)/Z$ with modulo a class of invertible morphisms;
see Remark \ref{ridmo}.

2. Let $B $ be an abelian category; let $F:A\to B$ be an additive
functor. Then any $g\in Z(X,Y)$ gives a zero morphism on the homology
of $F_*(X)$. It follows that the homology of $F_*(X)$ gives
well-defined functors $K_\w(A)\to B$. Besides, these functors are
easily seen to be homological, i.e., they send distinguished
triangles in $K_\w(A)$ into long exact sequences.

In particular, this is true for the 'universal' 
functor $A\to A_*'$
(recall that $A_*'$ is the full abelian subcategory of $A_*$
generated by $A$). Hence there are well-defined homology functors
$H_i:K_\w(A)\to A_*'$.

Conversely, suppose that for some $g\in K(A)(X,Y)$ and for any $F:A\to B$ (for an arbitrary $B$) $g$ induces the zero map on homology.
Then Theorem 2.1 of \cite{barrabs} easily yields that $g\in Z(X,Y)$.


3. Now suppose that for a triangulated $\cu$ we have a weakly exact
functor $u:\cu\to K_\w(A)$. Then the homology of $F_*(u(X))$ gives
well-defined functors $\cu\to B$. Again, distinguished triangles in
$\cu$ become long exact sequences.

In particular, this statement can be applied to the weight complex
$t:\cu\to K_\w(\hw)$ (whose functoriality is proved in  Theorem
\ref{fwc} below)(II). This concludes 
the proof
of (part IV of) Theorem
\ref{hsps}.

\end{rema}

Lemma \ref{lzm} immediately yields the following statement.

\begin{pr} \label{wkar}
1. The projection $p:K(A)\to K_\w(A)$ is conservative.

2. Let $A$ be an idempotent complete category. Then $K^b_\w(A)$ is idempotent
complete also. Besides, $K^b_w(A)^{[i,j]}$ is idempotent complete for any $i\le j\in \z$.
\end{pr}

\begin{proof}

1. Immediate from  Lemma \ref{lzm}(II1).

2. It is well known that $K^b(A)$ is idempotent complete; see, for
example, Theorem 2.8 of \cite{ba}. Hence  Lemma \ref{lzm}(II3) yields the
first part of the assertion.

The second part of the assertion can be easily proved using the method of the proof of Proposition 6.4.1 of \cite{mymot}.

\end{proof}

\subsection{The functoriality  of the weight complex}\label{sfwc}

We will use the following simple fact.

\begin{lem}\label{lends}

If $X\in \cu^{w\ge 0}$, $Y\in\cu^{w\le 0}$, then any $f\in
\cu(X,Y)$ can be factored through some morphism $X^0\to Y^0$
(of the zeroth terms of  weight complexes).\end{lem}

\begin{proof} Easy  from the equality $\cu(X^{w\ge
1}[-1],Y)=\cu(X^0,Y^{w\le -1}[1])=\ns$.\end{proof}

Now we prove that  the usual and  'infinite' weight
decompositions define certain functors. Let $X,X'$ denote arbitrary
objects of $\cu$.

\begin{theo}\label{fwc}
I1.  The (single) weight decomposition of objects and morphisms
gives a functor $WD:\cu\to K_\w^{[0,1]}(\cu)$ (certainly, here we only obtain two-term complexes here, and we put   them in degrees $[0,1]$).

2. Morphisms $g\in \cu(X,X')$, $h\in \cu(X^{w\le 0}, X'{}^{w\le 0})$
and $i:\cu(X^{w\ge 1},X'{}^{w\ge 1})$ give a morphism of weight
decompositions (of $X$ and $X'$) if and only if $(h,i)=WD(g)$ in
$K_\w(\cu)$.

3. The homomorphism $\cu(X,X')\to K_\w^{[0,1]}(\cu)(WD(X),WD(X'))$
is surjective.

4. For all $X,X'\in \obj \cu$ consider the groups $$T(X,X')=\ke
(\cu(X,Y)\to K_\w(\hw)(WD(X),WD(X')).$$ Then $T\nde \mo\cu$; $T^2=0$.

5. If $WD(X)\cong WD(X')$ in $K_\w(\cu)$ then $X\cong X'$ in $\cu$.

6. For any $X\in \obj \cu$, $p\in \cu(X,X)$, if $WD(p)$ is
idempotent then $WD(p)$ can be lifted to an idempotent $p'\in
\cu(X,X)$.

II  
The correspondence $X\to t(X)$ gives a functor $\cu\to
K_\w(\hw)$.

\end{theo}
\begin{proof}

1. By  Lemma \ref{fwd}(1), any morphism $X\to X'$  can be
extended to a    morphism of their (fixed) weight decompositions.
This extension is uniquely defined in $K_\w^{[0,1]}(\cu)$ by part 3
of loc. cit.  One can compose such homomorphisms in $K_\w(\cu)$
since one of the possible extensions of the composition of morphisms
$X\to X'\to X''$ (in $C(\cu)$) is the composition of (arbitrary)
extensions for the morphisms $X\to X'$ and $X'\to X''$.

It remains to check that the image of $X$ in $\obj
K^{[0,1]}_\w(\cu)$ does not depend on the choice of the weight
decomposition. Let $K,K'\in \obj K(\cu)$ be given by two weight
decompositions of $X$; $id_X$ induces $g\in  K(\cu)(K,K')$ and $h\in
K(\cu)(K',K)$. By  Lemma \ref{fwd}(3), $h\circ g-\id_K\in
Z(K,K)$ and $g\circ h-id_{K'}\in Z(K',K')$. It suffices to check
that $h\circ g$ and $g\circ h$ are invertible in $K(\cu)$; this
follows from  Proposition \ref{wkar}(1).

2. By definition of $WD$, the triple $(g, WD(g))$ gives a morphism
of weight decompositions.

Now suppose that $(h,i)=WD(g)$, i.e., $(h,i)\in C(\cu)
(WD(X),WD(X'))$ and $(h,i)\equiv WD(g)\mod T(WD(X),WD(X'))$. It
follows that $i\circ f=f'\circ h$ (in the notation of
(\ref{wdeck})). Besides, there exist $(h',i')$ that give a
morphism of weight decompositions; $h-h'=s\circ f$ and
$i-i'=f'\circ t$ for some $s,t\in \cu(X^{w\ge 1},X'{}^{w\le 0})$.
We obtain that $h\circ a =h'\circ a=a'\circ g$ and $b' \circ i=b'
\circ i'=g[1]\circ b$.

Hence $(g,h,i)$ give a morphism $T_0\to T'_0$.

3. By definition, any $h\in K_\w^{[0,1]}(\cu)(WD(X),WD(X'))$ comes
from some commutative square
$$\begin{CD}
X^{w\le 0}@>{f^0}>>X^{w\ge 1}\\
@VV{}V@VV{}V \\
X'{}^{w\le 0}@>{f'^0}>>X'{}^{w\ge 1}
\end{CD}$$
 Extending this square to a
morphisms of  triangles $T_0\to T'_0$ (i.e., of weight
decompositions of $X$ and $X'$) immediately yields the result.

4. Since $WD$ is a functor, $T$ is an ideal.

We prove that $T^2=0$  similarly to the proof of  Lemma
\ref{lzm}(I2).

 Let $X,X',X''$ be objects of $ \cu$; let $g\in T(X',X'')\subset \cu(X,X'),\
h\in T(X',X'')\subset \cu(X',X'')$.

We should check that $h\circ g=0$ (in $\cu$). We can choose any
weight decompositions of $X,X',X'$; denote them by $T,T,T''$
(similarly to (\ref{wdeck})).

Since $WD(g)=WD(h)=0$, by assertion I2 we obtain that $(g,0,0)$
 and $(h,0,0)$ give morphisms of weight decompositions. This means
 that $ a'\circ g=a''\circ  h=g[1]\circ b= h[1]\circ b'=0$.
 Hence $g$ can be presented as $b'[-1]\circ c$ for some
$c\in \cu(X, X'{}^{w\ge
 1}[-1])$. Then $h\circ g= (h[1]\circ b')[-1]\circ c=0$.


5. By assertion I3 any isomorphism $WD(X)\to WD(X')$ is induced by
some morphism $X\to X'$. Now by  Lemma \ref{lzm}(II1), $t$
is conservative (we apply assertion I4); this yields the result.

6. Immediate from  Lemma \ref{lzm}(II2).

II Exactly the same reasoning as in assertion I1 will prove the assertion
after we verify that morphisms in $\cu$ give well-defined morphisms
of weight complexes (in $K_\w(\hw)$).

A $g\in\cu(X,X')$ 
can be extended to a morphism $Po(X)\to Po(X')$ (see  Remark \ref{rpoto}(1));
hence we also obtain some morphism $t(g):t(X)\to t(X')$. It remains to
verify that for $g=0$ we have $t(g)\in Z(t(X),t(X'))$.
We use the notation of  Proposition \ref{mwdec}(2).

We study the possibilities for $g^i:X^i\to X'{}^i$ (that we choose to be compatible with (\ref{wdeck3}), without considering (\ref{wdeck4})). By
construction, $g^i$ depends  on the maps $r^k:X^{w\le k}\to X'{}^{w\le
k}$ only for $k=i,i-1$. This dependence is linear. Moreover, any
pair of $(r^i,r^{i-1})$ can be presented as $(0,r^{i-1})+(r^i,0)$.
Indeed, for $g=0$ any of $r^k$ may be zero, whereas distinct
$r^k$ are 'independent' by  Proposition \ref{mwdec}(1). Hence
it suffices to prove that $g^i$ can be presented as $(s^{i+1}\circ
h^{iX}+h^{i-1,X'}\circ t^i)$ for some $s^{i+1}\in
\hw(X^{i+1},X'{}^{i}),\ t^i\in \hw(X^i,X'{}^{i-1})$ in the following cases:
either $r^i$ or $r^{i-1}$ equals $0$ (recall that $h$ denotes the corresponding boundaries
 of  weight complexes).

In the case $r^i=0$ we can present $g^i[-1]$ as the second
component of $WD(0:X^{w\le i}[-1]\to X'{}^{w\le i}[-1])$. Hence
$g^i$ equals $c^{i-1,X'}\circ u^i$ for some $u^i\in \cu(X^i,\
X'{}^{w\le i-1})$ (by assertion I1). Note now that $u^i$ can be factored through
$X'{}^{i-1}$ (see Lemma \ref{lends}).

In the case $r^{i-1}=0$ we can present $g^i$ as the first
component of $WD(0:X^{w\ge i}\to X'{}^{w\ge i})$. Hence $g^i$
equals $v^{i+1}\circ x^{iX}[1]$ for some $v^{i+1}\in \cu(X^{w\ge
i+1},X'{}^{i})$. It remains to note that $v^i$ can be factored
through $X^{i+1}$.

Combining the two cases, we obtain our claim.


\end{proof}

\begin{rema}\label{isogt}
The functoriality of $t$ implies that for any $X\in\obj \cu$ any two
choices for $t(X)$ are connected by a (canonical) isomorphism in
$K_\w(\hw)$. Then  Lemma \ref{lzm}(II1) (combined with part I2
of the Lemma) implies that they are isomorphic (not
necessarily canonically) 
 in $K(\hw)$, i.e., they are homotopy equivalent (in $C(\hw)$).

\end{rema}

$WD$ and $t$ commute in the following sense.

\begin{pr}\label{sogl} Let $X,X'$ be objects of $\cu,\ g\in \cu(X,X')$.

1. Any  choice of $(t(i),t(l))$ for $(i,l)=WD(g)$  comes from a
truncation of $t(g)$ (here we fix some weight decompositions of
$X$ and $X'$ and consider all compatible lifts of $t(g)$ to $\mo
C(\hw)$).

2. Suppose that $(r',s')=(t(i'),t(l'))$ for some weight decomposition
$(i',l')$ of $g$, let $r+ s:t(X)\to t(X')$ be homotopic to $r'+
s'$ (here we consider sums of collections of arrows). Then
$(r,s)=(t(i),t(l))$ for some (other) weight decomposition $(i,l)$
of $g$.
\end{pr}

\begin{proof}

1. By the definition of $t(g)$ (see  Theorem \ref{fwc}(II))
any choice of $(t(i),t(l))$
 is a possible truncation of $t(h)$ over $C(\hw)$.

2. It suffices to prove the statement for $g=0$. Suppose that
$(r,s)$ can be obtained from some $WD(0)$ via $t$.
 Note  that (replacing $r,s$ by  equivalent
morphisms if needed) we can assume that $r=r^0$, $s=s^1$ (i.e., they
are concentrated in degrees $0,1$). Hence
 there exists some $l\in \hw (X^1,X'{}^{0})$ such that
$r^0=l\circ h^0$, $s^1=h'^0\circ l$.

Now it remains to note that the triple $(0,d'^0\circ l\circ
c^1,x'^0[1]\circ l\circ y^1)$ gives a weight decomposition of
$0:X\to X'$. This fact follows from the equalities $ d'^0\circ
l\circ c^1\circ a^0=0=b'^0\circ x'^0[1]\circ l\circ y^1$ (see
(\ref{wdeck3}) and (\ref{wdeck4})), whereas
$$f'^0\circ d'^0\circ l\circ c^1=x'^0[1]\circ l\circ c^1=
 x'^0[1]\circ l\circ y^1\circ f^0.$$
\end{proof}

\subsection{Main properties of the weight complex }\label{smpwc}

Now we  prove the main properties of the weight complex functor.

\begin{theo} \label{wecomp}[The weight complex theorem]

I Exactness.

$t$ is a weakly exact functor.

II Nilpotency.

$I(-,-)=\ke (\cu(-,-)\to K_\w(t(-),t(-)))$ defines an ideal in $\mo
\cu$. For any $i\le j\in\z$ the restriction $I^{[i,j]}$ of $I$ to
$\cu^{[i,j]}$ satisfies $(I^{[i,j]}){}^{j-i+1}=0$.

III  Idempotents.

If $X\in \cu^b$, $g\in \cu(X,X)$, $t(g)=t(g\circ g)$, then $t(g)$
can be lifted to an idempotent $g'\in \cu(X,X)$.

IV Filtration.

If $X\in \cu^{w\le i}$ (resp. $\cu^{w\ge i}$) for some $i\in\z$ then
$t(X)\in K_\w(\hw)^{\le i}$ (resp. $K_\w(\hw)^{\ge i}$), i.e., it is
homotopy equivalent to a complex concentrated in degrees $\le i$
(resp. $\ge i$).

If $X$ is bounded  from above (resp. from below) then the converse
implications are valid also.

V  Conservativity.

If $w$ is non-degenerate, then the functor $t$ is conservative on
$\cu^+$ and $\cu^-$. 

 VI If $X,Y\in \cu^{[0,1]}$ then $t(X)\cong t(Y)$ implies $X\cong
 Y$.

VII Let $X$ belong to $\cu^{w\ge a}$ for some $a\in \z$; consider the
homomorphism $t_*: \cu(X,X')\to K_\w(\hw)(t(X),t(X'))$. Then the
following statements are valid.

1. If $X'\in \cu^{w\le a}$ then $t_*$ is bijective.

2. If $X'\in \cu^{w\le a+1}$ then $t_*$ is surjective.

\end{theo}

\begin{proof}

I Let $C\stackrel{a}{\to} X\stackrel{f}{\to}
X'\stackrel{b}{\to}C[1]$ be a distinguished triangle. We should
prove that the triangle $t(C)\stackrel{t(a)}{\to}
t(X)\stackrel{t(f)}{\to} t(X')\stackrel{t(b)}{\to}t(C)[1]$  is
distinguished. It suffices to construct a triangle of morphisms
\begin{equation}\label{eqv}
V:t(X'[-1])\stackrel{m}{\to} t(C) \stackrel{n}{\to}
t(X)\end{equation}
 that
splits componentwisely (in $C(\hw)$) such that $m$ is some choice for
$t(b)[-1]$ and  $n$ is some choice for $t(a)$. Indeed, it is a well
known fact that any such $V$ gives a distinguished triangle in
$K(\hw)$. Hence any two sides of $t(V)$ can be lifted to two sides
of a distinguished triangle in $K(\hw)$; so $t(V)$ is distinguished
(see Definition \ref{dkw}).

In order to prove our claim we construct 'nice' weight
decompositions of  $C[i]$ for all $i\in \z$.
To this end we apply the method used in the proof of  Proposition \ref{trfun}(II1)
for all $k\in\z$.

We apply Lemma \ref{twd} for all
$i\in\z$. By the lemma, the triangles $C[i]\to C_i\to C_i'$
obtained from (\ref{dia3na3t}) by shifting the last row are weight
decompositions of $C[i]$ for all $i\in\z$. Hence the first two
 columns of (\ref{dia3na3t}) can be completed to  morphisms of the analogues of (\ref{wdeck3}) for  $X'$,  $C[1]$, and  
 $X[1]$ (for all $k\in \z$).

Now, in order to 'connect' shifted weight decompositions we use exactly
the same reasoning as in the proof of  Proposition \ref{trfun}(II1). 
 We obtain that the corresponding map of weight complexes
 splits componentwisely indeed.

II  $I$ is an ideal since $t$ is an additive functor.

 Obviously, it suffices to check that for $X\in
\cu^{[0,n]}$ the ideal $J=\{g\in \cu(X,X):\ t(g)=0\}$ of the ring
$\cu(X,X)$) satisfies  $J^{n+1}=0$. We will prove this fact by
induction on $n$. 
In the case $n=0$ we have $\cu^{[0,n]}=\hw$,
hence $J=\ns$.

To make the inductive step we consider $g_0\circ g_1\circ\dots g_n$,
$g_i\in J$, let $r=(g_0\circ g_1\circ\dots g_{n-1})[n-1]$,
$s=g_n[n-1]\circ r$. By Proposition \ref{sogl}, for all $i$ we can choose 
representatives $(h_i,l_i)$ of $WD(g_i[n-1])$  such that $t(h_i)=0$.
Then by the inductive assumption we have $WD(r)=(0,m)$ for some
$m:X^n\to X^n$. Considering the morphism of triangles corresponding
to $WD(r)$ we obtain that $r=b_{n-1}[-1]\circ q$ for some
$q:X[n-1]\to X^n[-1]$. Next, since $t(g_n)=0$, we can assume that
$t(g_n[n-1])=(u,0)$ for some $u$ (by Proposition \ref{sogl}). Hence
$g_n[n-1]=v\circ a^{n-1}$ for some $v\in \cu(X^{w\le n-1},X[n-1])$
and we obtain $s=v\circ (a^{n-1}\circ b^{n-1}[-1])\circ q =0$. The
assertion is proved.

III Follows from assertion II by a standard reasoning; see Remark
\ref{rlzm}.

IV By  Proposition \ref{fprwc}(2), if $X\in \cu^{w\le i}$
(resp. $X\in \cu^{w\ge i}$) then choosing $X^{w\ge i+1}=0$  (resp.
$X^{w\le i-1}=0$) we obtain that the corresponding choice of
$t(X)$ is concentrated in degrees $\le i$ (resp. $\ge i$). Now
note that all choices of $t(X)$ are homotopy equivalent by 
Proposition \ref{wkar}(1).

Conversely, let $w$ be non-degenerate, let $t(X)\in K_\w(\hw)^{w\le
i}$. We can assume that $i=0$; let $X\in \cu^{w\le n}$
(for some $n\ge 0$).
Then
$t(\id_X)$ is  equivalent to a morphism whose non-zero components are in
degrees $\le 0$. Hence Proposition \ref{sogl} implies that for
$WD(\id_X)=(l,m)$ we can assume that $t(m)=0$.

Then by assertion II we have $WD(\id_X^n)=(l^n,0)$. Considering
the distinguished triangle corresponding to $WD(\id_X^n)$ we
obtain that $id_X=id_X^n$ can be factored through $X^{w\le
0}$. Hence $X$ is a retract of $X^{w\le 0}$; since $\cu^{w\le 0}$
is Karoubi-closed in $\cu$ we obtain that $X\in \cu^{w\le 0}$.

The case $t(X)\in K_\w(\hw)^{w\ge i}$ is considered similarly.

V Since $t$ is weakly exact (see Definition \ref{dkw}), it
suffices to check that $t(X)=0$ implies $X=0$. This is immediate
from assertion IV.

VI Immediate from  Theorem \ref{fwc}(I5).

VII We can assume that $a=0$.

1.  The proof is just a repetitive application of axioms (of
weight structures).

Note first that $t_*$ is bijective for $X,X'\in \cu^{w=0}$. Next,
for $X\in \cu^{w=0}$ and any $X'$ we consider the distinguished
triangle $X'{}^{w\le -1}\to X'{}^0 \to X'\to X'{}^{w\le -1}[1]$.
Then orthogonality yields that any $h:\cu (X,X'{}^0)$ gives a
morphism $X\to X'$; hence $t$ is surjective in this case. We also
can apply this statement for $X''=X'{}^{w\le -1}$.  Hence
considering the diagram $$\begin{CD} (\cu(X,X'{}^{w\le -1})@>{}>>
\cu(X,X'{}^0)@>{}>> \cu(X,X')@>{}>> 0\\
@VV{}V@VV{}V@ VV{}V \\
K_\w(\hw)(t(X),t(X'{}^{w\le -1}))@>{}>>
K_\w(\hw)(t(X),t(X'{}^0))@>{}>> K_\w(\hw)(t(X),t(X'))@>{}>> 0
\end{CD}
$$ induced by $t$ we obtain that $t_*$ is bijective in this case.

Now considering the distinguished triangle $X^{w\ge 1}[-1]\to X\to
X{}^0 \to X^{w\ge 1}$ and applying the dual argument one can
easily obtain the claim.

2. Let $h$ belong to $K_\w(\hw)(t(X),t(X'))$. By definition, we can 'cut' $h$
to obtain a commutative diagram $$\begin{CD}
t(X^{0})@>{t(f^0)}>>t(X^{w\ge 1})\\
@VV{}V@VV{}V \\
t(X'){}^{w\le 0}@>{t(f'^0)}>>t(X'{}^{1})
\end{CD}$$

By assertion VII1, this diagram  corresponds to some homomorphism
$WD(X)\to WD(X')$.  It remains to apply Theorem
\ref{fwc}(I3). 
\end{proof}

\begin{rema}\label{virtri}

1. By assertions IV and V, 
$t$ is always conservative and strictly
respects weight filtration on $\cu^b$.

2. In fact, our restrictions on a distinguished triangle in $K_\w(A)$
(see Definition \ref{dkw}) are rather weak. Our definition
is similar to the
notion of  {\it exact triangle} in Definition 0.3 of
\cite{virttriang}. Since exact triangles are not distinguished in
general (see loc. cit.), part I of our theorem does not imply that
for a distinguished triangle $C\to X\to X'$ the triangle $t(C)\to
t(X)\to t(X')$ comes from some distinguished triangle in $K(\hw)$.

We will not prove the latter fact in detail, since we will only
need it in Remark \ref{rctst} below. Yet the proof is rather easy.
In the proof of part I of our theorem it suffices to check that
some choice of $t(f)$ in $K(\hw)$ yields the third side of the
distinguished triangle in question. Using (obvious) functoriality
properties of the construction in the proof, one can reduce the
latter claim to the case $X,X'\in \cu^{w=0}$. Certainly, the
statement is obvious in this case.

\end{rema}

Probably $t$ can be 
lifted to a certain 'strong'
weight complex functor.

\begin{conj}\label{ctst}
$t$ can be lifted to an exact functor $t^{st}: \cu\to K(\hw)$.
\end{conj}

\begin{rema}\label{rctst}

1. Suppose that $\hw$  is (fully) embedded into the subcategory $B= Proj\ A$
of projective objects  of an abelian category $A$ (a reasonable choice for $A$ is $\hw_*'$; cf. Lemma \ref{lreg}
below). Then we have a full embedding $K(\hw)\subset K(B)$ and a canonical functor $ K(B)\to
D(A)$. Suppose now that $A$ is of projective dimension $1$. Then the latter functor is an equivalence and any
complex over $A$ is quasi-isomorphic to a complex with zero
differentials; hence it can be presented in $D(A)$ as a direct sum
of 
some monomorphisms in $B$ (i.e., of complexes of the form $\dots
0\to X\stackrel{i}{\hookrightarrow}Y\to 0\to\dots$ placed in 
pairwise distinct
dimensions). We check that $K_\w(B)=K(B)$. Note that it suffices
to prove the corresponding fact for $K^b(-)$. Therefore it suffices
to check that $K_\w(B)(X,Y)=K(B)(X,Y)$ for $X,Y$ being
monomorphisms (as two-term complexes); let
$X=X^{-1}\hookrightarrow 
 X^0$. If
$Y\in C^{[-1,0]}(B) $ then
$K(B)(X,Y)=A(H^0(X),H^0(Y))=K_\w(B)(X,Y)$ (see Remark
\ref{rkw}(2)). If $Y\in C^{[-2,-1]}(B) $ then the equality
$K(B)(X,Y)=K_\w(B)(X,Y)$ is obvious (cf. Theorem
\ref{wecomp}(VII)). For $Y$ placed in all other positions we have
$K(B)(X,Y)=\ns=K_\w(B)(X,Y)$.

 We conclude that
$K_\w(\hw)=K(\hw)$. Therefore Remark \ref{virtri}(2) implies
that $t$ is exact (as a functor of triangulated categories).

 In particular, this reasoning can be applied if
$\hw=\ab_{fin.fr}$ or $\hw=\ab_{fr}$. 
 Hence this is
the case for all categories of spectra considered in \S\ref{sh}
below. Note here that Example 2.3 of \cite{barrabs} (that states that our abelian group observation is wrong) is erroneous, since the sequence of arrows $dh$ considered in it does not yield a morphism of complexes.

2. In  \S \ref{cwecom} below we will also verify the conjecture in
the case when $\cu$ has a {\it differential graded enhancement}.

3. Prof. A. Beilinson has kindly communicated to the author a
proof of the  conjecture in the case when $\cu$ has a {\it
filtered triangulated enhancement}; see  \S\ref{sftriang} below.
Probably, a filtered triangulated enhancement exists for any
'reasonable' triangulated category.
\end{rema}


\section{ Constructing weight structures and adjacent $t$-structures}\label{dualh}

In this section (especially in \S\ref{sadst}) we
prove that weight structures are closely
related to $t$-structures.

In \S\ref{dtst} we recall the definition of a
$t$-structure in a triangulated $\cu$.
In  \S\ref{colim} we recall the (standard) construction of countable
homotopy colimits in  triangulated categories and study its
properties.

In \S\ref{srecw} we show that in many cases a weight structure can
be described by specifying a {\it negative} $H('\approx' \hw)\subset
\cu$. In particular, this is the case for the category of finite
spectra ($\subset SH$), and for $\chow\subset \dmgm$ (that we will discuss later).

In \S\ref{sadst} we define the notion of (left or right)
{\it adjacent} weight and
$t$-structures for $\cu$; their hearts are dual in a very
interesting sense (see Theorem \ref{sdt}). It turns out that the
truncations of an object $Y$ with respect to a $t$-structure that is
adjacent to $w$ represent exactly the virtual $t$-truncations of the
functor $\cu(-,Y)$ (with respect to $w$, see Remark \ref{rtrfun}).
 Hence  spectral
sequences arising from adjacent weight and $t$-structures are closely
related.
Lastly, a functor of triangulated categories is $t$-exact
(with respect to some $t$-structures) if and only if its (left) adjoint is
 {\it weight-exact} with respect to  weight structures that are
(left) adjacent to these $t$-structures.

In \S\ref{eads} we study the conditions for adjacent weight and
$t$-structures to exist. We only consider in details the cases which
are relevant for our main examples (motives and spectra); other
possibilities are mentioned in Remark \ref{cocomp}. 

In \ref{sh} we apply the results of this section to the study of
$SH$. In particular, we construct a {\it spherical} weight structure
for it; it is adjacent to the Postnikov $t$-structure.

In \S \ref{ddme} below we will apply our results to $\dme$ (the
category of motivic complexes of Voevodsky).

\subsection{$t$-structures: reminder}
\label{dtst}

To fix the notation we recall the definition of a $t$-structure.

\begin{defi}\label{dtstr}

A pair of subclasses  $\cu^{t\ge 0},\cu^{t\le 0}\subset\obj \cu$
 will be said to define a
$t$-structure $t$ if 
they satisfy the
following conditions:

(i) $\cu^{t\ge 0},\cu^{t\le 0}$ are strict, i.e., contain all
objects of $\cu$ isomorphic to their elements.

(ii) $\cu^{t\ge 0}\subset \cu^{t\ge 0}[1]$, $\cu^{t\le
0}[1]\subset \cu^{t\le 0}$.

(iii) {\bf Orthogonality}. $\cu^{t\le 0}[1]\perp
\cu^{t\ge 0}$.

(iv) {\bf $t$-decompositions}.

For any $X\in\obj \cu$ there exists
a distinguished triangle
\begin{equation}\label{tdec}
A\to X\to B{\to} A[1]
\end{equation} such that $A\in \cu^{t\le 0}, B\in \cu^{t\ge 0}[-1]$.

\end{defi}

Non-degenerate and bounded (above, below, or both) $t$-structures
can be defined similarly to Definition \ref{d2}.

We will need some more notation for $t$-structures.

\begin{defi} \label{dt2}

1. A category $\hrt$ whose objects are $\cu^{t=0}=\cu^{t\ge 0}\cap
\cu^{t\le 0}$, $\hrt(X,Y)=\cu(X,Y)$ for $X,Y\in \cu^{t=0}$,
 will be called the {\it heart} of
$t$. Recall (cf. Theorem 1.3.6 of \cite{BBD}) that $\hrt$ is always
abelian; short exact sequences in $\hrt$ come from distinguished
triangles in
$\cu$.

2. $\cu^{t\ge l}$ (resp. $\cu^{t\le l}$) will denote $\cu^{t\ge
0}[-l]$ (resp. $\cu^{t\le 0}[-l]$).

\end{defi}

\begin{rema}\label{rts}
1. Recall (cf. Lemma IV.4.5 in \cite{gelman}) that (\ref{tdec})
defines additive functors $\cu\to \cu^{t\le 0}:X\to A$ and $C\to
\cu^{t\ge 1}:X\to B$. We will denote $A,B$ by $X^{t\le 0}$ and
$X^{t\ge 1}[-1]$, respectively. (\ref{tdec}) will be called the
{\it t-decomposition} of $X$.

More generally, the $t$-components of $X[i]$ will be denoted by
$X^{t\le i}\in
\cu^{t\le 0}$ and $X^{t\ge i+1}[-1]\in \cu^{t\ge 1}$, respectively.

2. 
The functor $X\mapsto X^{t\ge 1}[1]$ is  left adjoint to the inclusion $\cu^{t\ge 1}\to \cu$.
It follows that this functor
commutes with all  coproducts that exist in $\cu$. Besides,
if $\coprod X_i$, $\coprod X_i^{t\le 0}$, and $\coprod X_i^{t\ge 1}$
exist in $\cu$, then the distinguished triangle $\coprod X_i^{t\le
0}\to \coprod X_i\to \coprod X_i^{t\ge 1}[-1]$ (here we apply the dual of Proposition 1.2.1 of \cite{neebook}) yields that $(\coprod
X_i)^{t\le 0}=\coprod (X_i^{t\le 0})$. 

\end{rema}

We denote by $H^t_{0}$ the zeroth homology functor corresponding to
$t$ (cf.  \S IV.4(10)  of \cite{gelman}), i.e., $H^t_{0}(X)$ is
defined similarly to $X^{[0,1]}$ in Proposition
\ref{mwdec}(1). Shifting the $t$-decomposition of $X^{t\le 0}[-1]$ by
$[1]$ we obtain  a canonical and functorial (with respect to $X$)
distinguished triangle $X^{t\le -1}[1] \to X^{t\le 0}\to H^t_{0}(X)$.

Lastly, $\tau_{\le i}X$ will denote $X^{t\le i}[-i]$; $\tau_{\ge
i}X=X^{t\ge i}[-i]$; $H^t_i=H_0^t(X[i])$.

\subsection{Countable homotopy colimits in triangulated categories:
 the construction and
properties}\label{colim}

 The triangulated construction of countable (filtered)
homotopy colimits is fairly standard, cf. Definition 1.6.4 of
\cite{neebook}.

\begin{defi}\label{dcoulim}

Suppose that we have a sequence of objects $Y_i$ (starting from some
$j\in\z$)  and maps $\phi_i:Y_{i}\to Y_{i+1}$.  Suppose that there exists
$D=\coprod Y_i$ 
 in $\cu$. We consider the map $d:\oplus
\id_{Y_i}\bigoplus \oplus (-\phi_i): D\to D$ (we can define it since
its $i$-th component 
 can be easily factorized  as the composition
$Y_i\to Y_i\bigoplus Y_{i+1}\to D$).
 Denote  a cone of $d$ as $Y$. We will write $Y=\inli Y_i$ and
 call $Y$ the {\it homotopy colimit} of $Y_i$; we will not consider any other homotopy colimits in this paper.

 We will say that
 the colimit exists (in $\cu$) if the coproduct $D$ exists.
\end{defi}

\begin{rema}\label{rcoulim}

1. By Lemma 1.7.1 of \cite{neebook} the homotopy colimit of
$Y_{i_j}$ is the same for any subsequence of $Y_i$. In particular,
we can discard any (finite) number of first terms in $Y_i$.

2. By Lemma 1.6.6 of \cite{neebook} the homotopy colimit of
$X\stackrel{\id_X}{\to}X\stackrel{\id_X}{\to}
X\stackrel{\id_X}{\to} X\stackrel{\id_X}{\to}\dots$ is $X$. Hence
we obtain that $\inli X_i\cong X$ if for $i\gg 0$ all $\phi_i$ are
isomorphisms and $X_i\cong X$.

3. The construction of $\inli Y_i$ easily yields the following: if countable
coproducts exist in $\cu^{w\le 0}$ then $\cu^{w\le 0}$ is closed
(in $\cu$) with respect to (countable) homotopy colimits. Indeed, we have $D\in
\cu^{w\le 0}$; hence it suffices to recall that  $\cu^{w\le 0}$  is
extension-stable (see  Proposition \ref{bw}(3)). On the other
hand, it is easy to construct a counterexample to the similar
statement for $\cu^{w\ge 0}$ (though
  countable colimits of objects of $\cu^{w\ge 0}$
   always belong to $\cu^{w\ge -1}$). To settle this
problem below we will describe a 'clever' method for passing to the
colimit in $\cu^{w\ge 0}$.

\end{rema}

We study the behaviour of colimits  under (co)representable
functors.

\begin{lem}\label{coulim} 
1.  For any $C\in\obj\cu$ we have a natural surjection $
\cu(Y,C)\to \prli \cu(Y_i,C)$.

2. This map is bijective if all $\phi_i[1]^*:
\cu(Y_{i+1}[1],C)\to \cu(Y_i[1],C)$  are surjective for all $i\gg
0$.

3. If $C$ is compact then $\cu(C,Y)= \inli \cu(C,Y_i)$.
\end{lem}
\begin{proof}

 1. For any $C$ we have
 $\cu(D,C)=\prod \cu(Y_i,C)$.

 This yields a long exact sequence
  $$\dots\to \cu(D[1],C)\stackrel{a[1]^*} \to\cu(D[1],C)\to
 \cu(Y,C)\to \cu(D,C)\stackrel{a^*}{\to} \cu(D,C)\to\dots.$$

  It is easily seen that the kernel of $a^*$ equals
$$\{(s^i): s^i\in \cu (Y_i,C),\ s^{i+1}=s^i\circ \phi_i\}=\prli
  \cu(Y_i,C);$$ this yields the result.

2. By Remark \ref{rcoulim}(1),  we can assume that the
homomorphisms $\phi[1]^*$ are surjective for all $i$. In this case
$a[1]^*$ is easily seen to be surjective; this yields the result.

3. Similarly to the proof of assertion 1, we consider the long exact
sequence
$$\dots\to \cu(C,D)\stackrel{a_*}{\to} \cu(C,D)\to \cu(C,Y)\to
 \cu(C,D[1])\stackrel{a[1]_*} \to\cu(C,D[1]) \to\dots.$$

Since $C$ is compact, we have $\cu(C,D)=\bigoplus \cu(C,Y_i)$. Then it
is easily seen that  $a[1]_*$ is surjective, whereas the cokernel
of $a_*$ is $\inli \cu(C,Y_i)$. See also Lemma 2.8 of \cite{neeo}.

\end{proof}

Now we describe a 'clever' method for passing to the colimit in
$\cu^{w\ge 0}$. Since we will use it to prove that a certain
candidate for being a weight structure is a weight structure
indeed, we will describe it in a (somewhat) more general setting
than that of weight structures.

Suppose that we have a  full extension-stable
(see Definition \ref{exstab}) subcategory $D\subset \cu$.
Define a full subcategory $E\subset\cu$ by $\obj E={}^\perp (D[1])$.
Note that $E$ is also extension-stable, whereas the pair $(D,E)$
satisfies the conditions of Remark \ref{rtwd}.

\begin{lem}\label{gcoulim}

Let $\phi_i: Y_i\to Y_{i+1}$ be a sequence of $\cu$-morphisms;
denote $\co \phi_i$ by $C_i$; let the first of $Y_i$  be $Y_l$.
Suppose that  $Y_l$ and all $C_i$ have 'weight decompositions with
respect to $D,E$', i.e., that there exist distinguished triangles
$Y_l\to D_l\to H_l$ and $C_i\to F_i\to G_i$ with $D_l,F_i\in\obj
D$ and $H_l,G_i\in\obj E$. Suppose also that for any possible
choice of  'weight  decompositions' of $Y_i$ ($Y_i\to D_i\to E_i$
with $D_i\in\obj D$ and $E_i\in\obj E$) the coproduct $\coprod E_i$ 
exists. Then there exists a choice of $E_i$ and of the morphisms
$\phi'_i:E_i\to E_{i+1}$ compatible with $\phi_i$
such that $\inli E_i\in\obj E$ (note
that the colimit exists!).
\end{lem}
\begin{proof}
We  fix  weight decompositions for all $C_i$ and for   $Y_l$.

Next we fix $\phi'_i$ and the weight decompositions of  $Y_{i+1}$
starting from $i=l$ inductively.

Suppose that we have fixed some weight decomposition of  $Y_i$. By
Remark \ref{rtwd} we can construct $E_{i+1}$ and $\phi'_i$ that fit
into a distinguished triangle $E_i\stackrel{\phi'_i}{\to} E_{i+1}\to
G_i$, whereas $E_{i+1}[-1]\to Y_{i+1}$ yields a weight decomposition
of $Y_i$.

Now we check that passing to the limit of $E_i$ this way we obtain
an object of $E$. Let $Z$ be the limit of $E_i$. We should check
that $Z\perp  D[1]$. By Lemma
\ref{coulim}(2), to this end it suffices to check that all
$\phi'_i[1]^*: \cu(E_{i+1},C)\to \cu(E_{i}^{w\ge 1}[1],C)$ are
surjective. Indeed, then we will have $\cu(Z,C)=\prli
\cu(E_{i},C)=\ns$.
 Lastly, the surjectivity is immediate from the
long exact sequences (for all $i$) $$\dots\to \cu(E_{i+1}[1],C) \to
\cu(E_{i}[1],C)\to \cu(G_i,C)(=\ns)\to\dots.$$

\end{proof}

\begin{rema} 
Note that $t$-structures do not have a similar property
(since there is no flexibility for 
$t$-decompositions).
\end{rema}

Lastly we prove that $t$-truncations 'approximate' objects. We
will prove this statement in the form that is relevant for
\S\ref{ddme}; certainly, some other versions of it are valid for
similar reasons.

\begin{lem}\label{tcoulim}

Let $t$ be a non-degenerate $t$-structure.
Suppose that all
countable coproducts exist in $\cu^{t\le 0}$; suppose also that
(countable) coproducts respect $t$-decompositions of objects
 of $\cup_{i\in\z}\cu^{t\le i}$. 

Suppose that all $Y_i$ belong to $\cu^{t\le l}$ for some $l\in
\z$; let $\phi_i: Y_i\to Y_{i+1}$ be a sequence of
$\cu$-morphisms.
Suppose that there exists an $Y$ such that for any $j\in \z$ and all
$i\ge j$ we have $Y_i^{t\ge - j}\cong Y^{t\ge - j}$  and these
isomorphisms commute with ${\phi_i}_*:\ Y_i^{t\ge - j}\to Y_{i+1}^{t\ge
-j}$. Then  $\inli Y_i$ exists and $\cong Y$.

\end{lem}
\begin{proof}

Since countable coproducts exist in $\cu^{t\le 0}$, they also
exist in $\cu^{t\le l}$. This implies the existence of $\inli
Y_i$. We denote $\inli Y_i$ by $Z$.

We obviously have $Y\in \cu^{t\le l}$. Then the definition of
$\inli$ easily yields (at least, one) morphism $\inli
(Y\stackrel{\id_{Y}}{\to}Y
\stackrel{\id_{Y}}{\to}Y\stackrel{\id_{Y}}{\to}\dots)=Y\to \inli
Y_i=Z$; to this end one should apply Proposition \ref{3na3}.

Since $t$ is non-degenerate, it suffices to prove that
$H^t_k(Y)\cong H^t_k(Z)$ for any $k\in\z$. We fix $k$.

Since countable coproducts respect $t$-homology in question,
we have a long exact sequence (in $\hrt$):
$\dots \to \coprod H^t_k(Y_i)
 \to H^t_k(Z)
\to \coprod H_{k+1}^t(Y_i) \to\dots$. We can also obtain a similar
long exact sequence for $H_k^t(Z)$ (with $Y_i$ replaced by $Y$)
 if we present it as
$\inli
(Y\stackrel{\id_{Y}}{\to}Y
\stackrel{\id_{Y}}{\to}Y\stackrel{\id_{Y}}{\to}\dots)$.
Now,
by  Remark \ref{rcoulim}(1) we can assume that
all $H_{k}^t(Y_i)\cong H_{k}^t(Y)$ and
  $H_{k+1}^t(Y_i)\cong H_{k+1}^t(Y)$. This concludes the proof.

\end{proof}

\subsection{Recovering $w$ from its heart}\label{srecw}

It is often the case that instead of describing $\cu^{w\le 0}$ and $\cu^{w\ge
0}$, it is easier to specify only $\cu^{w=0}$. We describe some
conditions ensuring that $w$ can be recovered from $\hw$ (Theorem \ref{recw}(II) is especially easy to apply). We will need the following definitions.

\begin{defi}\label{negth}

Let $H$ be a 
full additive subcategory of $\cu$.

1.  We will say that $H$ is {\it negative} if
 $\obj H\perp (\cup_{i>0}\obj (H[i]))$.

2. We define the {\it small envelope} of an additive category $A$
as a category $A'$ whose objects are $(X,p)$ for $X\in\obj A$ and
$p\in A(X,X)$ such that  $ p^2=p$ and there exist $Y\in \obj A$ and
$q\in A(X,Y)$, $s\in A(Y,X)$ satisfying $sq=1-p$, $qs=\id_Y$. We
define
\begin{equation}\label{mthen}
A'((X,p),(X',p'))=\{f\in A(X,X'):\ p'f=fp=f \}.\end{equation}
\end{defi}

The small envelope of $A$ is (naturally) a full subcategory of the
idempotent completion of $A$ (cf. \S\ref{ridcomp} below).
One should think of $A'$ as of the category of $X\ominus Y$ for
$X,Y\in\obj A$, $Y$ is a retract $X$. Here $X\ominus Y$ is a
certain 'complement' of $Y$ to $X$.

 It can be easily checked
that the small envelope of an additive category is additive; $X\to
(X,\id_X)$ gives  a full embedding $A\to A'$.

\begin{theo}\label{recw}

I Let $A$ be a full additive subcategory of some triangulated
$\cu$. Then the embedding $A\to\cu$ can be extended to a full
embedding of the small envelope of $A$ into $\cu$.

II Assume that $H$ is negative and generates $\cu$. Then the following statements are valid.

1. There exists a unique weight structure $w$ for $\cu$ such that
$H\subset \hw$. Moreover, it is  bounded.

2.  $\hw$ equals the Karoubi-closure of $H$ in $
\cu$. The latter is equivalent to the small
 envelope of $H$.

III  Suppose that $H$ is negative and weakly generates $\cu$, and  that
for any $X\in \obj\cu$ there exists a $j\in\z$ such that
\begin{equation}\label{bouab} \obj H\perp\{ X[i],\ i>j\}.\end{equation}
 Let $H'\subset
H$ be an additive subcategory. Suppose that one of the following conditions is fulfilled.

(i) There exists an infinite cardinality $c$ such that  any
coproduct of $<c$ objects of $H$ exists and belongs to $H$,
whereas $\card\ H'<c$. For any $X\in \obj\cu$ and  any $Y\in\obj H'$
the group $\cu(Y,X)$ considered as a $\cu(Y,Y)$-module can be
generated by $<c$ elements. Any object of $H$ can be presented as
$\coprod_{i\in I}C_i$ 
for $C_i\in \obj H'$, $\card\ I< c$. For any
$I$ such that $\card (I)<c$, any $Y\in \obj H'$, $j\in\z$, and $X_i\in \obj H,\
i\in I$, we have
\begin{equation}\label{addc} \cu (Y,\coprod_{j\in I} X_j)=\bigoplus
\cu (Y,X_j)\end{equation}

 or

(ii) Arbitrary coproducts exist in $H$; all objects of $H'$ are
compact; $\obj H'$ is a set; any object of $H$ can be presented as
$\coprod_{i\in I}C_i$ 
for $C_i\in \obj H'$ and some set $I$.

Then there exists a 
weight structure $w$ for $\cu$ such
that $\hw$ is the Karoubi-closure of $H$ in $\cu$. Moreover, $w$ is non-degenerate
and bounded  above. Lastly, in case (ii) $\cu^{w\le 0}$ is closed with respect to all those coproducts that exist in $\cu$;  in case (i) it 
is closed with respect to coproducts of $< c$ objects.

IV Suppose that all conditions    of assertion III ((i) or (ii)) except
(\ref{bouab})  are fulfilled. Denote the set of objects of $\cu$
satisfying (\ref{bouab}) for some $j\in \z$ by $\cu^-$; denote the
class of objects of $\cu$ satisfying (\ref{bouab}) for a fixed
$j\in \z$ by $\cu^{w\le j}$. Then
 the category $\cu^-$ is triangulated and satisfies all
conditions of assertion III (we will identify the class  $\cu^-$ with
the corresponding full subcategory of $\cu$).

V Assume that the conditions of assertion III(ii) are fulfilled, and that for any $Z\in \obj \cu$ the colimit $\inli_{i\ge 0} w_{\ge -i}Z$ (the connecting morphisms here are 'induced' by $\id_Z$ via Lemma \ref{fwd}(2)) exists. Then the following statements are fulfilled.

1. For any $Z\in \obj \cu$ we have $Z\cong \inli_{i\ge 0} w_{\ge -i}Z$.

2. A compact $X\in \obj \cu$ belongs to $\cu^{[j,q]}$ (for $j\le q\in \z$) if and only if
$H'\perp X[i]$ for all $i>q$, and  $X[i]\perp H'$ for all $i<j$.

\end{theo}
\begin{proof}

I We map $(X,p)$ to (any choice of) $\co(q)$; we denote this
object by $Z$.

Now we define the embedding on morphisms. We note that in $A$ the
map $q$ is a projection of $X$ onto $Y$. Hence in $A'$ we have
$X\cong (X,p)\bigoplus Y$, the isomorphism is given by $(p,q)$. Since
$q$ has a section in $\cu$, we have a distinguished triangle $Z\to
X\stackrel{q}{\to}Y\stackrel{0}{\to} Z[1]$, i.e., we also have a
similar decomposition of $X$ in $\cu$. It is easily seen that
$\cu(Z,Z')$ is given exactly by the formula (\ref{mthen}) if we
assume that $Z$ is a subobject of $X$, i.e., if we fix the splitting
of the projection $X\to Z$. Hence if we fix the embedding $Z\to X$
for each $(X,p)$ then (all possible choices) of objects $\co(q)$
would give a subcategory that is equivalent to the small envelope
of $A$; it is obviously additive.

Alternatively, the statement can be easily deduced from the
 functoriality of
 the idempotent completion procedure (proved in \cite{ba}; see
\S\ref{ridcomp} below).

II 1. We define $\cu^{w\ge 0'}$ as  the smallest
extension-stable (see Definition \ref{exstab})
 subclass of $\obj\cu$ that contains $\obj H[i]$ for $i\le 0$;
for $\cu^{w\le 0'}$ we take
a similar 'closure' of the set $\cup \obj H[i]$ for $i\ge 0$.

Obviously, $\cu^{w\ge 0'}$ and $\cu^{w\le 0'}$ satisfy
property (ii) of Definition \ref{dwstr}; we define $\cu^{w\ge
i'}$ and $\cu^{w\le i'}$ for $i\in\z$ in the usual way.

 If we have a distinguished triangle $X\to Y\to Z\to X[1]$ with
 $\cu(X,A)=\cu(Z,A)=\ns$ for some $X,Y,Z,A\in\obj \cu$, then
 $\cu(Y,A)=\ns$; the same statement is valid for a functor of the
 type $\cu(B,-)$.
  Hence from the fact that $ H[i]\perp H[j]$ for all $i<0\le j$,
we obtain (by induction) that
 $\cu^{w\ge 1'}\perp \cu^{w\le 0'}$.

Now we verify that any  $X\in\obj\cu$ has a '{}weight
decomposition{}' (with respect to $\cu^{w\ge 0'}$ and
$\cu^{w\le 0'}$). We prove this by induction on the 'complexity'
of $X$, i.e., on the number of distinguished triangles that we have to
consider to obtain $X$ from objects of $H[i],\ i\in\z$ (by considering cones of morphisms).

 For
$X$ of 'complexity' zero (i.e., for $X\in\obj H[i]$) we can take a
trivial weight decomposition, i.e., put $X^{w\le 0}$ equal to $X$ for
$i\ge 0$ and to $0$ otherwise; $X^{w\ge 0}$ will be $0$ and $X$,
respectively.

Suppose now that $X\cong\co(Y\stackrel{d}{\to} Z)$ for $Y,Z$ of
'complexity' less than that of $X$.   By the inductive assumption
there exist '{}weight decompositions{}' of $Y[1]$ and $Z$, i.e.,
distinguished triangles $Y[1]\stackrel{a}{\to}A\to B$ and
$Z\stackrel{a'[-1]}{\to}A'[-1]\to B'[-1]$ for $A \in
\cu^{w\le 0'}$, $A' \in \cu^{w\le -1'}$, $B \in \cu^{w\ge
1'}$, $B' \in \cu^{w\ge 0'}$. We apply Remark \ref{rtwd} for
$D=\cu^{w\le 0'}$ and $E=\cu^{w\le 0'}$. It yields a 'weight
decomposition' of $X$.

Now we take for $\cu^{w\ge 0}$ and $\cu^{w\le 0}$ the
Karoubi-closures of $\cu^{w\ge 0'}$ and $\cu^{w\le 0'}$,
respectively. By Lemma \ref{lsimple}(4), they satisfy the
orthogonality axiom of weight structures. Hence they define a weight
structure $w$ for $\cu$.

Now, since any object of $\cu$ can be obtained by a finite
sequence of considerations
 of cones of morphisms from objects of $\hw$, we obtain that
 $w$ is bounded.

It remains to check that $w$ is the only weight structure such
that $H\subset \hw$. By Proposition \ref{bw}(3), for any
weight structure $u$ satisfying $H\subset {\underline{Hu}}$ we have
$\cu^{w\ge 0'}\subset \cu^{u\ge 0}$ and $\cu^{w\le 0'}\subset
\cu^{u\le 0}$. Since $\cu^{u\ge 0}$ and $\cu^{u\le 0}$ are
Karoubi-closed in $\cu$, we also have
$\cu{}^{w\ge 0}\subset \cu^{u\ge 0}$
and $\cu{}^{w\le 0}\subset \cu^{u\le 0}$. Now Lemma \ref{lodn}
implies our claim immediately.


2. By assertion I, $\cu$ contains the small envelope of $H$. To
check that this envelope is actually contained in $\hw$ it suffices
to note that the object $X\ominus Y$ can be presented both as a
cone of the 'embedding' $Y\to X$ and of the projection $X\to Y$.


Now we verify the inverse inclusion; let $X\in \cu^{w=0}$.

 We apply the weight complex functor $t$. 
 We
obtain that $t(X)=X$   is a $K_\w^{b}(\hw)$-retract both of a certain $A\in  C^{b, \ge 0}(H) $ and of some $B\in 
C^{b, \le 0}(H) $ (here we use the descriptions of $\cu^{w\le 0}$ and $\cu^{w\ge 0}$ given above, and the weak exactness of $t$). Next, applying Lemma \ref{lzm} we obtain
that the same is true in $K^b(\hw)$. Since $t$ is also conservative, we can replace $\cu$ by $K^b(\hw)$ (with the same $H$); the assertion in this case can be verified 'by hand' (cf. the proof of Proposition 6.4.1 of \cite{mymot}).

III  Again, for $\cu^{w\ge 0}$ we take the smallest Karoubi-closed
extension-stable
subset of $\obj\cu$ that  contains $H[i]$ for $i\le 0$.

We take $\cu^{w\le 0}= (\cup_{i<0} \obj H[i])^{\perp}.$ 
Since any object of $H$ is a coproduct of objects of $H'$, we could have replaced $H$ by $H'$ in this description; it follows immediately  that $\cu^{w\le 0}$ is closed with respect to the corresponding coproducts.

We verify the orthogonality axiom using  the same argument as in the proof of assertion II1. We have
$Y\perp \cu^{w\le -1}$
for  any $Y \in H[i]$ for $i\le 0$.  Now 
 using the fact that
all (co)representable
functors are homological on $\cu$, we obtain that the same is true
for any $Y\in \cu^{w\ge
0}$.

 $\cu^{w\le 0}$ is  Karoubi-closed and extension-stable by  Lemma \ref{lsimple}(1);
$\cu^{w\ge 0}$ also fulfills these properties. 

To prove that $w$ is a weight structure, it remains to prove the
existence of weight decompositions (see  Lemma
\ref{lsimple}(2)). We will construct $X^{w\le 0}$ and $X^{w\ge 1}$ for
a fixed $X\in\obj\cu$ explicitly. The construction may be called
a {\it weight resolution}, cf. the proof of Proposition
\ref{madts} below and Proposition 7.1.2 of \cite{axstab}.

First we treat case (i).  For each object $Y$ of $H'$ any
$Z\in\obj\cu$ we choose a set of $f_i(Y,Z)\in \cu (Y,Z)$ of
cardinality $<c$ so that $f_i(Y,Z)$ are $\cu (Y,Y)$-generators of
$\cu (Y,Z)$. Assume that  (\ref{bouab}) if fulfilled for some $j\in\z$.


Now we construct a certain sequence of $X_k$ for $k\le j$ starting
from $X_j=X$. For $k=j$ we take $P_j=\coprod_{Y\in \obj
H',f_i(Y,X_j[j])} Y$. 
Note that the number of summands is $<c$,
hence the sum exists and belongs to $\obj H$. Then we have a
morphism $f_j:P_j\to X_j[j]$ given by $\prod f_i(Y,X[j])$. Let
$X_{j-1}[j]$ denote a cone of $f_j$. Repeating the construction
for $X_{j-1}$ instead of $X_j$ and with $k=j-1$ we get an object
$P_{j-1}\in \obj H'$, $f_{j-1}:P_{j-1}\to X_{j-1}[j-1]$; we denote
a cone of $f_1$ by $X_{j-2}[j-1]$. Proceeding,  we get an infinite
sequence of $(P_i,f_i,X_i)$. Note that we have $P_i\in \cu^{w\ge
0}$.

We denote the maps $X_i\to X_{i-1}$ given by the construction by
$g_{i}$,  $h_i= g_{j}\circ \dots \circ g_{i+2}\circ g_{i+1}:X\to
X_i$. We denote a cone of $h_i$ by $Y_i[-1]$; the map $Y\to X_i$
given by the corresponding distinguished triangle by $r_i$.
 Then Remark \ref{rpoto}(3) yields that $Y_i[i]\in \cu^{w\ge 0}$
 for all $i\le j$
(see the definition of $\cu^{w\ge 0}$).

Now we  denote  $Y_0$ by $Y$ and $X_0$ by $Z$.  $Y,Z$ will be our
candidates for $X^{w\ge 0}$ and $X^{w\le 0}$.

It remains to prove that $Z\in \cu^{w\le 0}$. We should check that
$\cu(C,Z[k])=\ns$ for all $k>0$, $C\in \obj H$.  Since $\cu(-,Z)$
transforms arbitrary coproducts into products, it suffices to consider
$C\in\obj H'$.

First we prove that $C\perp X_{k-1}[k]$ for all $ k\le j$.

 We use the distinguished triangle
\begin{equation}\label{dist}
V_k: P_k\to X_{k}[k]\to X_{k-1}[k]\to P_k[1].\end{equation} Using
(\ref{addc}), we obtain $\cu(C,P_k[k])= \bigoplus_{Y\in\obj H', f_i(Y,
X_{k}[k])}\cu (C,Y).$ By the definition of $f_i(C,Y)$ we obtain that
this  group surjects onto $\cu(C,X_{k}[k])$. Moreover,
$\cu(C,P_k[1]=\bigoplus_{Y\in\obj H', f_i(Y,
X_{k}[k])}\cu(C,Y[1])=\ns$. We obtain  $\cu(C,X_{k-1}[k])=\ns$.

Now we use distinguished triangles $V_l$  for all $l< k$. Again
(\ref{addc})  yields $\cu (C, P_l[1])= \cu (C, P_l[2])=\ns$. Hence
$\cu(C, X_{l-1}[k])= \cu(C,X_{l}[k])=\ns$ for all $l< k$.

Hence $C\perp Z[k]$ for all $j\le k>0$.

Lastly, the distinguished triangles $V_k$  easily yield by
induction that $C\perp X_l[k]$ for all $l\le j$ and $k>j$.

The proof in case (ii) is almost the same; one should only always
replace some choice of generators $f_i(Y,Z)\in\cu(Y,Z)$ by all
elements of $\cu(Y,Z)$.

$(\cu,w)$ is obviously bounded  above by (\ref{bouab}).

Now we check that $(\cu,w)$ is non-degenerate. The condition
(\ref{bouab})  implies that $\cap \cu^{w\ge i}=\ns$. Next, for any
$X\in\obj\cu\setminus\ns$ there exists an $f\in \cu(Y[i],X)$ for
some $Y\in H$ and $i\in\z$ such that $f\neq 0$. Hence such an $X$
does not belong to $\cu^{w\le -1-i}$ (see the definition of
$\cu^{w\le i}$ in the proof of assertion III).

It remains to consider $\hw$. It contains $H$ by the definition of $w$. The method used in the proof of assertion II2 also easily yields that any object of $\hw$ is a retract of an object of $H$.

IV Everything is obvious except that a cone of a morphism of
objects of $\cu^-$ belongs to $\obj\cu^-$. This fact is easy also
since the functors $\cu(Y,-)$ are homological.

V1. By Lemma \ref{coulim}(1), there exists a morphism $c:\inli_{i> 0} w_{\ge -i}Z\to Z$ compatible with the corresponding weight decompositions.  By part 3 of loc. cit., $c$ yields an isomorphism $\cu(T[r],\inli_{i> 0} w_{\ge -i}Z)\to \cu(T[r],Z)$ for all $T\in \obj H'$, $r\in \z$ (note that $T[r]$ is compact, whereas $\cu(T[r],w_{\ge l}Z)\cong \cu(T[r],w_{\ge l-1}Z)\cong \cu(T[r],Z)$ for any $l<-r$). Since $H$ weakly generates $\cu$, the same is true for $H'$; hence we obtain that $c$ is an isomorphism.
 
 2. If $X\in \cu^{[j,q]}$, then the orthogonality axiom of weight structures immediately yields that is satisfies the orthogonality conditions required.

We prove the converse implication. 
 
 Since $\cu(-,X)$ is additive and converts (small) coproducts into products, we also have $\cu^{w=i}\perp X$ for all $i>q$. 
 It follows that $\cu^{[q+1,r]}\perp X$ for any $r>q$ (since any $Y\in \cu^{[q+1,r]}$ possesses a weight Postnikov tower with $Y^i=0$ for $i\le q$ and $i>r$).
 Since $\cu$ is bounded above, 
 we also obtain  that $\cu^{w\ge q+1}\perp Y$. 
 Then 
Proposition \ref{bw}(1)  yields that $X\in \cu^{w\le q}$.

It remains to verify that $X\in \cu^{w\ge j}$. 
By Proposition \ref{bw}(2), to this end
it suffices to prove that $X\perp Z$ for any $Z\in \cu^{w\le j-1}$. We fix some $Z$ 

Since $X$ is compact, we have 
$X\perp \cu^{w=i}$ for all $i<j$. Considering weight Postnikov towers  again, we obtain that
 $X\perp \cu^{[l,j-1]}$ for any $l<j$. In particular, $X\perp w_{\ge l}Z$ 
 (here we have $w_{\ge l}Z\in \cu^{[l,j-1]}$ by Proposition \ref{mwdec}(1)). 

 Since $X$ is compact, Lemma \ref{coulim}(3) yields that  $X\perp \inli_{i> 0} w_{\ge -i}Z$. 
  Hence the previous assertion yields  the result.

\end{proof}

\begin{coro}\label{wsh}
It is well known  that there are no morphisms of positive degree
between (copies) of the sphere spectrum $S^0$ in the stable
homotopy category $SH$; cf. \S\ref{sh} below. 

Hence 
Theorem \ref{recw}(II) immediately implies that the category of finite
spectra $SH_{fin}$ (i.e., the full subcategory of $SH$ generated by
$S^0$) has a bounded weight structure $w$.
 Its heart can be
described as a category $H$ of finite sums of (copies of) $S^0$
(since any retract of $(S^0)^n$, $n>0$, is trivial, no new objects
appear in
the small envelope of $H$). Since $SH(S^0,S^0)=\z$, $\hw$ is
equivalent to $\ab_{fin.fr}$ (the category of finitely generated
free abelian groups).

This weight structure obtained is {\it adjacent} to the
Postnikov $t$-structure for $SH$; see Definition \ref{deadj}  and \S
\ref{sh} below.

\end{coro}

\begin{rema}\label{exoticcond}

1. Recall that in the proof of Theorem \ref{recw} it was specified
explicitly how to recover $w$ from $\cu^{w=0}$.

2. The conditions of parts III and IV of the theorem may seem to be rather
exotic. Yet they can be easily verified for a certain subcategory of
{\it quasi-finite} objects in $SH$, see \S\ref{sh}.


3. Obviously, if for any $X\in \obj\cu$ and  any $Y\in\obj H'$ the
group $\cu(Y,X)$  is generated by  $<c$ elements as a group, then
it is also generated  by $<c$ elements as a $\cu(Y,Y)$-module. In
particular, this is the condition which we will actually check for
the category $SH_{fin}$.

4. $\cu^{w\le 0}$ described in the proof of  Theorem
\ref{recw}(III) is often a $\cu^{t\le 0}$-part of a certain
$t$-structure; then this $t$-structure is {\it left adjacent} to $w$
(see Definition \ref{deadj} below). Yet in order for the
$t$-decompositions to exist when we take  the only possible
candidate for $\cu^{t\ge 0}$ (cf. Proposition \ref{lads} below) the
homotopy  colimit of all $Y_i$
(defined as in the proof of part III) 
should
exist for all $X\in\obj \cu$ (cf. the proof of Proposition
\ref{madts}). Note that this is not true for the category
$SH_{fin}$ (see
 Corollary \ref{wsh}). For example, one can note that
Eilenberg-MacLane spectra do not belong to $SH_{fin}$. Another
example is the following one:  the Chow weight structure is defined on $\dmge$,
whereas the corresponding $t$-structure is only defined on $\dme$;
see \S\ref{remmymot}  and \S\ref{ddme}.

This shows that weight structures exist 'more often' than
$t$-structures, whereas Theorem \ref{sdt} below shows that they
contain (more or less) the same information as the corresponding
$t$-structures. This evidence  supports the author's belief that
weight structures are more relevant for 'general' triangulated
categories than $t$-structures. See also (part 2 of) Remark
\ref{rtrfun} for more facts supporting this opinion.

Besides, it can be easily seen that the natural 'almost dual' to
the statement of part II (i.e., we take a positive generating
subcategory $H$ and ask whether a $t$-structure with $H\subset \hrt$
exists) is false.

\end{rema}

\subsection{Adjacent weight and $t$-structures;
 weight-exact functors}\label{sadst}

\begin{defi}\label{deadj} We say that a weight structure $w$ is
left (resp. right)
{\it adjacent} to a $t$-structure $t$ if $\cu^{w\le 0}=\cu^{t\le 0}$
(resp. $\cu^{w\ge 0}=\cu^{t\ge 0}$).

In this situation we will also say that $t$ is 
adjacent to $w$.
\end{defi}

A simple example is given by $\cu$ being  (an appropriate version of) $D(A)$
 for an abelian $A$. Then (under certain conditions ensuring
that $D^?(A)$ is isomorphic to $K^?(Proj\,A)$)
the stupid weight structure coming from $K^?(Proj\,A)$
is left adjacent to the canonical $t$-structure for $D^?(A)$.
Dually, if $D^?(A)\cong K^?(Inj\,A)$, then the canonical $t$-structure
for it is right adjacent to the stupid weight structure coming from
$K^?(Inj\,A)$. We will construct more (interesting) examples of adjacent
structures below.

The following theorem describes several properties of
adjacent structures. Parts 1,2 of the theorem were already proved (above);
parts 7,8 of it
show that for $Y\in\obj \cu$ the virtual $t$-truncations of the functor
$\cu(-,Y)$ (defined in Remark \ref{rtrfun}) are represented
by actual $t$-truncations of $Y$ with respect to $t$.
We also prove that $\hrt$ and $\hw$
are connected by a natural generalization of the relation between
the categories $A$ and $Proj\, A$ (for an abelian $A$). Still note that
in the latter case we also have $\hw\subset \hrt$; 
this is a rather
non-typical situation.

\begin{theo}[Duality theorem]\label{sdt}

Assume that $w$ is left adjacent to $t$. Then the following statements are
valid.

1. $\cu^{w\ge 0}=
{}^\perp\cu^{w\le -1}$; $\cu^{w\le 0}=(\cu^{w\ge 1})^{\perp}$.

2. $\cu^{t\ge 0}=(\cu^{t\le - 1})^{\perp}$; $\cu^{t\le 0}
={}^\perp\cu^{t\ge 1}$.

3. The functor $\cu(-,\hrt): \hw \to \hrt^*$ (see the Notation) that
sends $X\in \cu^{w=0}$ to $Y\to \cu(X,Y),\ (Y\in \cu^{t=0}),$ is a
full embedding of $\hw$ into the full subcategory
$Ex(\hrt,\ab)\subset \hrt^*$ which consists of exact functors.

4. The functor $\cu(\hw,-): \hrt \to \hw_*$ that sends $X\in
\cu^{t=0}$ to $Y\mapsto \cu(Y,X),\ (Y\in \cu^{w=0}),$ is a full exact
embedding of $\hrt$ into the abelian category $\hw_*$.

5. Suppose that $t$ is non-degenerate. Then $\cu^{t=0}$ equals 
$S=\{X\in \obj \cu:\ \cu^{w=0}\perp X[i]\ \forall \ i\neq 0
\}.$

6.  Let $i\in\z$; let $Y\in\obj \cu$ be fixed. Then  for the functor
$F(X)=\cu(X,Y)$ we have $W^i(F)(X)=\imm (\cu(X,\tau_{\le i}Y)\to 
\cu(X,Y))$ for any $X\in\obj\cu$.

7. For any $i,j$ we have a functorial isomorphism $$\cu(X,Y^{t\le
i}[j])\cong \imm(\cu(X^{w\le -j},Y[i])\to \cu(X^{w\le
1-j},Y[i+1])).$$  

8. For any $i,j$ we have a functorial isomorphism $$\cu(X,Y^{t\ge
i}[j])\cong \imm(\cu(X^{w\ge -j},Y[i])\to \cu(X^{w\ge
-1-j},Y[i-1])).$$

9. For any $X,Y\in \obj \cu$  let $\dots\to Y^{-1}\to Y^0\to
Y^1\to\dots$ denote an arbitrary choice of the weight complex for
$Y$ (in its homotopy equivalence class). Then we have
\begin{equation}\label{nfrestrw} \cu(Y,X^{t=0})
=(\ke ( \cu(Y^{0},X)\to
 \cu(Y^{-1},X)) /\imm
(\cu(Y^1,X)\to \cu(Y^0,X)).\end{equation} 

10. We have $\cu^{b,w\ge 0}=\obj \cu^b
\cap ({}^\perp \cup_{i<0} \cu^{t=i})$ and
$\cu^{b,w\le 0}=\obj \cu^b \cap ({}^\perp \cup_{i>0} \cu^{t=i})$.

\end{theo}
\begin{proof}

1. This is just
 Proposition \ref{bw}(1,2). 

2. A well-known property of $t$-structures (certainly, it doesn't
depend on $w$).

3. First we note that for any $X\in \cu^{w=0}$ orthogonality for $w$
implies $X\perp \cu^{w\le -1}(=\cu^{t\le -1})$,
whereas orthogonality for $t$ gives $X\perp \cu^{t\ge 1}$. In particular,
\begin{equation}\label{orth}
X\perp (\cup_{i\neq 0} \cu^{t=i}).
\end{equation}

Now, short exact sequences in $\hrt$ give distinguished triangles in
$\cu$. Hence for any homological functor $F:\cu\to \ab$ and for
$0\to A\to B\to C\to 0$ being a short exact sequence in $\hrt$
we have a long
exact sequence $\dots \to F(C[-1])\to F(A)\to F(B)\to F(C)\to
F(A[1])\to\dots$. If $F=\cu(X,-)$ then $F(C[-1])=F(A[1])=\ns$ (as
was just noted). Hence objects of $\hw$ induce exact functors on
$\hrt$.

To prove that the restriction $\hw \to \hrt^*$ is a fully faithful
functor it suffices to prove the following:  the restriction of the functor
$\cu(X,-)$ to $\hrt$ for $X\in \hw$ determines $X$ in a functorial way.
Using Yoneda's lemma, we see that it suffices to recover $\cu(-,X)$
from its restriction.

We prove that
\begin{equation}\label{frestrt} \cu(X,Y)\cong \cu(X,H^t_{0}(Y))\
\forall X\in \cu^{w=0},\ Y\in\obj\cu.\end{equation}

We apply $t$-decompositions  (i.e., (\ref{tdec})) twice.

First we obtain a distinguished triangle
$Y^{t\ge 1}[-1]\to Y^{t\le 0}\to
Y\to Y^{t\ge 1}$. Since $\cu(X,Y^{t\ge 1}[-1])=\cu(X,Y^{t\ge 1})$,
we obtain $\cu(X,Y)=\cu(X,Y^{t\le 0})$.

Next, we have a distinguished triangle $Y^{t\le -1} \to Y^{t\le
0}\to H^t_{0}(Y) \to Y^{t\le -1}[1].$ Since $\cu(X,Y^{t\le -1})=
\cu(X,Y^{t\le -1}[1])=\ns$, we obtain (\ref{frestrt}).

4. Again, it suffices to prove that the restriction of the functor
$\cu(-,X)$ to $\hw$  for $X\in \cu^{t=0}$  determines $X$
functorially.

We note that for any $X\in \cu^{t=0}(\subset \cu^{w\le 0})$
the orthogonality axiom for
$t$ implies $(\cu^{w\le -1})=\cu^{t\le
-1}\perp X$, whereas orthogonality for $w$ gives $\cu^{w\ge 1}\perp X$.

 Now we prove that
\begin{equation}\label{frestrw} \cu(Y,X)\cong (\ke ( \cu(Y^{0},X)
\to \cu(Y^{-1},X)) /\imm
(\cu(Y^1,X)\to \cu(Y^0,X)).\end{equation} 

Indeed, consider the (infinite) weight decomposition of $Y$ that
gives our choice of the weight complex and apply Theorem
\ref{csps} to the functor $\cu(-,X)$. The spectral sequence
obtained converges since it satisfies condition II(ii) of Theorem
\ref{csps} (it has only one non-zero column!).
 It remains to
note that this (possibly) non-zero column 
of $(E_1^{pq}(T(\cu(-,X)(Y)))=(\cu(Y^{-p},X[q]))$
 is exactly $(E_1^{p0}(T))=\dots\to \cu(Y^1,X)\to
\cu(Y^0,X) \to \cu(Y^{-1},X)\to\dots$. 

We obtain (\ref{frestrw}).

5. By assertion 4, an object of $\hrt$ is non-zero if and only if
it represents a
non-zero functor on $\hw$. Hence applying (\ref{frestrt}) we obtain
that $S$ is exactly the class of objects that satisfy $H_{i}^t(X)=0$
for all $i\neq 0$. It remains to note that for a non-degenerate $t$
this class is exactly $\cu^{t=0}$.

6. We can assume that $i=0$. We should check that $g\in \cu(X,Y)$
lifts to some $h\in \cu(w_{\le 0}  X,Y)$ if and only if it lifts to some $l\in
\cu(X,\tau_{\le 0}Y)$. Now, the equality $$\cu(w_{\le 0}
X,\tau_{\ge 1} Y)=\cu(w_{\ge 1}X,\tau_{\le 0} Y) =\ns$$
 yields that any morphism of any  of the two morphism
groups in question can be
lifted to some $m\in \cu(w_{\le 0} X,\tau_{\le 0}Y)$. Hence if one
of $(h,l)$ exists, then the other one can be constructed from the
corresponding $m$ 
using the commutativity of the diagram
$$\begin{CD}
\cu(w_{\le 0} X,\tau_{\le 0}Y)@>{}>>\cu(w_{\le 0}X,Y)\\
@VV{}V@VV{}V \\
\cu(X,\tau_{\le 0}Y)@>{}>>\cu(X,Y)
\end{CD}$$

7. Shifting $X,Y$ we can easily reduce the statement to the case
$i=j=0$.

The $t$-decomposition of $Y$ yields  exact sequences
$\ns=\cu(X^{w\le 0},Y^{t\ge 1}[-2]) \to  \cu(X^{w\le 0},Y^{t\le
0})\to \cu(X^{w\le 0},Y) \to \cu(X^{w\le 0},Y^{t\ge 1}[-1])
 =\ns$ and
$$\ns=\cu(X^{w\le 1},Y^{t\ge 1}[-1]) \to  \cu(X^{w\le 1},Y^{t\le
0}[1])\to \cu(X^{w\le 1},Y[1]) \to \cu(X^{w\le 1},Y^{t\ge 1})\to
\dots$$

Next, weight decompositions of $X$ and $X[1]$  similarly yield that
the obvious homomorphism $\cu(w_{\le 0}X, Y^{t\le 0})\to \cu(w_{\le
1}X,Y^{t\le 0})$ is surjective whereas $\cu(w_{\le 1} X,Y^{t\le
0})\cong \cu(X,Y^{t\le 0})$.

We obtain a commutative diagram
$$\begin{CD}
\cu(w_{\le 0} X,Y^{t\le 0})@>{f}>>\cu(w_{\le
0}X,Y)\\
@VV{g}V@VV{h}V \\
\cu(w_{\le 1} X,Y^{t\le 0})@>{p}>>\cu(w_{\le 1}X,Y)
\end{CD}$$
with $f$ being bijective, $g$ being surjective, and $p$ being
injective. Hence $\cu(X,Y^{t\le 0})\cong \cu(w_{\le 1} X,Y^{t\le
0})\cong \imm g\cong \imm h$.

Note that the isomorphism constructed  is obviously natural in $Y$
whereas it is natural in $X$ by  \ref{fwd}(2). 

8. This assertion is exactly the dual of the previous one (see
Remark \ref{rdd}).

9. Immediate from (\ref{frestrw}) and (\ref{frestrt}).

10. Let $X$ be an object of $\cu^b$.

If we also have $X\in \cu^{w\ge 0}$ (resp. $X\in \cu^{w\le 0}$) the
orthogonality statements desired are valid by  assertion 1.

Now we prove the converse implication. Assume that
 $X\perp \cup_{i<0} \cu^{t=i}$. We should check that
$X\perp Z$ for all
$Z\in \cu^{t\le -1}$. We have $X\perp (H^t_{i}(Z)[-i])$ for all $i$.
Besides, since $X$ is bounded, we have $X\perp (\tau_{\le j} Z)$
for some $j$ (that is small enough). Hence considering the
$t$-decompositions of $Z^{t\le k}[-1]$ for all $k>j$ one can easily
obtain the orthogonality statement required.

The case $ X \perp \cup_{i>0} \cu^{t=i}  $ is
considered similarly.


\end{proof}

\begin{rema}\label{rsdt}

1.  One can often describe the images of embeddings in assertions 3
and 4.

2. Dually to assertion 6, for $w$ right adjacent to $t$ and
$F(X)=\cu(Y,X)$ we obtain $W_i(F)(X)=\imm(\cu(\tau_{\ge i} Y,X)\to
\cu(Y,X))$. 

3. Assertion 6 and its dual show that weight truncations are
'almost adjoint' to the corresponding $t$-truncations.
These statements are counterparts to the fact that
(for an arbitrary $t$-structure)
 any morphism $X\to Y$ for $Y\in\cu^{t\ge 0}$ can be uniquely
factored through $X^{t\ge 0}$ (and to its dual).

4. Assertions 7 and 9 (almost) imply that the 
derived exact couple  for the spectral sequence
$\cu(X^{-p}[-q],Y)\implies \cu(X[-p-q],Y)$ (as in Theorem \ref{csps})
can also be described in terms of $\cu(X[i],Y^{t\le j})$ and
$\cu(X[i],Y^{t= j})$; see \S2.6 of \cite{bger} for the
complete proof of this fact. 
This is
no surprise by  Remark \ref{rtrfun}(3,4).
It follows that the spectral
sequence $S$ converging  to $\cu(X,Y)$ corresponding to the
$t$-truncations of $Y$ can be 'embedded' into our $T$ (i.e., for
all $i>0$ any $E_i^{pq}(S)\cong
E_{i+1}^{p'q'}(T)$ for  $p'=q+2p,q'=-p$; 
these isomorphisms
respect the structure of spectral sequences).

In algebraic topology, this statement corresponds to the fact (and
implies it) that the Atiyah-Hirzebruch spectral sequence for the
cohomology of a space $X$ with coefficients in a spectrum $S$ can
be obtained either by considering the cellular filtration of $X$ or
the Postnikov $t$-truncations of $S$.

Note that  our method describes all terms of the corresponding exact couples (in
contrast to \cite{paran}, for example). The advantage of this is
that the $D$-terms of those may be very interesting; see Proposition
\ref{trcoh} and Corollary \ref{cmot} below.

5. In fact, one can extend the notion of  adjacent structures to
the case when there are two distinct triangulated categories $\cu$
with a weight structure $w$ and $\du$ with a $t$-structure $t$. They
should equipped with a {\it duality} $\Phi:\cu^{op}\times\du\to A$,
$A$ is an abelian category (that generalizes
$\cu(-,-):\cu^{op}\times\cu \to \ab$) along with a (bi)natural
transformation $\Phi(X,Y)\cong \Phi (X[1],Y[1])$; $\Phi$ is
homological with respect to both arguments (see Definition 
2.5.1 of \cite{bger}).
 $\du^{t \ge 1}$ should
 annihilate  $\cu^{w\le 0}$ with respect to $\Phi$, whereas
$\du^{t \le -1}$ should  annihilate  $\cu^{w\ge 0}$; cf. (part 1 of)
Remark \ref{rtrfun}. 
In   Theorem 2.6.1 of \cite{bger} the author proves
natural analogues of  Theorem \ref{sdt}(6--9). 
This 
yields another
proof of the comparison of spectral sequences statement
of \S\ref{twmot} 
 (in the general case; see Remark \ref{tmotco}); cf.
also (parts 3 and 4 of) Remark \ref{rtrfun}.

See also \S\ref{scompt} below for further ideas in this direction.

 6. Even more generally, for parts 7,8 of the theorem it suffices
for $X,Y$ (lying either in the same category or in different ones)
to have Postnikov towers whose terms satisfy   the same orthogonality
conditions as the ones provided by the definition of the
weight structure. 
Indeed, the same proofs work!

Furthermore, it is sufficient for the orthogonality conditions to be
satisfied 'in the limit' for a directed system of Postnikov towers
(for $X$). 
Again, it is no problem to generalize the proof to this
case. Still, in order to make the statement easier to understand,
the author chose to formulate it in \S\ref{coniv} below only for a
partial (yet very important!) case corresponding to  coniveau
spectral sequences.

\end{rema}


Recall now that if an exact functor $\cu\to \cu'$ is $t$-exact with
respect to some $t$-structures on these categories, its
(left or right) adjoint is usually not $t$-exact (it is only left
or right $t$-exact, respectively). This situation
can be described much more precisely if there exist
 adjacent weight structures for these $t$-structures.

Till the end of this subsection, $\cu$ and $\cu'$ will
be triangulated categories, that are possibly endowed with
$t$-structures $t$ and $ t'$, and with weight structures $w$ and
 $w'$, respectively. $F:\cu\to \cu'$ and $G:\cu'\to \cu$ will
be exact functors.

\begin{defi}\label{dadad}

1. $F$ is called $t$-exact (with respect to $t,t'$) if it maps
$\cu^{t\le 0}$ to $\cu'^{t'\le 0}$, and
maps $\cu^{t\ge 0}$ to $\cu'^{t'\ge 0}$.

2.
 $F$ is called weight-exact (with respect to $w,w'$) if it maps
$\cu^{w\le 0}$ to $\cu'^{w'\le 0}$, and
maps $\cu^{w\ge 0}$ to $\cu'^{w'\ge 0}$.
\end{defi}

Certainly, $t$- and weight-exact functors also
map hearts of the corresponding structures to hearts.

Now we describe the main properties of the notions introduced.

\begin{pr}\label{padad}

1.  Any composition of $t$-exact functors is $t$-exact.

2. Any composition of weight-exact functors is weight-exact.

3. Suppose that $(w,t)$ and $(w',t')$ are left adjacent and that $G$ is the left
adjoint to $F$. Then $G$ is weight-exact if and only if $F$ is $t$-exact.

4. Suppose that $(w,t)$ and $(w',t')$ are right adjacent and  that $G$ is right
adjoint to $F$. Then $G$ is weight-exact if and only if $F$ is $t$-exact.

\end{pr}
\begin{proof}

1. Obvious from the definition of $t$-exact functors, and also
well-known.

2. Obvious from the definition of weight-exact functors.

3. We only need the orthogonality properties of weight and $t$-structures here.

  By  Theorem \ref{sdt}(1), $G(\cu'^{w'\ge 0})\subset
 \cu^{w\ge 0}$  if and only if
 $G(\cu'^{w'\ge 0})\perp \cu^{w\le -1} $. The adjunction (along with shifting by $[-1]$) yields
the following:
this is
equivalent to $\cu'^{w'\ge 1}\perp F(\cu^{w\le 0}).$
Applying  Theorem \ref{sdt}(1) again, we obtain that this is
equivalent to $F(\cu^{w\le 0})\subset \cu'^{w'\le 0}$. Now
 recall that  $\cu^{w\le 0}=\cu^{t\le 0}$ and
$\cu'^{w'\le 0}=\cu'^{t'\le 0}$. Hence we proved that
$G(\cu'^{w'\ge 0})\subset
 \cu^{w\ge 0}$  if and only if $F(\cu^{t\le 0})\subset
 \cu'^{t'\le 0}$.

It remains to check that  $G(\cu'^{w'\le 0})\subset
 \cu^{w\le 0}$  if and only if $F(\cu^{t\ge 0})\subset
 \cu'^{t'\ge 0}$. Essentially, this is a well-known
property of $t$-structures (that we already mentioned above).
The first inclusion formula can be rewritten as follows:
$G(\cu'^{t'\le 0})\subset
 \cu^{t\le 0}$. Next, Theorem \ref{sdt}(2)
yields that this happens if and only if
$G(\cu'^{t'\le 0})\perp\cu^{t\ge 1}$.
Applying the adjunction and shift by $[1]$ we conclude the proof.

4. This is exactly the dual of the previous assertion
(see Remark \ref{rdd}).

\end{proof}

\begin{rema}
Certainly, one can also introduce the notions of right weight-exact and
left weight-exact functors, and prove the natural versions of the
assertions above for them. In fact, the arguments needed for this
 are essentially
contained in the proof of Proposition \ref{padad}.

\end{rema}

Theorem \ref{sdt} also yields a simple description of adjacent
 structures (of any type) when they exist.

\begin{pr}\label{lads}

1. Let $w$ be a weight structure for $\cu$. Then there exists a
$t$-structure which is left (resp. right) adjacent to $w$ if and only
 if for
$\cu^{t\le 0}=\cu^{w\le 0}$ and $\cu^{t\ge 0}=(\cu^{w\le -1})^\perp$
(resp. $\cu^{t\ge
0}=\cu^{w\ge 0}$ and $\cu^{t\le 0}={}^\perp \cu^{w\ge 1}$), and any
$X\in \obj \cu$ there exists a
$t$-decomposition (\ref{tdec}) of $X$. In this case our choice of
($\cu^{t\le 0}$, $\cu^{t\ge 0}$) is the only one possible.

2. Let $t$ be a $t$-structure for $\cu$. Then there exists a weight
structure which is left (resp. right) adjacent to $w$ if and
only if for
$\cu^{w\le 0}=\cu^{t\le 0}$ and $\cu^{w\ge 0}={}^\perp \cu^{t\le -1}$
(resp. $\cu^{w\ge
0}=\cu^{t\ge 0}$ and $\cu^{w\le 0}=(\cu^{t\ge 1})^\perp$), and any
 $X\in \obj \cu$ there exists a
weight decomposition (\ref{wd}) of $X$. In this case our choice of
 ($\cu^{w\le 0}$, $\cu^{w\ge 0}$) is the only one possible.
\end{pr}

\begin{proof}

First we note that by  Theorem \ref{sdt}(1,2) our
choices of the
structures are the only ones possible. Hence it suffices to check
when these choices indeed give the corresponding structures.

1. We only consider the left adjacent structure case; the 'right
case' is similar (and, in fact, dual; see Remark \ref{rdd}).

The class
$\cu^{t\ge 0}$ is automatically strict;
since $\cu^{w\le
0}[1]\subset \cu^{w\le 0}$, we have
$\cu^{t\ge 0}\subset \cu^{t\ge 0}[1]$.

Hence we obtain a $t$-structure if and only if there always exist
   $t$-decompositions.

2. As for assertion 1, we consider only the 'left' case (for the same
reason).

It is well known that $\cu^{t\le 0}$ is Karoubi-closed; hence both
$\cu^{w\le 0}$ and $\cu^{w\ge 0}$ are Karoubi-closed also. Again
 $\cu^{t\le
0}[1]\subset \cu^{t\le 0}$ implies
$\cu^{t\ge 0}\subset \cu^{t\ge 0}[1]$.

Hence we obtain a weight structure if and only if
 there always exist   weight decompositions.
\end{proof}

\subsection{Existence of adjacent structures}\label{eads}

Now   we study certain sufficient conditions for adjacent weight
and $t$-structures to exist.
We will  prove a
 statement that is  relevant for Voevodsky's $\dme$ and for $SH$.

First we describe a certain version of the compactly generated category
notion; $\dme$ and $SH$ will satisfy our conditions.

\begin{defi}\label{wg}

 We will say that  a set of objects $C_i\in \obj\cu$, $i\in
I$ ($I$ is a set) {\it negatively well-generates} $\cu$
  if

(i) $C_i$ are compact; they weakly generate $\cu$ (cf. the
Notation).

(ii) For all $j>0$ we have
$\{C_i\}\perp \{C_i[j]\}$ (i.e., the set $\{C_i\}$ is negative).

(iii) $\cu$ contains the category $H$ whose objects are arbitrary
(small) coproducts of $C_i$; $\cu$   also contains  all homotopy
colimits of $X_i\in\obj \cu$ (see Definition \ref{dcoulim}) such that
$X_{-1}=0$ and  $\co(X_i\to X_{i+1})\in \obj H[i]$.

\end{defi}

\begin{theo}\label{madts}

I1. Suppose that $C_i\in \obj\cu$, $i\in I$,  {\it negatively
well-generate} $\cu$. For $H$  described in (iii) of Definition
\ref{wg} we consider a full subcategory $\cu^-\subset \cu$ whose
objects are
\begin{equation}\label{bouab1} X\in \obj \cu:\
\forall\ Y\in\obj H\ {\text there\
exists}\ j\in\z\ {\text such\ that}\ Y\perp X[i]\ \forall
i>j.\end{equation}

Then there exist 
  a weight structure $w$ on $\cu^-$ and a $t$-structure $t$ on
$\cu$ such that $H\subset \hw$, $t$ restricts to a $t$-structure on
$\cu^-$, and $\cu^{t\le 0}=\cu^{-,w\le 0}$.

(Note that $w$ and $t$ restricted to $\cu^-$ are adjacent by
definition.)

2. If $\cu$ also admits arbitrary countable coproducts,  then $w$
can be extended to the whole $\cu$.

 II Let $C_i,\cu, w,t$ be either as in assertion I2 or as in I1 with
the additional
 condition $\cu=\cu^-$ (i.e., $w$ is defined on $\cu$) fulfilled. Then the
 following statements are valid.

1.   $\hw$ is the idempotent completion of the category $H$
(whose objects
are coproducts of $C_i$) in $\cu$.

2. Restrict the functors from $\hrt$ (considered as a subset of
$\hw_*$ by Theorem \ref{sdt}(4)) to the full additive
subcategory $C\subset \hw$ consisting of finite direct sums of
$C_i$. Then this restriction functor gives an equivalence of $\hrt$
with $C_*$.

3. For any object of $Y\in \cu^{t=0}$ and any $X\in\obj \cu$ we
have $\cu(X,Y)=(\ke ( \cu(X^0,Y)\to \cu(X^{-1},Y)) /\imm
(\cu(X^{1},Y)\to \cu(X^0,Y))$ 
where
 $\dots\to X^{-1}\to X^0\to
X^1\to\dots$ 
is an arbitrary choice of a weight complex for $X$.

\end{theo}
\begin{proof}

I1. The existence of $w$ on $\cu^-$ is immediate from part III
(version (ii)) of Theorem \ref{recw}.

Now define  $\cu^{t\ge 0}= (\cu^{w\le -1})^\perp$. In order to prove
that $t$ is a $t$-structure it suffices (cf. Proposition \ref{lads})
to check that for any $X\in \obj \cu$ there exists a
$t$-decomposition (\ref{tdec}).

We will construct $X^{t\le 0}$ and $X^{t\ge 1}$ explicitly.
Our construction uses almost the same argument as the one in the proof
of  version (ii) of Theorem \ref{recw}(III). It may also be
thought about as being a triangulated version of the construction of
Eilenberg-MacLane spaces (this construction really allows
constructing  Eilenberg-MacLane spectra from $S^0$ in $SH$, see
\S\ref{sh} below!).

 We take
$P_0=\coprod_{i\in I,s\in \cu(C_i,X)} C_i$. Then we have a
morphism $f_0:P_0\to X$ whose component that corresponds to
$(C_i,s)$ is given by $s$. Let $X_0$ denote a cone of $f_0$.
Repeating the construction for $X_0[-1]$ instead of $X$ we get an
object $P_1$ being a coproduct of certain $C_i$, $f_1:P_1\to
X_0[-1]$; we denote a cone of $f_1$ by $X_1[-1]$. Proceeding (with
$X_i[-1-i]$), we get an infinite sequence of $(P_i,f_i,X_i)$. We
denote the map $X\to X_0$ given by the construction by $g_0$,
$g_i:X_{i-1}\to X_{i}$, $h_i= g_{i}\circ \dots \circ g_1\circ
g_0:X\to X_i$. We denote a cone of $h_i$ by $Y_i[1]$; the map
$Y_i\to X[1]$ given by the corresponding distinguished triangle by
$r_i$. We have  $P_i\in \cu^{w= 0}$ by  definition.

We have $Y_0=P_0$. Then Remark \ref{rpoto}(2) yields that $Y_i\in
\cu^{w\le 0}$ for all $i$.

Now we  consider the homotopy colimit of $Y_i$; cf. Definition
\ref{dcoulim}. 
By  Lemma \ref{coulim}(1), the sequence $r_i$ can be
lifted to some morphism $f:Y\to X$. We denote its cone as $Z$.
$Y,Z$
will be our candidates for $X^{t\le 0}$ and $X^{t\ge 1}$.

By  Remark
\ref{rcoulim}(3) we have $Y\in \cu^{w\le 0}$ ($=\cu^{t\le 0}$).

We verify that $Z\in \cu^{t\ge 1}$. First we  check that
$C_i[j]\perp Z$ for all $i\in I,\ j\ge 0$.
 This is equivalent to the fact that the map
 $f_*:\cu(C_i[j],Y) \to
\cu(C_i[j],X) $ is an isomorphism for all $i\in I,\ j\ge 0$ and is
injective for $j=-1$ (cf.  Remark \ref{rthick}(1)).

By   Lemma \ref{coulim}(3),  for any compact $C$  we have
$\cu(C,Y)=\inli \cu(C,Y_i)$. Moreover, we have $C_i[-1]\perp Y$
since $Y\in \cu^{w\le 0}$ and  $C_i[-1]\in \cu^{w=1}$. Hence it
suffices to verify that $C_i[j]\perp X_l$ for $l>j\ge 0$ (this
gives $\cu(C_i[j],Y_l)\cong \cu(C_i[j],X)$).

 We apply the distinguished triangle
$P_j[j]\to X_{j-1}\to X_j\to P_j[j+1]$.
Since $C_i[j]$ is compact, we easily obtain $$\cu(C_i[j],P_j[j])=
\bigoplus_{m\in I,s\in \cu(C_m[j],X_{j-1})}\cu (C_i,C_m).$$ Hence
this group has an element for each morphism $ C_m[j]\to X_{j-1}$; it
follows that the map $\cu(C_i[j], P_j[j])\to \cu(C_i[j],X_{j-1})$ is
surjective. Next (since $C_i[j]$ is compact and $C_i[j]\in
\cu^{w=j}$), $\cu (C_i[j], P_j[j+1])=\cu(C_i,P_j[1])$ equals the
direct sum of the corresponding $\cu(C_i,C_m[1])$; hence it is zero by
the orthogonality property for $w$ (cf. Definition \ref{dwstr}). We
obtain $C_i[j]\perp X_{j}$.

Now we consider distinguished triangles $P_l[l]\to X_{l-1}\to X_l\to
P_j[l+1]$ for $l>j$. Again compactness of $C_i[j]$ yields $\cu
(C_i[j], P_l[l+1])= \cu (C_i[j], P_l[l])=\ns$. Hence
$\cu(C_i[j],X_l)= \cu(C_i[j],X_{l-1})=\ns$ for all $l>j$.

 It remains to check that for any $T\in \obj\cu$ the condition
 $\{C_i[j]\}\perp T$
 for all $j\ge 0$ implies that $\cu^{w\le 0}\perp T$.
 This follows immediately from part III (version (ii)) of Theorem
 \ref{recw}.

Moreover, loc. cit. also implies that
 on $\cu^-$ (defined as in part IV of loc. cit.) there
exists a weight structure such that $H\subset \hw$. Note that
$\cu^-$ also satisfies the conditions of the theorem. The
description of $\cu^{w\le 0}$ in the proof of Theorem \ref{recw}
shows that $w$ is left adjacent to $t$ on $\cu^-$.

2. We should check that $w$ can be extended to the whole $\cu$. We
define $X^{w\ge 1}$ 
using the orthogonality axiom (of weight
structures). 

 For any $X\in\obj\cu$ we denote $\tau_{\le i}X$ by
$X_i$ for $i>1$ and take $Y $ being the homotopy colimit of  $
Y_i$ for $Y_i=X_i^{w\ge 1}$ (see Definition \ref{dcoulim}). Here
the morphisms $Y_{i}\to Y_{i+1}$ are obtained by applying 
Lemma \ref{fwd}1) to the natural morphisms $X_i\to X_{i+1}$. By
Lemma \ref{gcoulim} we can assume that $Y\in \cu^{w\ge 0}$.

 $Y$ will be our
candidate for $X^{w\ge 1}$ (cf. the proof of assertion I1). By  Lemma \ref{coulim}(1) the system of   composed maps $Y_i[-1]\to
X_i\to X$ can be lifted to some $f\in \cu( Y[-1], X)$.

Now we show that $f$ extends to a weight decomposition of $X$ using
 Remark \ref{rthick}(1). We should check that
$C_k[j]\perp \co(f)$ for all $k\in I$ and $j<0$ (see the
description of $\cu^{w\le 0}=\cu^{-,w\le 0}$ in the proof of Theorem \ref{recw}(III)). Since all $C_k$ are compact, as in the
proof of assertion I1 we obtain that $\cu(C_k[j],Y)=\inli
\cu(C_k[j],Y_i)$. Moreover, $\cu(C_k[j],X_i)=\cu(C_k[j],X)$ for
$i>-j$. Hence it suffices to note that the direct limit of
isomorphisms is an isomorphism, whereas the  direct limit of surjections
is surjective if the targets stabilize (obvious!).

Now it remains to apply Lemma \ref{lsimple}(3).

II1. Obviously, $\hw$ contains  $\obj H$. Since $\hw$ is
Karoubi-closed in $\cu$, it also contains all retracts of objects
of $H$.

We check that any idempotent $h\in \cu(X,X)$ yields an object of
 $\cu$ if $X\in \obj H$. We apply
 Neeman's argument (see Proposition 1.6.8 of \cite{neebook}).
 One can easily check that the (formal) image of $h$ can be presented as
 $\co f:\coprod_{i\ge 0} X_i\to \coprod_{i\ge 0} X_i$, where all
 $X_i\cong X$; $f_{i,i}=\id_{X_i}$, $f_{i,i+1}=-h$, all other
components of $f$ are zero.

It remains to check that any object of $\hw$ is  a
retract of some object of $H$.

We consider the 'weight resolution' of $X\in \cu^{w=0}$  constructed
 as in the proof of  Theorem \ref{recw}(III)
(in fact, it suffices to consider a few last terms). We obtain that
the weight complex of $X$ can be presented (in $K_\w(\hw)$) by
$\dots \to P_1\to 
P_0$, where $P_i\in \obj H$. Since it is homotopy equivalent to $X$,
 we obtain that $X$
is a  retract of $P_0$ (see  Proposition \ref{bw}(8)).
The assertion is proved.

2. By assertion II1, the restriction  of representable functors to
the category of all coproducts of $C_i$ is  a fully faithful functor on
$\hrt$ (see Theorem \ref{sdt}(4)); 
we can also fully faithfully restrict these functors further to $C$. So
it remains to compute the categorical image of this restriction.

Now, $\cu$ contains all
coproducts of $C_i$. Since all objects of $C$ are compact, these coproducts
represent functors $\bigoplus C(-,C_i)$ on $C$. Since $\hrt$ is
abelian,
 its image also contains all cokernels of morphisms of objects that
can be presented as  $\bigoplus C(-,C_i)$.

It remains to note that cokernels of morphisms of objects of the
type  $\bigoplus C(-,C_i)$ give the whole
 $C_*$. This fact was mentioned in the Notation, see also
  Lemma 8.1 of \cite{vbook}.
  In fact, this is very easy since every $F:C\to\ab^{op}$ can be
presented as a factor of the natural
  $$h:\sum_{i\in I,x\in F(C_i)}C_i\to F,$$
and the same may be said about the kernel of $h$.

3. This is just the formula (\ref{frestrw}).

\end{proof}

\begin{rema}\label{cocomp} 

1. Dualizing assertion I2, one obtains certain sufficient conditions for
 right adjoint weight and $t$-structures to exist. Unfortunately,
this requires 'positive' products and cocompact weak cogenerators
which do not usually exist (yet see \S4.7 of \cite{bger} and \cite{bgern}). 

2. Since $C_i$ are compact, $H$ can be described as the idempotent completion of the category
of 'formal' coproducts of $C_i$, 
i.e.,
$\cu(\coprod_{l\in L}C_{i_l},\coprod_{j\in J}C_{i_j})
=\prod_{l\in L}(\bigoplus_{j\in J}\cu(C_{i_l},C_{i_j}))$;
here $i_j,i_l\in I$, $L,J$ are index sets.

3. If $\cu$ is endowed with a $t$-structure then the question of
existence of an adjacent weight structure seems to be difficult in
general; cf. Remark \ref{gersten} below. Yet see Theorem 4.1 of
\cite{konk} for an interesting result in this direction  (though in
rather restrictive conditions).

\end{rema}

\subsection{The spherical weight structure for the stable homotopy
category}\label{sh}

We consider the (topological) stable homotopy category $SH$. Recall some of its
basic properties.

The objects of $SH$ are called {\it spectra}. $SH$ contains the {\it
sphere spectrum} $S^0$ that is compact and weakly generates it.

The groups $A_i=SH(S^0[i],S^0)$ 
 are called the stable homotopy groups of spheres. We have $A_i=0$
 for $i<0$, $A_i=\z$ for $i=0$; $A_i$ are finite 
  for
 $i>0$. For an arbitrary $A\in\obj SH$ the groups $SH(S^0[i],A)$
 are called the homotopy groups of $A$ (they are denoted by $\pi_i(A)$)).

The category $SH_{fin}\subset SH$ of finite spectra was defined in
Corollary \ref{wsh}. We will also consider the category
$SH_{qfin}\subset SH$ of {\it quasi-finite spectra}. Its objects
are described by the following conditions: all $\pi_i(A)$ are
finitely generated and $\pi_i(A)=0$ for all $i>j$ for some
$j\in\z$. Lastly, we will also mention the full subcategory
$SH^-\subset SH$ whose objects are spectra with (stable) homotopy groups
that are zero for $i>j$ (for some $j$ that depends on the spectrum
chosen). Obviously, all categories mentioned  are triangulated
subcategories of $SH$.

 We see that $SH_{fin}$ and $SH_{qfin}$ satisfy the assumptions of
part III (version (i))
  of
Theorem \ref{recw} if we take $H=H'$ equal to the category of
finite coproducts of $S^0$ and $c=\omega$. Indeed, in this case
we only need finite sums and their properties which are valid for
arbitrary $\cu$.

Hence we obtain a certain non-degenerate weight structure $w$ on
$SH_{fin}\subset SH_{qfin}$. It is   bounded  above for $SH_{qfin}$,
whereas $SH_{fin}$ is bounded since it is generated by $S^0$. Recall
that $S^0$ is compact, hence all objects of $H'$ also are. Hence
using part III (version (ii)) of Theorem \ref{recw} we can extend $w$
to $SH^-$.

Now we  describe the hearts 
all the versions of $w$. Since $SH(S^0,S^0)=\z$, we obtain that $H'\cong
\ab_{fin.fr}$ (the category of finitely generated free abelian
groups); note that
$H=H'$ in this case. Since $H$ it is idempotent complete,
 Theorem \ref{recw}(III(i)) implies
that $\hw_{SH_{fin}}=\hw_{SH_{qfin}}\cong \ab_{fin.fr}$.

In $SH^-$ we have $H\cong\ab_{fr}$ (the category of all  free
abelian groups). Since $\ab_{fr}$ is idempotent complete, we
obtain $\hw_{SH^-}\cong \ab_{fr}$.

Now recall that $SH$ admits countable (and also, in fact, arbitrary small)
coproducts. Hence by Theorem \ref{madts}(I2) we can extend
$w$ to the whole $SH$. This certainly means that
$\hw_{SH}\cong\ab_{fr}$. Hence the functor $t$ is actually 'strong' for all
categories of spectra mentioned, see  Remark \ref{rctst}(1).

Note that any object of $SH^{w=0}$ is isomorphic to a coproduct
of spherical spectra. Hence 
weight Postnikov towers  in
this case become  {\it cellular towers} for spectra in $SH$ in the sense described in (the beginning of) \S6.3 of \cite{marg} (since $\inli X^{w\ge i}=X$; see Proposition 3.3 of ibid. and Proposition \ref{sinhomol}(1) below; more details are given in \S2.4 of \cite{bkw}). 
Their construction and the functoriality properties (in the
topological case) are described in \S6.3 of \cite{marg};
certainly, the results of loc. cit. are parallel to ours.
It follows that the corresponding weight spectral sequences
(for (co)homological functors defined on $SH$) are actually
Atiyah-Hirzebruch spectral sequences, i.e., they relate (co)homology
of an arbitrary spectrum $X$ to the cohomology of $S^0$
(in a way that depends   on the cellular
tower of $X$).

Now we  describe the relation of the weight complex functor for
this weight structure to
singular homology and cohomology of spectra.

To this end  we recall that $SH$ supports a non-degenerate {\it
Postnikov} $t$-structure $t_{Post}$; the corresponding cohomology
functor to $\hrt_{Post}\cong \ab$ is given by $SH(S^0,-)$. We obtain that $SH^{-,w\le
0}=SH^{t_{Post}\le 0}$. Hence $t_{Post}$ in this case is exactly the
$t$-structure described in  Theorem \ref{madts}(I1).
 Besides by  Theorem \ref{sdt}(5), any Eilenberg-MacLane spectrum
belongs to $SH^{t_{Post}=0}$.
 Recall that the singular (i.e., the  $H\z$-) cohomology theory $\hsin^i$ for spectra
is represented by the Eilenberg-Maclane spectrum $H\z$ that
corresponds to $\z$, whereas the singular (the   $H\z$-) homology of $X$ 
is calculated as $SH(S^0,H\z\wedge X[i])$; we will denote it by
$\hsing_i(X)$ (see the notation of Definition \ref{numbers}). 

We identify $\hw=H$ with $\ab_{fr}$ using the functor $H(S^0,-)$; so the
target of $t$ is $K(\ab_{fr})$.

\begin{pr}\label{sinhomol}

Let $X$ be a spectrum.

1. $H_i(t(X))\cong \hsing_i(X)$.

2. $H_0(\ab(X^{-*},\z))\cong \hsin^0(X)$.

\end{pr}
\begin{proof}

1. We apply Theorem \ref{hsps} to the functor $\hsing_0$. We have
$E_1^{pq}=\hsing_q(X^p)$, where each $X^p$ is a (possibly,
infinite) coproduct of copies of $S^0$.
 Now, since the only non-zero homology group of
$S^0$ is $\z$ placed in dimension $0$ and  the functor $Y\to
H\z\wedge Y$ commutes with (small) coproducts,  the spectral sequence $T(\hsing,X)$ reduces to the weight
complex of $X$ (considered as a complex of free abelian groups). 

By
the convergence condition II(iii) of loc. cit. we have
$T(\hsing,X)\implies \hsing(X)$ if $X\in \obj SH^-$. Now, an arbitrary $X\in \obj SH$ can be presented as the homotopy colimit of $\tau_{\le i}X\in SH^{w\le i}$ (see 
the proof of Theorem \ref{madts}(I2)). Since  both the left and the right hand side of our assertion yield homological functors from $SH$ into $\ab$ 
that commute with all small coproducts (cf. Theorem 2.2.6(II.6) of \cite{bgern}), they also respect (countable) homotopy colimits and we obtain the result. 

2. Part II3 of Theorem \ref{madts} calculates $SH(X,Y)$ for any
Eilenberg-MacLane spectrum $Y$. In particular, taking $Y=H\z$ we
obtain the claim.

\end{proof}

\begin{rema}

1. If   we take an Eilenberg-MacLane spectrum $HI$
corresponding to some injective group $I$ instead, we will get
$SH(X,HI)=\ab(H^0(t(X)),I)$.

2. Note also that $S^0{}^{t\ge 0}$ is exactly $H\z$. Hence
$H\z$ can be obtained by applying the construction described in
the proof of Theorem \ref{madts}  to $S^0$.

3. The proof of Proposition \ref{sinhomol}(1) shows that
the weight filtration given by the spherical weight structure on
singular homology coincides with the canonical filtration. This is
not the case for (stable) homotopy groups of spectra.

\end{rema}

\section{Idempotent completions; $K_0$ of categories with bounded
weight structures}\label{idemp}

In \S\ref{ridcomp} we recall that an idempotent completion of a
triangulated category is triangulated. In \S\ref{sidcw} we prove
that  a bounded $\cu$  is idempotent complete (i.e., pseudo-abelian)
 if and only if $\hw$ is; in
general, the idempotent completion of a bounded $\cu$ has a weight
structure whose heart is the idempotent completion of $\hw$.

In \S\ref{sksws} we prove the following: if $\cu$  is bounded and idempotent
complete, then the embedding $\hw\to \cu$ induces an isomorphism
$K_0(\cu)\cong K_0(\hw)$. It is a ring isomorphism if $\hw\subset
\cu$ are endowed with compatible tensor structures. In
\S\ref{skzen} we study a certain  Grothendieck group of
endomorphisms in $\cu$. Unfortunately, it is not always isomorphic
to $K_0(\en \hw)$; yet it is if $\hw$ is {\it regular}; see
Definition \ref{dreg}. Besides, we can still say something about
it in other cases. In particular, this allows us to generalize
Theorem 3.3 of \cite{bloesn} to arbitrary endomorphisms of motives
(in Corollary \ref{genbloesn}); see also \S8.4 of \cite{mymot}.

In \S\ref{skzsh} we calculate explicitly the groups $K_0(SH_{fin})$
and $K_0(\en SH_{fin})$. It turns out that the classes of $[X]$ and
$[g:X\to X]$ are easily recovered from the rational singular
homology of $X$; see Proposition \ref{pkzsh}. More generally, one
can calculate certain groups $K_0(\en^n SH_{fin})$ for $n\in\n$ in
a similar way, see Remarks \ref{rdense} and \ref{rrepr}.

\subsection{Idempotent completions: reminder}\label{ridcomp}

We recall that an additive category $A$ is said to be {\it
idempotent complete} (or pseudo-abelian) if for any $X\in\obj A$
and any idempotent
$p\in A(X,X)$ there exists a decomposition $X=Y\bigoplus Z$ such that
$p=i\circ j$, where $i$ is the inclusion $Y\to X(\cong Y\bigoplus Z)$, $j$ is the projection $X(\cong Y\bigoplus Z)\to Y$. 

Any additive $A$ can be canonically idempotent completed. Its
idempotent completion is (by definition) the category $A'$ whose
objects are $(X,p)$ for $X\in\obj A$ and $p\in A(X,X):\ p^2=p$; we
define
$$A'((X,p),(X',p'))=\{f\in A(X,X'):\ p'f=fp=f \}.$$ It can be easily
checked
that this category is additive and idempotent complete, and for
any idempotent complete $B\supset A$ we have a natural unique
embedding $A'\to B$.

The main result of \cite{ba} (Theorem 1.5) states that an
idempotent completion of a triangulated category $\cu$ has a
natural triangulation (with distinguished triangles being all
direct summands of distinguished triangles of $\cu$).

In this section $\cu'$ will denote the idempotent completion of
$\cu$, $\hw'$ will denote the idempotent completion of $\hw$.

Note that if $\cu$ is idempotent complete then $\hw$  also is, since
$\hw\subset \cu$ and  is Karoubi-closed in it.

\subsection{Idempotent completion of a triangulated category
with a weight structure}
\label{sidcw}

 We prove that $\cu^b$ is idempotent complete if $\hw$ is.

\begin{lem}\label{lbeil}

If $w$ is bounded, $\hw$ is idempotent complete, then $\cu$ also
is.

\end{lem}
\begin{proof}

We prove that all $\cu^{[i,j]}$ are idempotent complete by
induction on $j-i$. The base of induction is the fact that $\cu^{[i,i]}=\cu^{w=0}[-i]$ is
idempotent complete.

To make the inductive step it suffices to prove that $\cu^{[-i,1]}$
if idempotent complete if $\cu^{[-i,0]}$ is (for $i>0$). For $X\in
\cu^{[-i,1]}$ and an idempotent $p\in \cu(X,X)$ we consider the
functor $WD$ (see  Theorem \ref{fwc}(I)). We obtain an
idempotent  $q=WD(p)\in K^{[0,1]}_\w(\cu)(WD(X),WD(X))$, whereas
$Y=WD(X)$ has the form $(Z\to T)$ for some $Z,T \in \cu^{[-i,0]}$. Since $\cu^{[-i,0]}$ is idempotent complete,
$K^b_\w(\cu^{[-i,0]})$ also is by  Proposition \ref{wkar}(2).
Moreover, loc. cit. yields  the existence of a morphism
$Z'\to T'$ and idempotent endomorphisms $r,s$ of $Z'$ and $T'$,
respectively, such that $(Y,q)$ can be presented by the diagram
$$\begin{CD}
Z'@>{}>>T'\\
@VV{r}V@VV{s}V \\
Z'@>{}>>T'
\end{CD}$$
 (in
$K^{[0,1]}_\w(\cu^{[-i,0]})$).

By part I5 of Theorem \ref{fwc}, $Z',T'$ come from a certain
weight decomposition of $X$. Then any corresponding weight
decomposition of $p$ is homotopy equivalent to $(r,s)$.
 Then  
Theorem \ref{fwc}(I2) yields that $(r,s)$ also give a weight
decomposition of $p$. Hence
 the object $(X,p)\in\obj \cu'$ (see \S \ref{ridcomp}) can
be presented as a cone of a certain map  $(Z',r)\to (T',s)$ in
$\cu'$; whereas $(Z',r),(T',s)\in \obj\cu$ by the inductive
assumption.

\end{proof}

Now we prove that in the general (bounded) case a weight structure
can be extended from $\cu$ to its idempotent completion $\cu'$.

\begin{pr}\label{picompl}

Suppose that $w$ is bounded.
Then the following statements are valid.

1. $w$ extends to a bounded weight structure $w'$ for $\cu'$.

2. 
$\hw'$ (i.e., the idempotent completion of $\hw$) generates
 $\cu'$.

3. The   heart of $w'$ is equivalent to $\hw'$. 

\end{pr}

\begin{proof}

1,2. By part II1 of Theorem \ref{recw}, we have a bounded weight
structure that extends $w$ on the subcategory  $D\subset \cu'$
generated by $\hw'$. Hence  it suffices to recall that $D$ is
idempotent complete; see Lemma \ref{lbeil}.

3. Since $\hw'$  is idempotent complete, the assertion follows from
 Theorem \ref{recw}(II2).

\end{proof}

\begin{rema}\label{roblom}

Possibly  the boundedness condition on $w$
 in Proposition \ref{picompl} can be weakened. However this does not
seem to be actual
  since usually (for the triangulated categories that interest mathematicians) either $(\cu,w)$ is bounded or $\cu$
admits countable coproducts (at least, 'positive' or 'negative' ones).
 In the latter case $\cu$ is idempotent
complete, see Proposition 1.6.8 of \cite{neebook}.
\end{rema}

\subsection{$K_0$ of a triangulated category with a bounded
 weight structure}
\label{sksws}

 We recall some standard definitions (cf. 3.2.1 of
\cite{gs}). We define the Grothendieck group of an additive
category $A$  as the Abelian group whose generators are of the form $[X],\
X\in\obj A$; the relations are 
$[B]=[C]+[D]$ if 
 $B\cong C\bigoplus D$ for 
$B,C,D\in\obj A$. The $K_0$-group of a triangulated category $T$ is
defined as the Abelian group whose generators are $[t],\ t\in \obj T$; if
$D\to B\to C\to D[1]$ is a distinguished triangle then we set
$[B]=[C]+[D]$. Note that  $X\bigoplus 0\cong X$ implies that
$[X]=[Y]$ if $X\cong Y$ (in $A$ or in $T$).

For an additive $A$ we define $K_0(K^b_\w(A))$ similarly to
$K_0(K^b(A))$; hence it equals $K_0(K^b(A))$ (see Definition
\ref{dkw}).

 The existence of a  bounded $w$
allows us to calculate $K_0(\cu)$ easily.

\begin{theo}\label{kzw}
Suppose that $(\cu,w)$ are  bounded and that $\hw$ be idempotent complete. Then
the inclusion $i:\hw\to\cu$ induces an isomorphism $K_0(\hw)\to
K_0(\cu^b)$.
\end{theo}

\begin{proof}
 Since $t$ is an weakly exact functor  (see Definition
\ref{dkw}), it gives an abelian group homomorphism $a:K_0(\cu)\to
K_0(K_\w^b(\hw))=K_0(K^b(\hw))$. By Lemma 3 of 3.2.1 of \cite{gs},
there is a natural isomorphism $b:K_0(K^b(\hw))\to K_0(\hw)$. The
embedding $\hw\to \cu$  gives a homomorphism $c:K_0(\hw)\to K_0(\cu)$.
The definitions of $a,b,c$ imply immediately that $b\circ a\circ
c=\id_{K_0(\hw)}$. Hence $a$ is surjective, $c$ is injective.

It remains to verify that $c$ is surjective. This follows
immediately from the fact that $\hw$ generates $\cu$, see
  Corollary \ref{gcu}.

\end{proof}

\begin{rema}\label{rthenzkz}
Obviously, if $\cu$ is a tensor triangulated category then
$K_0(\cu)$ is a ring. If the tensor structure on $\cu$ induces a
tensor structure on $\hw$, then $K_0(\hw)$ is a ring also and $c$ is
a ring isomorphism.

\end{rema}

For the convenience of citing  we  concentrate certain assertions
relevant for motives in a single statements.

\begin{pr}\label{unext}

Suppose that $\cu$ contains an additive  negative (see Definition
\ref{negth}) subcategory $H$ such that $H$ is idempotent complete
and $\cu$ is the idempotent completion of $\lan H\ra$. Then the
following statements are valid.

1. $\lan H\ra=\cu$.

2. There exists a conservative  weight complex functor $\cu\to
K_\w^b(H)$ which sends $h\in \obj H$ to $h[0]\in \obj K^b_\w (H)$.
It can be lifted to an exact functor $t^{st}:\cu\to K^b(H)$ in the
case when $\cu$ has a differential graded enhancement  (see  Definition \ref{dprtp}(3) below).

3. $K_0(\cu)\cong K_0(H)$.

\end{pr}
\begin{proof}

1.  By  Theorem  \ref{recw}(II) there exists a bounded weight
structure $w'$ on $\lan H\ra$ whose heart is equal the small envelope
of $H$, i.e., to $H$ itself. Next, by Proposition \ref{picompl} $w$ extends
to some bounded $w$ on $\cu$ whose heart is equal to the idempotent
completion of $H$, i.e., to $H$ again. Hence 
Proposition \ref{picompl}(2) immediately
yields assertion 1.

2. The weight complex functor $t:\cu\to K_\w(H)$ can be factored
through $K_\w^b(H)$ since $w$ is bounded. $t$ is conservative by
 Theorem \ref{wecomp}(V). If $\cu$ has a differential graded
enhancement then $t$ can be lifted to $t^{st}$ by  Remark
\ref{rwndg}(3) below.

3. Immediate from Theorem \ref{kzw}.

\end{proof}

\subsection{$K_0$ for categories of endomorphisms}\label{skzen}

Now we define various    Grothendieck groups of endomorphisms in an
additive category $A$. Our definitions are similar to those of
\cite{alm}.

\begin{defi}\label{kzen}

1. The  generators of $K_0^{add}(\en A)$ are endomorphisms of
objects of $A$; we impose the relations   $[g]=[f]+[h]$ if
$(f,g,h)$ give an endomorphism of a split short exact sequence.

2. If $A$ is also abelian then we also consider the group
$K_0^{ab}(\en A)$. Its  generators again are endomorphisms of
objects  of $A$; we set  $[g]=[f]+[h]$ if
$(f,g,h)$ give an endomorphism of an arbitrary short exact
sequence.

3.  For a triangulated $A$ we  consider the group
$K_0^{tr}(\en A)$.

Its  generators are endomorphisms of objects
of $A$ again; we set $[g]=[f]+[h]$ if $(f,g,h)$ give an
endomorphism of a distinguished triangle in $A$.

\end{defi}

Note that  $K_0^{ab}(\en A)$ and  $K_0^{tr}(\en A)$   are natural
factors of $K_0^{add}(\en A)$ (when they are defined).
Indeed, $K_0^{ab}(\en A)$ and  $K_0^{tr}(\en A)$   have the same
generators as  $K_0^{add}(\en A)$ and more relations.

Suppose that $\cu$ is bounded. We provide some sufficient conditions for
$K_0(\en\cu)$ to be isomorphic to   $K_0(\en \hw)$. We need a
notion of  {\it regular} additive category $A$. Recall that
$A_*'$ is the full abelian subcategory of $A_*$ generated by $A$.

\begin{defi}\label{dreg}

An  additive category  $A$ will  be called {\it regular} if it
satisfies the following conditions.

1. $A$ is isomorphic to its small envelope (see  Definition
\ref{negth}(2)), i.e., if  $X,Y\in \obj A$, $X$ is a retract of $Y$, is
 then $X$ has a complement to $Y$ (in $A$).

2. Every object of $A_*'$ has a finite resolution by objects of
$A$.

\end{defi}

The most simple examples of regular categories are abelian
semisimple categories and the category of finitely generated
projective modules over a Noetherian (commutative) local ring all
of whose localizations are regular local; cf. the end of \S1 of
\cite{alm}.

We will need the following technical statement. Let $R$ be an
associative  ring with a unit.

\begin{lem}\label{lreg}

1. If $A$ is regular then $K_0^{add}(\en A)\cong K_0^{ab}(\en
A_*')$.

2. If $A$ is the category of finitely generated projective modules
over $R$ then $A_*$ is the category of all (left) modules over
$R$.

\end{lem}
\begin{proof}

1. We apply the method of the proof of Proposition 5.2 of
\cite{alm2}. First we consider the obvious category $\en \hw_*'$
and note that it is abelian. Next, the objects of $\hw$ become
projective in $\hw_*'$. Hence all 3-term complexes in $\hw$ that
become exact in $\hw_*'$ do split in $\hw$. Therefore we can define
$K_0^{add}(\en \hw_*)$ as the Grothendieck group of an exact
subcategory of $\en \hw_*'$.

Condition 1 of Definition \ref{dreg} ensures that for any short
exact sequence $0\to G'\to G\to G''\to 0$ in $\en \hw_*'$ if  $G,
G''\in \en \hw$ then $G'\in \en \hw$ (i.e., $G'$ is an endomorphism of
an object of $\hw$). Lastly, condition 2 of Definition \ref{dreg}
easily implies that any $G\in \en \hw_*'$ has a finite resolution by
objects of $\en \hw$ (again note that objects of $\hw$ become
projective in $\hw_*'$!). Hence applying Theorem 16.12 of \cite{swan}
(page 235) we obtain the result.

2. The equivalence is given by sending a functor $F$ to $F(R)$
(here $R$ 
is considered as right $R$-module; this endows $F(R)$ with a
 left $R$-module structure);
 and
an $R$-module $Q$ to $P\mapsto \homm_R(P,Q)$. Note that all $F(P)$
can be uniquely recovered from $F(R)$ since all finitely
generated
projective modules are direct summands of $R^m$ (for some $m>0$). 

\end{proof}

\begin{pr}\label{pkzen}

1. There exist natural homomorphisms $K_0(\en \hw)\stackrel{c}{\to}
K^{tr}_0(\en \cu) \stackrel{d}{\to} K^{ab}_0(\en \hw_*') $; $c$ is
a surjection.

2. $c$ is an isomorphism if $\hw$ is regular. 

\end{pr}
\begin{proof}

1. $c$ is induced by $i:\hw\to \cu$. For $g:X\to X$ we define

\begin{equation}\label{kzco}
d(g)=\sum (-1)^i[
g_{i*}:H_i(t(X))\to H_i(t(X))].\end{equation}
 Here
$H_*(t(X))\in\obj \hw_*'$ is the homology of the weight complex;
see  Remark \ref{rkw}(2). We obtain a well-defined
homomorphism since $t$ is a weakly exact functor (see Definition
\ref{dkw}); see Remark \ref{rkw}(3).

$c$ is surjective since for $g:X\to X$ we have the equality
$[g]=\sum (-1)^i[g^i:X^i\to X^i]$. This equality follows easily from
the fact that a repetitive application of  the (single, shifted)
weight decomposition functor to a morphism yields its infinite
weight decomposition (see Theorem \ref{fwc}; note that $X$ is
bounded).

2. In the case when $\hw$ is abelian semi-simple we have
$\hw=\hw_*'$. Hence the equality $d\circ c=\id_{K_0(\en \cu)}$
yields the assertion (in this case).

Now, in the general (regular) case it suffices to apply the
equality $K_0(\en \cu)=K_0(\en \hw_*')$ (this is  Lemma
\ref{lreg}(1)).

\end{proof}

\begin{rema}\label{rtkzen}

1. Unfortunately, $c$ is not an isomorphism in the general case. To
see this it suffices to consider the example described in 
Remark \ref{bmorwc}(3) for $\cu=K^b(Z)$, where $Z$ is the category
 of free $\z/4\z$-modules. This fact is also related to the
observation in the end of \S1 of \cite{alm}. Certainly, $Z$ is not
regular.

2. Certainly, if $i:\hw\to \cu$ is a tensor functor then  $c,d$ are
ring homomorphisms, cf. Remark \ref{rthenzkz}.


\end{rema}

The surjectivity of $c$ immediately implies the following fact.

\begin{coro}\label{genbloesn}
Let $r:\cu\to D^b(R)$ and  $s:\cu\to D^b(S)$ be  exact functors for
 abelian $R,S$; let $r_*:K_0(\cu)\to K_0(D^b(R))$ and
$s_*:K_0(\cu)\to K_0(D^b(S))$ be the induced homomorphisms. Let
$u:K_0(\en D^b(R))\to K_0(\en R)$ and $v:K_0(\en D^b(S))\to K_0(\en
S)$ be defined as $(g:X\to X)\to [g_{i*}:H_i(X)\to H_i(X)]$.
Let $T$
be an abelian group; $x:K_0(\en R)\to T$ and $y:K_0(\en S)\to T$ be
group homomorphisms. Then the equality $x\circ u\circ r_*\circ
c=y\circ v\circ s_*\circ c$ implies $x\circ u\circ r_*=y\circ v\circ
s_*$.

In particular, one can take $\cu=\dmge\q$, $\hw=\chowe\q$ (see
Remark \ref{rperf} below), $r,s$  given by $l$-adic cohomology
realizations (for two distinct $l$'s), and $x,y$  given by traces of
endomorphisms. It follows that the alternated sum of traces of maps
induced by a $g\in \dmge\q(X,X)$ on the cohomology of $X$ does not
depend on $l$. We also obtain the independence from $l$ of
$n_{\lambda}(H)=(-1)^i n_{\lambda} g^*_{H_i(X)}$; here $n_{\lambda}
g^*_{(H_i(X))}$ for a fixed algebraic $\lambda$ denotes the
algebraic multiplicity of the eigenvalue ${\lambda}$ for the
operator $g^*_{H_i(X)}$.

This generalizes Theorem 3.3 of \cite{bloesn} to arbitrary
morphisms of motives; see \S8.4 of \cite{mymot} for more
details.

\end{coro}

Lastly, we consider some more general $K_0$-groups.

\begin{rema}\label{rrepr}
1. For an additive $A$ instead of $\en A$ one can for any $n\ge
0 $ consider the category $\en^n A$ whose objects are the
following $n+1$-tuples: $(X\in\obj A;g_1,\dots, g_n\in A(X,X))$.
We have $\en^0(A)=A$, $\en^1(A)=\en A$. Generalizing  Definition
\ref{kzen} in an obvious way one defines $K_0^{add}(\en^n A)$,
$K_0^{ab}(\en^n A)$, and $K_0^{tr}(\en^n A)$ (for $A$ additive,
abelian or triangulated, respectively). Next, one can define $c,d$
as in Proposition \ref{pkzen}; exactly the same argument as in the
proof of the Proposition shows that $c$ is always surjective and
it is also injective if $\hw$ is regular. In particular, this is
true for $\cu=SH_{fin}$; see Proposition \ref{pkzsh} below.

2. Even more generally, for any ring $R$ one may consider the
category $\en(R,A)$ of $R$-{\it representations} in $A$, i.e., of
pairs $(X,H:R\to A(X,X))$; here $X\in\obj A$, $H$ is a unital 
homomorphism of rings. In particular, we have  $\en(R,A)=A$ for
$R=\z$, $=\en A$ for $R=\z[t]$, and $=\en^n A$ for $R=\z\lan
t_1,\dots, t_n\ra$ (the algebra of non-commutative polynomials).
Again one defines $K_0^{add}(\en (R,A))$, $K_0^{ab}(\en (R,A))$, and
$K_0^{tr}(\en (R,A))$, $c$ and $d$. Yet the method of the proof of
Proposition \ref{pkzen} fails for a general $R$; one can only note
that $d\circ c$ is an isomorphism if $\hw$ is abelian semi-simple.
\end{rema}

\subsection{An application: calculation of $K_0(SH_{fin})$ and
$K_0(End\, SH_{fin})$}\label{skzsh}

Now we calculate explicitly the groups $K_0(SH_{fin})$ and
$K^{tr}_0(\en SH_{fin})$. The author doesn't think that (all of)
these results are new; yet they illustrate our methods  very well.

 We will  need the following simple observation: $K_0(A)$ is
naturally a direct summand of $K_0(\en A)$ (both in the
'triangulated' and in the 'additive' case). The splitting is
induced by $[f:X\to X]\to [X]\to [0:X\to X]$; see \S1 of 
\cite{alm}.

We define the group $\lamb$ as a subgroup  of the multiplicative
group $\lamb(\z)=\{1+t\z[[t]]\}$ that is generated by polynomials
(with constant  term $1$). $\lamb$ and  $\lamb(\z)$ are also
rings; see Proposition 3.4 of \cite{alm2} for  $\lamb$ and
\cite{haz} for $\lamb(\z)$.

\begin{pr}\label{pkzsh}

1. $K_0(SH_{fin})\cong \z$ with the isomorphism sending $X\in \obj
SH_{fin}$ to $[X]=\sum (-1)^i\dim_\q (\hsing_i(X)\otimes \q)$ (the
rational singular homology of $X$).

2. $K^{tr}_0(\en SH_{fin})\cong \z\bigoplus \lamb$ with the
isomorphism sending $g:X\to X$ to $[X]\bigoplus \prod_i (\det_{\q[t]}
(id-g_it\otimes\q))^{(-1)^i}$; here $g_it\otimes\q$ is the map
induced by $g\otimes t$ on $\hsing_i(X)\otimes_\z \q[t]$.
\end{pr}
\begin{proof}
1. We have $\hw=\ab_{fin.fr}$ for the spherical weight structure
$w$ on $SH_{fin}$; see \S\ref{sh}. Hence $K_0(SH_{fin})\cong
K_0(\ab_{fin.fr})=K_0(\z)=\z$.

 The second assertion can easily be
deduced from Proposition \ref{sinhomol}(1). Note that
$K_0(SH_{fin})$ is a direct summand of $K^{tr}_0(\en SH_{fin})$;
hence $$[X]=\sum (-1)^i[H_i(t(X))]=\sum (-1)^i[\hsing_i(X)]$$ by
(\ref{kzco}). We also use the fact that $K_0(\z)$ injects into
$K_0(\q)$, so $[\hsing_i(X)]$ can be computed rationally.

2. By  Lemma \ref{lreg}(2) we have $\hw_*\cong \ab_{fin.fr}$
(the category of finitely generated abelian groups). Hence $\hw$ is
regular (see Definition \ref{dreg}). Therefore by 
Proposition \ref{kzen}(2) we have $K^{tr}_0(\en SH_{fin})\cong
K_0^{add}(\en \ab_{fin.fr})$. Then the Main Theorem in \S1 of
\cite{alm} implies that $K^{tr}_0(\en SH_{fin})\cong \z\bigoplus
\lamb$.

Next, (\ref{kzco}) implies $[g]=(-1)^i[g_{i*}]$. Now note that
$\lamb(\q)\to \lamb(\z)$ is injective; so it suffices to calculate
$[g_{i*}]$ rationally. Lastly, the equality
$$[g_{i*}\otimes\q]=\dim_\q (\hsing_i(X)\otimes \q)\bigoplus
\det{}_{\q[t]} (id-g_it\otimes\q)$$ follows from the formula at
the bottom of p. 376 of \cite{alm}.

\end{proof}

\begin{rema}\label{rdense}
1. Note that the isomorphisms described are compatible with the
natural ring structures of $K_0$-groups involved.

2. Assertion 1 doesn't seem to be new; yet the author doesn't know
of any paper that contains assertion 2 in  its current form.

3. One also has $K^{tr}_0(\en^n SH_{fin})\cong K^{add}_0(\en^n
\ab_{fin.fr})$; see Remark \ref{rrepr}.

\end{rema}

\section{Twisted complexes over negative differential graded
categories; Voevodsky's motives}\label{dgmot}

The goal of this section is to apply our theory to triangulated
categories that have differential graded enhancements (as considered in \cite{mymot}); this will allow us to apply it for
motives.

In \S\ref{sbdedg} we recall the definitions of differential graded
categories and twisted complexes over them. In \S\ref{snegdg} we
consider negative differential graded categories; we obtain a weight
structure on the category of twisted complexes (over them). In
\S\ref{cwecom} we
 construct the 
 {\it truncation
functors} $t_N$; $t_0$ is the {\it strong weight complex functor}
for this case (see Conjecture \ref{ctst}).

In \S\ref{twmot} we recall the spectral sequence $S(H,X)$
constructed in \S7 of \cite{mymot} for $H$ having a differential
graded enhancement, and prove that it can be obtained from $T(H,X)$
by means of a certain shift of indices. In particular, this shows that $S$
does not depend on the choice of enhancements. We also prove that {\it
truncated realizations} for representable realizations are
represented by the adjacent $t$-truncations of representing objects
(see also \S\ref{strfun} and \S\ref{ddme}).

In \S\ref{remmymot} we apply our theory to Voevodsky's motivic
categories $\dmge\subset \dmgm$.

We calculate the heart of the {\it Chow} weight
structures obtained in \S\ref{hchow}. We also recall that we can
apply this theory with rational coefficients over a perfect field
$k$ of arbitrary characteristic.

\subsection{Basic definitions}\label{sbdedg}

 We recall relevant definitions for differential graded categories as
they were presented in \cite{mymot}; cf. also \cite{bev} and
\cite{bk}. 

Categories of {\it twisted complexes}
 were first considered in \cite{bk}. However
our notation differs slightly from that of \cite{bk}; some of the
signs are also different.

An additive category $C$ is called graded if for any $P,Q\in\obj C$
there is a canonical decomposition $C(P,Q)\cong \bigoplus_i C^i(P,Q)$
defined;
 this decomposition should satisfy
 $C^i(*,*)\circ C^j(*,*)\subset C^{i+j}(*,*)$.
 A differential graded category (cf. \cite{bk} or \cite{bev})
is a graded category
 endowed with an additive operator
$\delta:C^i(P,Q)\to C^{i+1}(P,Q)$ for all $ i\in \z, P,Q\in\obj C$.
$\delta$ should satisfy  the equalities $\delta^2=0$ (so $C(P,Q)$ is
a complex of abelian groups); $\delta(f\circ g)=\delta f\circ
g+(-1)^i f\circ \delta g$ for any $P,Q,R\in\obj C$, $f\in C^i(Q,R)$,
$g\in C(P,Q)$. In particular, $\delta (\id_P)=0$.

We denote $\delta$ restricted to morphisms of degree $i$ by
$\delta^i$.

For any additive category $A$ one can construct the following
differential graded categories.


We denote the first one by $S(A)$. We take $\obj S(A)=\obj A;\
S(A)^i(P,Q)=A(P,Q)$ for $i=0$;
 $S(A)^i(P,Q)=0$ for $i\neq 0$.
We take $\delta=0$.

We also consider the category $B^b(A)$ whose objects are the same as
for $C^b(A)$, whereas for $P=(P^i)$, $Q=(Q^i)$ we define
$B^b(A)^i(P,Q)=\bigoplus_{j\in \z} A(P^j, Q^{i+j})$.  Obviously
$B^b(A)$ is a graded category. $B(A)$ will denote the unbounded
analogue of $B^b(A)$.

We take $\delta f=d_Q\circ f-(-1)^i f \circ d_P$, where $f\in
B(A)^i(P,Q)$, $d_P$ and $d_Q$ are the differentials in $P$ and $Q$.
Note that the kernel of $\delta^0(P,Q)$ coincides with $C(A)(P,Q)$
(the morphisms of complexes); the image of $\delta^{-1}$ are the
morphisms homotopic to $0$.

$B^b(A)$ can be obtained from $S(A)$  by means of the category
functor $\prt$ described below.

For any differential graded $C$ we define a category $H(C)$; its
objects are the same as for $C$; its morphisms are defined as
$$H(C)(P,Q)=\ke \delta^0_C(P,Q)/\imm \delta^{-1}_C(P,Q).$$

Having a differential graded category $C$ one can construct
another differential graded category $\prt(C)$ as well as a
triangulated category $\tr(C)$. The simplest example of these
constructions is $\prt(S(A))=B^b(A)$.

\begin{defi}
 The objects of $\prt(C)$ are  $$\{(P^i),\ P^i\in\obj C, i\in\z,
q^{ij}\in C^{i-j+1}(P^i,P^j)\};$$
 here almost all $P^i$ are $0$;
 for any $i,j\in \z$ we have
 \begin{equation}\label{dif}
\delta q^{ij}+\sum_{l\in \z} q^{lj}\circ q^{il}=0\end{equation}
 We call $q^{ij}$ {\it
arrows} of degree $i-j+1$. For $P=\{(P^i),q^{ij}\}$,
$P'=\{(P'^i),q'^{ij}\}$ we set
$$\prt(C)^l(P,P')=\bigoplus_{i,j\in\z}C^{l+i-j}(P^i,P'^j).$$ For
$f\in C^{l+i-j}(P^i,P'^j)$ (an arrow of degree $l+i-j$) we define
the coboundary of the corresponding morphism in $\prt(C)$ as
$$\delta_{\prt(C)}f=\delta_C f+\sum_m(q'^{jm}\circ
f-(-1)^{(i-m)l}f\circ q^{mi}).$$

\end{defi}

It can be easily seen that $\prt(C)$ is a differential graded
category (see \cite{bk}). So we denote $H(\prt(C))$ by $\tr(C)$.

There is  an obvious translation 
functor on $\prt(C)$. Note also that the terms of the complex
$\prt(C)(P,P')$ do not depend on $q^{ij}$ and $q'^{ij}$, whereas the
differentials certainly do.

We denote by $Q[j]$ the object of $\prt(C)$ that is obtained by
putting $P^i=Q$ for $i=-j$, all other $P^j=0$, all $q^{ij}=0$. We
will write $[Q]$ instead of $Q[0]$.

Immediately from the definition we have $\prt(S(A))\cong B^b(A)$.

A morphism $h\in\ke \delta^0$ (a closed morphism of degree $0$) is
called a {\it twisted morphism}; note that twisted morphisms are exactly the ones giving morphisms in  $\tr(C)$. For a twisted morphism
$h=(h^{ij})\in \prt((P^i,q^{ij}),(P'^i,q'^{ij}))$, $h^{ij}\in
C^{i-j}(P^i,P'^j)$ we define $\co(h)=(P''^i,q''^{ij})$, where
$P''^{i}=P^{i+1}\bigoplus P'^{i}$,
$$q''^{ij}=\begin{pmatrix}q^{i+1,j+1}&0\\
h^{i+1,j}&q'^{ij}\end{pmatrix}$$ 

We have a natural triangle of twisted morphisms
 \begin{equation}\label{dgtri}
P\stackrel{f}{\to} P'\to \co(f)\to P[1],
\end{equation}
the components of the second map are $(0, id_{P'^i})$ for $i=j$
and $0$ otherwise. 
This triangle induces a triangle in the
category $H(\prt(C))$.

Now we
list the main properties of categories of twisted complexes. 


\begin{defi}\label{dprtp}

1. For distinguished triangles in $\tr(C)$ we take  the triangles
isomorphic to those that come from  (\ref{dgtri}) for
$P,P'\in \prt(C)$.

2. $\trp(C)$ is defined as the full (strict) triangulated category of $\tr(C)$ generated by 
$[P]:\ P\in \obj C$; we denote
the corresponding full subcategory of $\prt(C)$ by $\prtp(C)$.

3. We will say that $\cu$ admits a {\it differential graded
enhancement} if it is equivalent to $\trp(C)$  for some differential
graded $C$.
\end{defi}

We summarize the properties of the categories defined that are most
relevant for the current paper; see \cite{mymot} and \cite{bk} for
the proofs.

\begin{pr}\label{mdg}

I For any additive category $A$ there are natural isomorphisms

1. $\prt (B(A))\cong B(A)$.

2. $\tr(B(A))\cong K(A)$.

3.  $\tr (S(A))\cong K^b(A)$

II
 1. There are natural embeddings of categories $i:C\to
\prtp(C)$ and $H(C)\to \trp(C)$ sending $P$ to $[P]$.


2.  $\prt$, $\tr$, $\prtp$, are $\trp$ are functors on the category
of differential graded categories, i.e., any differential graded
functor $F:C\to C'$ naturally induces functors
 ${\prt}F$, ${\tr}F$, $\prtp{}F$, and $\trp{}F$.

3. Let $F: \prtp(C)\to D$ be a differential graded functor. Then the
restriction of $F$ to $C\subset \prt(C)$ gives a differential graded
functor $FC:C\to D$. Moreover, since $FC=F\circ i$, we have
$\prtp(FC)= \prtp(F)\circ \prtp(i)$; therefore $\prtp(FC)\cong
\prtp(F)$.

4.  $\trp(C)\subset \tr(C)$ (with the distinguished triangles described in Definition \ref{dprtp}) are triangulated categories.

\end{pr}

For example, for $X=(P^i,q^{ij})\in \obj \prt(C)$ we have
${\prt}F(X)=(F(P^i),F(q^{ij}))$; for a morphism $h=(h^{ij})$ of
$\prt(C)$ we have ${\prt}F(h)= (F(h^{ij}))$. Note that the
definition of ${\prt}F$ on morphisms does not involve $q^{ij}$; yet
${\prt}F$ certainly respects differentials for morphisms.

\subsection{Negative differential graded categories; a weight
structure for
$\tr(C)$}\label{snegdg}

Suppose now that a differential graded category $C$ is {\it
negative}, i.e., for any $X,Y\in \obj C$ and $i>0$ we have 
$C^i(X,Y)=\ns$ (cf. Definition \ref{negth}).

For $\cu=\tr(C)$ 
we define $\cu^{w\le 0}$ as a class that contains all
objects  isomorphic to those that satisfy $P^i=0$ for $i>0$.
$\cu^{w\ge 0}$ is defined similarly by the condition $P^i=0$ for
$i<0$.

\begin{pr}\label{wndg}

1. $\cu^{w\le 0}$ and $\cu^{w\ge 0}$
yield a bounded weight structure for $\cu$.

2. $\hw$ is isomorphic to the small envelope of   $HC$ in $\cu$ (cf.
Definition \ref{negth}).

\end{pr}
\begin{proof}
The definition of morphisms in $\cu$ immediately yields that
 $\cu^{w\ge 0}\perp \cu^{w\le 0}[1]$. We obviously have
$\cu^{w\le 0}[1]\subset \cu^{w\le 0}$; $\cu^{w\ge 0}\subset
\cu^{w\ge 0}[1]$. The verification of the fact that $\cu^{w\le 0}$
and $\cu^{w\ge 0}$ are Karoubi-closed in $\cu$ is straightforward.
However we will never actually use this statement
below (so we can replace $\cu^{w\le 0}$ and $\cu^{w\ge 0}$ described
by their
Karoubi-closures in the definition of $w$). 

It remains to check that any object $X$ of $\cu$ admits a weight
decomposition. We follow  the  proof  of Proposition 2.6.1 of
\cite{mymot}.

We take  $(P^i,f^{ij},\  i,j\le 0)$ as $X^{w\le 0}$ and
$(P^i,f^{ij},\   i,j\ge 1)[1]$ as $X^{w\ge 1}$. We should verify
that $X^{w\le 0}$  and $X^{w\ge 1}$ are objects of $\cu$.

We have to check that the equality (\ref{dif}) is valid for $X^{w\le
0}$ (resp. $X^{w\ge 1}$).  All terms of (\ref{dif}) are zero unless
$i\le j \le 0$ (resp. $1\le i\le j$). Moreover, in the case $i\le j
\le 0$ (resp. $1\le i\le j$)  the terms of (\ref{dif}) are the same
as for $X$. Both of these facts follow immediately from the
negativity of $C$.

Now we verify that $(\id_{P^i},\ i\le 0)$ gives a morphism $X\to
X^{w\le 0}$ and $(\id_{P^i},\ i\le 1)$ gives a morphism $ X^{w\ge
1}[-1]\to X$. Indeed, for these morphisms the equality
 $\delta_{\prt(C)}f=0$ is  obvious
by the negativity of $C$.

Next we should check that $X\to X^{w\le 0}$ is the second morphism
of the triangle corresponding to $ X^{w\le 1}[1]\to X$; this
easily follows from (\ref{dgtri}).


2. Obviously, the objects of $HC$ belong to $\cu^{w=0}$. Next, the
definition of $\cu$ easily yields that $\hw(X,Y)\cong HC(X,Y)$ for
$X,Y\in \obj HC$.

Moreover, assertion 1 implies that any object of $\cu$ has a 'filtration'
by  subobjects whose 'successive factors' come from $HC$. By  Theorem \ref{recw}(II2) we obtain that $\hw$ is isomorphic to the
small envelope   of $HC$.
\end{proof}

Obviously, the same construction also gives weight structures for
all unbounded versions of $\tr(C)$.

\begin{rema}\label{rwndg}

1. Alternatively, Proposition \ref{wndg} can be deduced from 
 Theorem \ref{recw}(II). In particular, this method easily allow
 proving
that $\cu^{w\le 0}$ and $\cu^{w\ge 0}$ are Karoubi-closed
using the fact  that the small envelope of $HC$ lies in both of
them (cf. the beginning of the proof of part II2 loc. cit.).

2. Since all objects of $\tr(C)$ possess Postnikov towers whose 'factors' belong to $HC$, we obtain that $\tr(C)=\trp(C)$
(when $C$ is negative!).

3. Let $\cu=\trp(C)$ for some $C$ such that $H_iC(X,Y)=0$ for all
 $i>0$, $X,Y\in \obj
C$ (the homology of $C(-,-)$ is concentrated in non-positive
degrees). Let $C_-$ be a (non-full!) subcategory of $C$ with
 the same objects and $C_-(X,Y)=C(X,Y)^{t\le 0}$ (morphisms are the
zeroth canonical truncation of those of $C$). Then  by 
Remark 2.7.4(2) of \cite{mymot} the embedding $C_-\to C$ induces an
equivalence of triangulated categories $\trp(C_-)\to\trp(C)$.

 It
follows that if $\cu\cong \trp(C)$ and $C(-,-)$
 is acyclic in positive degrees then we can assume $C$ to be
  negative. In particular, the strong (i.e.,
exact) weight complex functor $\cu\to K^b(\hw)$ exists in this case
(see below).
\end{rema}

\subsection{Truncation functors; comparison of weight
complexes}\label{cwecom}

For $N\ge 0$, $P,Q\in \obj C$ (for a negative $C$) we denote the $-N$-th canonical
filtration of $C(P,Q)$ (i.e., $C^{-N}(P,Q)/d_PC^{-N-1}(P,Q)\to 
C^{-N+1}(P,Q)\to\dots \to C^0(P,Q)\to 0$) by $C_N(P,Q)$.

We denote by $C_N$ the differential graded category whose  objects
are the same as for $C$, whereas the morphisms are given by
$C^i_{N}(P,Q)$. The composition of morphisms is induced by that in
$C$. For morphisms in $C_N$ presented by $g\in C^{i}(P,Q)$, $h\in
C^{j}(Q,R)$, we define their composition as the morphism represented
by $h\circ g$ for $i+j\ge -N$ and zero for $i+j<-N$. Certainly, all
$C_N$ are negative (i.e., all morphisms of positive degree are
zero).

We have an obvious functor $C\to C_N$. By   Proposition
\ref{mdg}(II2), 
 this gives canonically a functor $t_N: \cu\to
\tr(C_N)$. We denote $\tr(C_N)$ by $\cu_N$.

Obviously,  objects of $\cu_N$ can be represented as certain
$(P^i, f^{ij}\in C_N^{i-j+1}(P^i,P^j), i<j\le i+N+1)$, the morphisms
between $(P^i, f^{ij})$ and $(P'^i, f'^{ij})$ are represented by
certain $g^{ij}\in C_{N}^{i-j}(P^i,P'^j),\ i\le j\le i+N$, etc. The
functor $t_N$ 'forgets' all elements of $C^m([P],[Q])$ for $P,Q\in
\obj C$, $m<-N$, and factorizes $C^{-N}([P],[Q])$ modulo
coboundaries. In particular, for $N=0$ we get ordinary complexes
over $HC$, i.e., $\cu_0=K^b(HC)$.

$t_0$ will be called the {\it strong weight complex} functor.

One can easily verify that the strong weight complex functor
constructed  is a lift of the weight complex functor $t$
corresponding to  the weight structure $w$ to an exact functor
$t^{st}$ (as in Conjecture \ref{ctst}). This  follows immediately
from the explicit description of $X^{w\le 0}$ and $X^{w\ge 1}$ for
any $X\in \obj\cu$ (in the proof of Proposition \ref{wndg}).

\begin{conj}\label{cortrf}
1. For a general $(\cu,w)$ there also exist  certain exact {\it higher truncation
functors} $t_N$ such that $t_0$ is the 'strong' weight complex
functor; cf. Conjecture \ref{ctst}. Their targets $\cu_N$ should
satisfy the following conditions:  if $X,Y\in \cu^{w=0}$ then
$\cu_N(t_N(X),t_N(Y)[-i])=\cu(X,Y)$ for $0\le i\le N$ and $=\ns$
otherwise. These categories should admit full embeddings
$i_N:\cu^{[0,N]}\to \cu_N$;  distinguished triangles of $\cu$
consisting of elements of $\cu^{[0,N]}$ should be mapped to
distinguished triangles by $i_N$.

2. Let $I:\cu\to D(A)$ be an exact functor, where $\cu,w$ is a
triangulated category with a weight structure, $A$ is an abelian
category. If $I(\cu^{w=0}) \subset D_{[0,N]}(A)$ (i.e., acyclic for
degrees outside $[0,N]$) then $I$ can be  factored
through $t_N$.
\end{conj}

\subsection{Weight spectral sequences
for enhanced realizations}\label{twmot}

The method of construction of  weight spectral sequences in
\cite{mymot}
 was somewhat distinct
from the method we use here. In \cite{mymot} we have used a certain
filtration on the complex that computes cohomology; this filtration
can be obtained from the filtration corresponding to our current
method by Deligne's decalage (see \S1.3 of \cite{de2} or
\cite{paran}). So the spectral sequence there was 'shifted one level
down' (in particular, it was functorial in $X$ starting from
$E_1$ already). We
compare the methods here.

Let $J$ be some negative differential graded category, let
$\hk=\tr(J)$, $J'=\prt(J)$. Below we will use the same
 notation for Voevodsky's motives
(which are the most important example of this situation).

In \cite{mymot} weights were constructed only for (co)homological
 functors that admit an enhancement, i.e.,
those that can be factored through $\tr(F)$ for a differential
 graded functor $F:J\to C$.
 Here  we consider only $C=B(A)$ for an abelian $A$ and
homological functors of the form $H_{K_A}\circ \tr(F)$
 (here $H_{K_A}$ denotes the zeroth homology functor for $C(A)$).
 The cohomological functor case (also with $C=B(A)$) was considered
in \S7.3 of
 \cite{mymot} (certainly, reversing the arrows is no problem).
  Note still that for those realizations for which  $C\neq B(A)$ one can
  sometimes
 reduce the situation to the case $C=B(A')$ for a 'large' $A'$
 (in particular, this seems to be the case for the singular realization of motives with the values
 in the category of mixed Hodge complexes).

Now we recall the formalism of \cite{mymot} (modified for the
homological functor case).

We denote the functor $\prt(F):J'\to B(A)$ by $G$, denote
$\tr(F):\hk\to K(A)$ by $E$. It turns out that the
 virtual $t$-truncations
of $F$ (see Remark \ref{rtrfun}) have nice differential
graded enhancements.

We recall that for a complex $Z$ over  $A$, $b\in\z$,
 its $b$-th canonical truncation from above is the complex
 $ \dots \to Z^{b-1}\to \ke (Z^b\to Z^{b+1}), $ here $\ke
(Z^b\to Z^{b+1})$ is put in degree $b$.

For any $b\ge a\in \z$ we consider the following functors. By
$F_{\tau_{\le b}}$ we denote the functor that sends $[P]$ to
 $\tau_{\le b}(F([P]))$. These functors are differential graded;
hence they extend to $G_i=\prt(F_{\tau_{\le - i}}):J'\to B(A)$. Note
that we consider the $-i$-th filtration here in order to make the
filtration decreasing (which is usual when the decalage is applied);
this is another minor distinction of the current exposition from
that of \cite{mymot}.  The functors $\tr(F_{\tau_{\le - i}})$ were
called truncated realizations also in loc. cit.

Let $X=(P^i,q^{ij})\in \obj J'$.
 The complexes $G_b(X)$ give a
 filtration of $G(X)$; one may also consider
$G_{a,b}(X)=G_b(X)/G_{a-1}(X)$.
 We obtain the spectral sequence of a filtered complex
(see \S III.7.5 of \cite{gelman})
  \begin{equation} S:E_1^{ij}(S)\implies H^{i+j}(G(X)).\label{spectr}\end{equation}
 Here $E_1^{ij}(S)=H^{i+j}(G_{1-j}(X)/G_{-j}(X))$.

All $G_b(X)$  are $J'$-contravariantly functorial with respect to
$X$. Besides, starting from $E_1$ the terms of $S$ depend only on
the homotopy classes of $G_b(X)$. Hence starting from $E_1$ the
terms of $S$ are functorial with respect to $X$ (considered as an
object of $\hk$).

Now we compare  spectral sequences obtained using this method with
those provided by Theorem \ref{csps}. In fact, the comparison
statement can be proved by considering the derived exact couple
for $T$; see (part 3 of) Remark \ref{trfun}. Alternatively, we could
have extended
 Theorem
\ref{sdt}; see  Remark \ref{rsdt}(5,4).
Instead, here we give a proof in terms of filtered complexes.

To this end we compare the filtrations of $G(X)$ corresponding to
$T$ and $S$. Fortunately, we don't have to write down the
differential in $G$; it suffices to recall that
$G^j(X)=\bigoplus_{k+l=j} F_k(P^l)$.

The method of Theorem \ref{hsps} gives the following filtration on
$G(X)$: $Q_i G^j(X)=\bigoplus_{k+l=j, l\ge i} F_k(P^l)$.

Now we apply decalage to this filtration. It is easily seen that we
obtain the filtration given by $G_i$, i.e., $$(Dec
Q)_i(G^j(X))=\bigoplus_{k+l=j, l\ge j+i+1} F_k(P^l) \bigoplus \ke
(F_{-b}(P^{l+i})\to F_{-b+1}(P^{l+i})).$$

Hence $T_{n+1}^{pq}=S_n^{-q,p+2q}$ for all integral $i,j$ and
$n>0$; the corresponding filtrations on the limit (i.e., on
$H^{i+j}(E(X)))$
coincide up to a certain shift of indices.

 In \S7.3 of \cite{mymot}  so-called {\it truncated realizations }
 were
considered. They were defined as $\tr(F_{\tau\le b})$ and
$\tr(F_{\tau\le b}/F_{\tau\le a-1}) :\hk\to K(A)$ (for $a\le
b\in\z$). The formula (15) of \cite{mymot} computes all
$E_n^{ij}(S)$ for $n\ge 1$ in terms of the weight filtration of
truncated realizations  of $X$; this description is
$\hk$-functorial.

\begin{rema}\label{tmotco}

1. Suppose now that there exists a differential graded functor
$F^1:J\to B(A)$ and a differential graded transformation $F^1\to
F$ such that the induced homology functor morphisms are isomorphisms
in degrees $\le b$ and are zero in degrees $>b$. Let $F^2$ denote
$F^1_{\tau\le b}$. We have a natural transformation $F^2\to F^1$
which is an isomorphism on homology. Hence by 
Corollary 2.7.2(2) of \cite{mymot}, the  transformation of functors
$\tr(F^2\to F^1)$ induces quasi-isomorphisms of their values.
Next, the transformation $F^1\to F$ induces a transformation
$F^2\to F_{\tau\le b}$. Applying  Corollary 2.7.2(2) of
\cite{mymot} again we obtain that $\tr(F_{\tau\le b})\approx
\tr(F^2)$; hence both of them are quasi-isomorphic to $\tr(F^1)$.

In particular, let $A=\ab$; let $F$ be (contravariant) representable
by some $Y$ in some differential graded $K\supset J$ such that $\tr
K$ possesses a weight structure extending $w$ and its adjacent
$t$-structure $t$. Then our reasoning shows that the objects
$t_{\le i}Y$ represent the truncated realizations for
$\tr(K(-,Y))$ (in $\tr(K)\supset \hk$; up to quasi-isomorphism, i.e.,
they give the homology groups required). This is a differential
graded version of  Theorem \ref{sdt}(6). Besides, in this case
the
 fact that the
filtrations induced by the morphisms $X\to w_{\le i}X$ and by
$t_{\le i}Y\to Y$ coincide also follows from 
Theorem \ref{sdt}(6).

So, the results of \S\ref{ddme} below yield that the truncated
realizations both for the 'classical' (Weil) realizations of motives
and for motivic cohomology are representable by objects of $\dme$.
This fact seems to be far from obvious.

2. As was mentioned in  Remark \ref{rsdt}(4), one can deduce
the comparison of spectral sequences statement in the representable
case from the remark above and Theorem \ref{sdt}. Moreover, one
can construct a {\it nice duality} of $\cu,w$ with the $t$-structure on the
category $\du=\tr(DG-Fun(J,B(A)))$ (differential graded functors)
that corresponds to the canonical truncation of $A$-complexes; see
Definition 2.5.1 of \cite{bger} and 
 Remark \ref{rsdt}(5) below.  Also, $\du$ could be called 'a
category of functors $\cu\to A$'; see 
Remark
\ref{rtrfun}. The realizations of the type considered above
correspond to some objects of this category; truncations of a
realization with respect to this $t$-structure would be exactly its
truncated realizations (cf. Proposition 2.5.4 of \cite{bger}).

\end{rema}

\subsection{ $\smc$, $\dmge$ and $\dmgm$; the Chow weight
structure}\label{remmymot}

We recall some definitions of \cite{1}.

$k$ will denote our perfect ground field; we will mostly assume that
the characteristic of $k$ is zero. $\pt$ is a point, $\af^n$ is the
$n$-dimensional affine space (over $k$), $x_1,\dots,x_n$ are the
coordinates, $\p^1$ is the projective line.

$\var\supset \sv\supset \spv$ will denote the class of all
varieties over $k$, resp. of smooth varieties, resp. of smooth
projective varieties.

We define the category of smooth correspondences as follows: $\obj \smc=\sv$,
$\smc (X,Y)=\bigoplus_U\z$ for all $U\subset X\times Y$ that are integral
closed finite subschemes which are surjective over a connected
component of $X$. The elements of $\smc(X,Y)$ are called {\it
finite correspondences} from $X$ to $Y$.

\begin{rema}\label{multva}
The composition of $U_1\subset X \times Y$ with $U_2\subset Y \times
Z$ as in the definition of finite correspondences is defined as
always in the categories of motives, i.e., one considers the obvious
scheme-theoretic analogue of $\{(x\in X,z\in Z):\exists y\in Y:\
(x,y)\in U_1,\ (y,z)\in U_2\}$.
Note that the composition
is well-defined without any factorization by equivalence relations
needed. Next one extends composition to all $\smc(-,-)$ by
linearity.

Note that this definition is compatible with the naive notion of
composition of multivalued functions. Now, to $\sum c_iU_i\in
\smc(X,Y)$, $c_i\neq 0$ one can associate a multi-valued function
whose graph is $\cup U_i$. Applying this definition, one can
define images and preimages of finite correspondences (and their
restrictions). Below we will assume that images are closed integral
subschemes of the corresponding $Y$'s.

\end{rema}

$\smc$ is additive; the addition of objects is given by the disjoint
union operation for varieties. It is also a tensor category; the
tensor product operation is given by the Cartesian product of
varieties.

$\ssc$ is the abelian category of those additive  cofunctors
$\smc\to\ab$
that are sheaves in the 
Nisnevich topology.

$\dme\subset \dmk$ is defined as the subcategory defined by the
condition that
 the cohomology  sheaves are homotopy invariant
(i.e., $S(X)\cong S(X\times \af^1)$ for any $S\in\sv$).

There is a natural functor $RC\circ L: K^b(\smc)\to \dme$ (cf.
Theorem 3.2.6 of \cite{1}) given by Suslin complexes (see below);
it can be factorized as the composition of the 'localization by
homotopy invariance and Mayer-Viertoris' and a full embedding; it
categorical image will be denoted by $\dms$. 
One can restrict $RC\circ L$  to obtain a functor $\mg:\sv\to \dms$.
Moreover, one can extend $\mg$ to $\var$ (see \S4.1 of \cite{1});
unfortunately, in the case $\cha k>0$ one would have to take $\dme$
as the target of this (extended) $\mg$. Therefore, cohomology of
varieties can be expressed in terms of cohomology of motives.

$\dme$ is idempotent complete; hence it contains the idempotent
completion of $\dms$ which is
 Voevodsky's $\dmge$ (by definition; see \cite{1}).

Now we define a differential graded category $J$ with $\obj J=\spv$
(the addition of objects is the same as for $\smc$). The morphisms
of $J$ are given by {\it cubical Suslin complexes}
 $J^i(Y,P)\subset
\smc(\af^{-i}\times Y,P)$ consisting of correspondences that 'are
zero if one of the coordinates is zero'. Being more precise, we consider
$C'^i(P,Y) =\smc(\af^{-i}\times Y,P)$ for all $P,Y\in\sv$; note that
$C'^i$ are zero for positive $i$.
 For all $1\le j\le -i$,
$x\in k$, we define $d_{ijx}=d_{jx}:C'^i\to C'^{i+1}$ as
$d_{jx}(f)=f\circ g_{jx}$, where $g_{jx}:\af^{-i-1}\times
Y\to\af^{-i}\times Y$ is induced by the map $(x_1,\dots,x_{-1-i})\to
(x_1,\dots,x_{j-1},x,x_j,\dots,x_{-1-i})$. We define $J^i(Y,P)$ as
$\cap_{1\le j \le -i}\ke d_{j0}$. The boundary maps
$\delta^i:J^i(-,-)\to J^{i+1}(-,-)$ are defined as $\sum_{1\le j\le
-i}(-1)^jd_{j1}$.

The composition of morphisms in $J$ is induced by the obvious
composition $C'^i(Y\times \af^{-j},X\times\af^{-j})\times
C'^j(Z,Y)\to C'^{i+j}(Z,X)$ combined with the embedding of
$C'^i(Y,X)$  into $C'^i(Y\times \af^{-j},X\times\af^{-j})$ via
'tensoring' its elements
 by $\id_{\af^{-j}}$; here  $X,Y,Z\in\spv,\ i,j\le 0$.

It was checked in \S2 of \cite{mymot} that $J$ is a differential
graded category. It is negative by definition.

We denote $\tr(J)$ by $\hk$. $\hk$ is equivalent to $\dms$ (if
$\cha k=0$) by Theorem 3.1.1 of \cite{mymot}.

By Proposition \ref{wndg} we obtain that there exists a weight
structure $w$ in $\hk$; hence it also gives a weight structure for
$\dms$. We have $\hw=J_0'$, where $J_0'$ is  the small envelope of
$J_0=HJ$ (cf. Definition \ref{negth} and  Theorem
\ref{recw}(II2)).

\begin{rema}\label{mgc}

1. In \S4.1 of \cite{1} the {\it motif with compact support} for any
$X\in \var$ was defined as the Suslin complex of a certain sheaf
$L^c(X)$. For a proper $X$ we have $\mg^c(X)=\mg(X)$. However, in
order to increase the chances to obtain a geometric motif (with
compact support) one can define $\mg^c(X)$ using Poincare duality;
see Appendix B of \cite{hubka}. In the case $\cha k=0$ these
definitions coincide and yield an object of $\dms$ for any $X\in
\var$.

In Theorem 6.2.1 of \cite{mymot} it was proved that for a smooth $X$
we have $\mg(X)\in \dms^{w\ge 0}$, $\mg^c(X)\in \dms^{w\le 0}$.
Using the blow-up distinguished triangle (see Proposition 4.1.3 of
\cite{1}) one can also show that for a proper $X$ we have
$\mg(X)=\mg^c(X)\in \dms^{w\le 0}$.

2. As in  Remark \ref{semim}(3) one may consider
semi-motives $W_i(\mg(X))$ and $W_i(\mg^c(X))$ for all $i\in \z$,
$X\in \var$; they lie in $\dmge_*$. We obtain that
$W_0(\mg(X))=\mg(X)_*$ for proper $X$, whereas $W_{-1}(\mg(X))=0$
for  $X\in\sv$. Recall that (\ref{contrvar}) allows us to express the
weight filtration on the cohomology of (the motif of) $X$ in terms
of $W_i(\mg(X))$.
\end{rema}

In $\dmge$ we have a decomposition $[P^1]=[pt]\bigoplus \z(1)[2]$ for
$\z(1)$ being the {\it Tate motif}. Moreover, $\dmge$ is a tensor
category with $\otimes \z(1)$
 being a full embedding of $\dmge$ into itself
(the {\it Cancellation Theorem}, see Theorem 4.3.1 of \cite{1}
and \cite{voevc}). Hence one can define Voevodsky's $\dmgm$ as
the direct limit of $\dmge$
 with respect to tensoring by $\z(1)$;
  it also may be described as the union of $\dmge(-i)$
   (whereas each $\dmge(-i)$ is isomorphic to $\dmge$).


\begin{pr}\label{wsm}

$w$ extends to a weight structure for $\dmge$ and  $\dmgm$.
\end{pr}
\begin{proof}

I Extending $w$ to $\dmge$.

We define $\dmge{}^{w\le 0}$ as the set of retracts of
$\dms{}^{w\le 0}$ in $\dmge$; the same for $\dmge{}^{w\ge 0}$. By
Proposition \ref{picompl}, this gives a weight structure on
$\dmge$.

II Extending $w$ to $\dmgm$.

We note that tensoring by $\z(1)[2]$ sends $[P]$ to a retract of
$[P\times \p^1]$. Hence $\otimes\z(1)[2]$ maps $\dmge^{w\le 0}$ and
$\dmge^{w\ge 0}$
 into themselves. It follows that one can define $\dmgm^{w\le 0}$
 and $\dmgm^{w\ge 0}$
as $\cup \dmge^{w\le 0}(-i)[-2i]$ and $\cup\dmge^{w\ge 0}(-i)[-2i]$,
respectively. Indeed, the Cancellation Theorem gives us
orthogonality; since each object of $\dmgm$ belongs to
$\dmge(-i)=\dmge(-i)[2i] $ for some $i\in\z$, we also have the
weight decomposition property.

\end{proof}

\begin{rema}\label{wild} 

Note that (for any $\cu$)  if $w$ is bounded then $\cu^{w\le 0}$
consists exactly of objects that can be 'decomposed' into a
weight Postnikov tower 
(see Definition \ref{dpoto}) with $X^k=0$ for $k>0$; for
$X\in \cu^{w\ge 0}$ we can assume that $X^k=0$ for $k<0$.

Besides (see Proposition \ref{wndg}) for $\cu=\dms$ we can assume
that all $X^k$ can be presented as $\mg(P^k)$ for $P^k\in \spv$.
For $\cu=\dmge$ or $\cu=\dmgm$ we have $P^k\in \obj \chow \subset
\obj \dmgm$ (see \S\ref{hchow} below).

\end{rema}

We call the weight structure constructed  the {\it Chow} weight
structure (for any of $\hk$, $\dms$, $\dmge$, $\dmgm$, and also
for $\dme$ considered below).

\subsection{The heart of the Chow weight structure}\label{hchow}

Now we calculate the hearts of $w$ in each of the categories constructed.

First we recall the definition of (homological) Chow motives. In
\cite{1} it was proved that $\chow$ can be described in the
following way. One considers $\cho=J_0$; this is (essentially) the
usual category of rational correspondences. The category $\chowe$
is the idempotent completion of $\cho$; it was shown in Proposition
2.1.4 of \cite{1} that $\chowe$ is naturally isomorphic to the usual
category of effective homological Chow motives. Moreover, the
natural functor $\chowe\to \dmge$ is a full embedding (of additive
categories). Note that $\chowe$ is a tensor category.

$\chow$ will denote the whole category of Chow motives, i.e.,
$\chow[\z(-1)[-2]]$. 

So, the heart of the Chow weight structure for $\dms$ is the small
envelope of  $\cho$ (note that it contains $\z(1)[2]$, whereas $\cho$
does not). Now Proposition \ref{picompl} implies that the heart of
$\dmge$ is the idempotent completion of $\cho$, i.e., the whole
category $\chowe$. Lastly, we obtain that the heart of $\dmgm$
equals $\chow$.

We obtain that for {\bf any} (co)homological functor from $\dmge$
(or $\dmgm$)
 there exist (Chow-)weight spectral sequences and weight
filtrations. Note that we don't need any enhancements here (in
contrast to  \cite{mymot})! Moreover, Chow-weight spectral sequences
are functorial with respect to all natural transformations of
(co)homological functors (so, we also do not need any transformations for
enhancements).

\begin{rema}\label{rperf}
The same arguments as above also prove the existence of weight
structures on rational hulls of $\dms$, $\dmge$ and $\dmgm$ (i.e., we
tensor the groups of morphisms by $\q$) as well as on their
idempotent completions (which do not coincide with $\dmge\otimes \q$
and $\dmgm\otimes\q$). If we denote the latter by $\dmge\q\subset
\dmgm\q$, then their hearts will be $\chowe\q\subset \chow\q$ (i.e.,
the idempotent completions of rational hulls). Note that in these
statements one can take $k$ being an arbitrary perfect field (of any
characteristic, since the use of the resolution of singularities in
the proofs can be replaced by de Jong's alterations, see \cite{dej}).
 See \S8.3 of
\cite{mymot} and \S\ref{motperf} below for more details.

Moreover, recent (unpublished) results of O. Gabber imply that for
$\operatorname{char}k=p$ one can  prove our results
for  motives with $\z[1/p]$-coefficients.  The author 
will treat this
matter in a forthcoming paper.

Alternatively, one may  consider motives with
$\z/n\z$-coefficients (for some $n>1$ and prime to $p$).
\end{rema}

Lastly we note (as we also did in \cite{mymot}) that the results obtained (the existence of weight filtrations and of Chow-weight spectral sequences) 
also concern
motivic cohomology of motives; cf.   Remark
\ref{rwfiltcoh}(2).

\section{Some new results on motives}\label{othermot}

The first  subsection is dedicated to the study of $\dme$. We prove
that the Chow weight structure extends to it;  $\dme$ also
supports a {\it Chow} $t$-structure $\tcho$ that is (left) adjacent  to it. It
follows that the Chow $t$-truncations of those objects that represent the
'classical' realizations of motives (or motivic cohomology) represent
their  {\it truncated realizations}; see Remark \ref{tmotco}.

In \S\ref{smm} we note that any construction of the  motivic
 $t$-structure on
 $\dmge\q$ would automatically
yield a canonical weight filtration for
the objects of its heart (i.e., for mixed motives).

In \S\ref{motperf} we prove that a certain (possibly, 'infinite')
weight complex functor can be defined for motives over any perfect
field (without any resolution of singularities assumptions).

In \S\ref{coniv} we apply the philosophy of adjacent structures to
express the cohomology of a certain motif $X$ with  coefficients in the homotopy
($t$-structure) truncations of any $H\in \obj\dme$ in terms of the
limit of $H$-cohomology of certain 'submotives' of $X$. 
Luckily, to this end (instead of the {\it Gersten weight structure} that is
 constructed in \cite{bger} only in the case of a countable $k$)
it
suffices to have Gersten resolutions for homotopy invariant
pretheories (constructed in \cite{3}). 

In \S\ref{motet} we recall that (by the Beilinson-Lichtenbaum
conjecture which was recently proved) torsion motivic cohomology is
the (homotopy) truncation of the (torsion) \'etale one. Hence one can express
torsion motivic cohomology of certain motives in terms of \'etale
cohomology of their 'submotives'. In particular, we obtain a formula
for (torsion) motivic cohomology with compact support of a
smooth quasi-projective variety.

In \S\ref{sunram} we calculate $\tcho$ in certain 'simple' cases; it
turns out that it is closely related to unramified cohomology!


\subsection{Chow weight and $t$-structures for $\dme$}
\label{ddme}

We recall (see \S3 of \cite{1}) that for any $S\in\dme$ and
$X\in\sv$ we have $\dme(\mg(X), S)=\h_0(S)(X)$ (here $S$ is
considered as a complex of sheaves). It follows (cf. \S 
\ref{remmymot}) that $\mg(X)$ for $X\in\spv$ weakly generate
$\dme$.

Now we take $\{C_i\}=\obj\chow\subset\obj \dme$ (we can assume
that $\obj\chow$ is a set). We obtain that $(\dme,\{C_i\})$
satisfy the conditions of  Theorem \ref{madts}(I1). Hence
 it has a $t$-structure whose heart is $\chowe_*$; we will denote
it by $\tcho$.
Unfortunately, it seems that  $\tcho$ cannot be restricted to
$\dmge$ (i.e., it is not 'geometric').

\begin{rema}\label{rihom}
Recall that for any $X\in \obj \dmge\subset \obj \dme$ there is
 a natural exact functor $\ihom(X,-):\dme\to \dme$ defined
(see Remark 14.12 of \cite{vbook}). Since   tensor products
(in $\dme$) of  Chow motives are Chow motives also, we obtain the following
 (cf. Proposition \ref{padad}):
 for any $X\in \obj \chow(\subset \obj \dmge)$ we can restrict
$\ihom(X,-)$ to an exact functor $H\tcho\to H\tcho$.
In particular, $H\tcho$ admits certain 'negative Tate twists'
(defined by $\ihom(\z(1)[2],-)$).

Note (in contrast) that
 the homotopy $t$-structure for $\dme$ (see below) is
respected by $\ihom(\z(1)[1],-)$.

\end{rema}

Now we check that the Chow weight structure of $\dmge$ can be
extended to $\dme$. Till the end of this section $\cu=\dme$, $t$
will denote the homotopy $t$-structure of $\dme$ (defined as in
\cite{1}). This is the $t$-structure corresponding to Nisnevich
hypercohomology i.e., $X\in \cu^{t\le 0}$ (resp. $X\in \cu^{t\ge 0}$)
if and only if its Nisnevich hypercohomology is concentrated in
non-positive (resp. non-negative) degrees. Note that in all
$\cu^{t\le i}$ arbitrary coproducts exist.

We define $\cu^{w\le 0}$ as the Karoubi-closure in $\cu$ of the
closure of $\dmge{}^{w\le 0}$ in $\cu$ with respect to arbitrary
coproducts and  'extensions'
(as in Definition \ref{exstab}). Note that
$\cu^{w\le 0}\subset \cu^{t\le 0}$. We recover $\cu^{w\ge 0}$ from
$\cu^{w\le 0}$ via the orthogonality condition (in the usual way, see
 Lemma \ref{lsimple}(3)).  $\cu^{w\ge 0}$ is extension-stable
 (see 
Lemma \ref{lsimple}(1)). Besides, it contains arbitrary coproducts of
objects of $\dmge{}^{w=0}$ (here we apply the compactness of objects
of $\dmge$ in  $\cu$).

As usual, the only non-trivial axiom check here is the
verification of the existence of weight decompositions. Recall
that any object of $\ssc$ has a certain canonical resolution by direct
sums of $L(X)=\smc(-,Y)$ for $Y\in\sv$ (placed in degrees $\le 0$;
see \S 3.2 of \cite{1}). Hence any object $X$ of $\dme$ is a
homotopy colimit of certain $X_i$ ($i\in\z$)
 such that a cone of $X_i\to X_{i+1}$
is a coproduct of  certain $\mg(Y_{ij})[i]$; $X_l=0$ for some
$l\in\z$. The  limit of $X_i$ equals $X$ indeed by Lemma
\ref{tcoulim}.

We construct $Z=X^{w\ge 1}$ as a homotopy colimit of $X_l^{w\ge
1}$ (see Definition \ref{dcoulim}). Note that 
a weight
decomposition of $X_l$ may be constructed using any possible
weight decompositions of $\coprod\mg(Y_{ij})[i]$ (see Remark
\ref{rtwd}).

We should check that the colimit exists in $\dme$.  For any $Y\in\sv$, $i>0$, we have $(\mg(Y)[i])^{w\ge
1}\in \cu^{t\le 0}$ (for any choice of $(\mg(Y)[i])^{w\ge 1}$).
This is easy since $\mg(Y)[i]\in \cu^{t\le -i}$ and
$(\mg(Y)[i])^{w\le 0}\in \cu^{w\le 0}\subset \cu^{t\le 0}$.
Combining these statements for all $Y_{ij}$ and $i$ yields the
boundedness required.

We have the composed morphisms $X_l^{w\ge 1}\to X_l[1]\to
X[1]$; by Lemma \ref{coulim}(1) this system of morphisms
can be lifted to some morphism  $Z\to X[1]$. We should check
that it yields a weight decomposition (if we make the choices in
the construction in a 'clever' way). By Lemma \ref{gcoulim} we can
assume that $Z\in \cu^{w\ge 0}$. We denote a cone of $Z[-1]\to X$
by $Y$.

Now, it suffices (see  Remark \ref{rthick}(1)) to check that
the induced map $\cu(C,Z)\to\cu(Z,X[1])$ is an isomorphism for any
$C\in\cu^{w\ge 1}$ and is surjective for $C\in\cu^{w\ge 0}$.

First suppose that for some $i\in \z$ we have $C\perp R$ for any
$R\in\cu^{t\le i} $. Then the sequence $\cu (C[1],X_i)$ stabilizes;
this yields the result required by  Lemma \ref{coulim}(2).
Hence for any such $C$ we have $C\perp Y[1]$.

We denote ${}^\perp(Y[1])\cap \cu^{w\ge 0}$ by
$S$. We should check that $S=\cu^{w\ge 0}$. Certainly,
$S$ is extension-stable and closed with respect to arbitrary coproducts
(in $\cu$).  

We have $\dmge^{w\ge 0}\subset S$. Indeed, any $C\in \obj \dmge$ is
a retract of an object that can be obtained from (a finite
number of) motives of smooth varieties by considering cones of
morphisms; whereas for $X\in\sv$ we have $\dme (\mg(X),R)=\ns$ for
any $R\in\cu^{t\le -\dim X-1} $ (since the Nisnevich cohomological
dimension of a scheme is not greater than its dimension).
Next, all coproducts of objects of $\dmge^{w\ge 0}$ (belonging to
$\dme$) also belongs to $S$. Therefore, it suffices to prove that
any object of $\cu^{w\ge 0}$ can be 'approximated' by such
coproducts.

By the same method as above, we present $C\in \cu^{w\ge 0}$ as a
homotopy colimit of certain $C_i$ for a cone of $C_i\to C_{i+1}$
being a coproduct of  some $\mg(E_{ij})[i]$; $C_l=0$ for some
$l\in\z$.

Since any coproduct of distinguished triangles is a distinguished
triangle, we can construct distinguished triangles $(\coprod
\mg(E_{ij})[i])^{\ge 0}\to A_i\to B_i$ for $A_i \in \cu^{w\le 0}$
and  $B_i \in S$ (they will be coproducts of objects of $\dmge$).
Next,
applying Remark \ref{rtwd} for $D=\cu^{w\le 0}$, $E=S$,
 we can (starting from $C_l$) inductively construct
distinguished triangles $C_i\to F_i\to G_i$ for $F_i\in \cu^{w\le
0}$, $G_i\in S$. We also construct distinguished triangles
$C_i[-1]\to L_i\to M_i$ for $L_i\in \cu^{w\le 0}$, $M_i\in S$.

By Definition \ref{dcoulim}, we have a distinguished triangle
$\coprod C_i\to  \coprod C_i\to C$. Now note that $\coprod F_i$,
$\coprod G_i$, $\coprod L_i$ and $\coprod M_i$ exist in $\cu$
(since
by the same argument as the one used above all of the summands
belong to $\cu^{t\le l}$ for some $l\in\z$). Since any coproduct
of distinguished triangles in $\cu$ is a distinguished triangle
and $\cu^{w\le 0}$ and $S$ are closed with respect to all coproducts,
we obtain
 distinguished triangles $\coprod C_i\to \coprod
F_i\to \coprod G_i$ and
$\coprod C_i[-1]\to \coprod L_i\to \coprod
M_i$ with $\coprod F_i, \coprod L_i\in \cu^{w\le 0}$, $\coprod
G_i,\coprod M_i\in S$.

 Applying
 Remark \ref{rtwd} again we obtain a distinguished triangle
$C[-1]\stackrel{f}{\to}U\to V$.
 for some  $U\in \cu^{w\le 0}$ and $V\in S$.
Hence $f=0$; therefore $C$  is a retract of $V$. Thus $C\in S$.


\begin{rema}\label{gersten} 
 Unfortunately,
one cannot define a weight structure for $\dme$   that would be
left adjacent
to the homotopy $t$-structure. Indeed, for an object $X$ of this
heart the
functor $\dme(X,-)$ should be exact on the category of homotopy invariant
sheaves with transfers. So the heart should contain
'motives of points' i.e., motives of local smooth $k$-algebras; the
latter are (usually) pro-$k$-varieties and not varieties.
Yet
in \S4.1 of \cite{bger} we
define the corresponding {\it Gersten} weight structure in a certain
triangulated
$\gds\supset \dmge$; see also  Remark \ref{rws}(4) below.

\end{rema}

\subsection{Weight filtration for (conjectural) mixed
motives}\label{smm}

Suppose now that there exists a so-called  motivic
$t$-structure $t_{MM}$ on $\dmge$ or  on $\dmge\q$ (then one can extend
 it to
$\dmgm$  or to $\dmgm\q$, respectively). We will not discuss  its
properties here (until \S\ref{stfiltr}); however it would automatically
 induce a homological
functor $H^{MM}:\dmge\to MM$ for some abelian category $MM$ (of
so-called mixed motives) that is the heart of the $t$-structure.
Hence for any $X\in\dmge$ there will be a certain (weight)
filtration on $H_i^{MM}(X)$. 
This
filtration would be trivial (i.e., 'canonical') when $X$ is smooth
projective. It can be easily checked that there can exist only one
filtration on $H^{MM}_i(X)$ which is $\dmge$-functorial and
satisfies this property.

 Moreover, any transformation $H_{MM}\to H$ for
$H$ being a realization (of $\dmge$) with values in an abelian
category would induce the transformation of the weight filtration
for $H_{MM}$ to the weight filtration of $H$. Here the weight
filtration of $H$ is defined by the weight structure method, yet it
coincides with the 'classical' one (cf.  Remark
\ref{rwfiltcoh}(2)).

Therefore we obtain that our results will give a certain weight
filtration
for $H_i^{MM}(X)$ (and the corresponding Chow-weight spectral sequence)
automatically when $H_{MM}$ will be defined. Note we don't need any
information on $H_{MM}$ for this! However this construction does not
yield automatically
that  the filtration on $H_i^{MM}(X)$ obtained depends only on the
object $H_i^{MM}(X)$ and does not depend on the choice of
$X$.

Alternatively, one can obtain weights for $X\in MM$
by presenting it as $H_0^{MM}(X)$ (so we use the embedding
$MM=\hrt_{MM}\to\dmgm\q$). Then one obtains a weight
filtration for $X$ that does not depend on any choices. This filtration
certainly should coincide with the one given by the previous method;
yet in
order to prove this one needs to know that 
the weight spectral sequence $T_w(H_0^{MM},X)$ degenerates at $E_2$ (see Proposition 3.5(1) of \cite{btrans}).
The latter fact would follow from the conservativity and $t$-exactness of the \'etale realization of motives; see Corollary 3.2.4 of \cite{bmsh}.
Now, the assumptions mentioned 
are certainly expected to be
 true.
 We will say more on weights for mixed motives in
\S\ref{stfiltr} below; a detailed discussion may be found in \cite{btrans} and \cite{bmsh}.

\subsection{Motives over perfect fields of finite
characteristic}\label{motperf}

In our study of motives (here and in \cite{mymot}) we applied
several results of \cite{1} that use resolution of singularities. So
we had to assume that the characteristic of the ground field $k$ is
$0$. In \S8.3 of \cite{mymot} it was shown that using de Jong's
alterations one can extend most of our results to motives with
rational coefficients over an arbitrary perfect $k$.

In this subsection (and also in all remaining parts of this section)
we consider motives with integral coefficients over a perfect  field
$k$ of characteristic $0$. Our goal is to justify a certain claim
made in \S8.3.1 of \cite{mymot}.

In \cite{bev}  it was proved unconditionally that $\dms$ has a
differential graded enhancement. In fact, this fact can be
easily obtained by applying  Drinfeld's description of
localizations of enhanced triangulated categories.
 Moreover, Proposition 5.6 of
\cite{bev} extends the Poincare duality for Voevodsky motives to
our case. Therefore for $P,Q\in\spv$  we obtain $$\dms(\mg(P),
\mg(Q)[i])=\cho([P],[Q])\text { for }i=0;\ 0\text{ for }i>0 .$$
Hence the triangulated subcategory $DM_{pr}$ of $\dms$ generated
by $[P],\ P\in\spv$ can be described as $\tr(I)$ for a certain
negative differential graded $I$. In particular, we obtain the
existence of a conservative weight complex functor $t_0:DM_{pr}\to
K^b(\cho)$. Moreover, for any  realization of $DM_{pr}$  and any
$X\in \obj DM_{pr}$ one has the Chow-weight spectral sequence $T$.

The problem is that (to the knowledge of the author) at this
moment there is no way known to prove that $DM_{pr}$ contains the
motives of all smooth varieties (though it contains the motives of
those varieties that admit 'nice compactifications').

Instead we will prove that the weight structure on $DM_{pr}$ can
be extended to a weight structure on a larger category containing
all $\mg(X)$.

Recall that $\mg$ is a full embedding of $\dmge\supset DM_{pr}$ into
$\dme$, whereas $\dme\subset D(\ssc)$ ($\mg$ is denoted by $i$ in
Theorem 2.3.6 of \cite{1}). We denote by $D\subset D(\ssc)$ the full
category of complexes with homotopy invariant homology sheaves. We
have a full embedding $\dme\to D$.

We can extend to $D\subset D(\ssc)$ the assertion of Proposition
3.2.3 of \cite{1} (i.e., construct a projection $D(\ssc)\to D$
which is left adjoint to the embedding) using the fact that
$D(\ssc)^{t\le 0}\perp D(\ssc)^{t\ge 1}.$ Here $t$ denotes
the usual  $t$-structure of $D(\ssc)$ (corresponding to the homotopy
$t$-structure for $\dme$).
It follows that all objects of $\mg(\dmge)$ are compact. Indeed, it
is sufficient to prove this for $\mg([X])$, where  $X\in\sv$;
Proposition 3.2.3 of \cite{1} implies that $D(\mg(X),-)$ is the
corresponding hypercohomology functor which commutes with arbitrary
coproducts.

Consider $S=\{X\in \obj D:\ Y\perp X\ \forall Y\in\obj
DM_{pr}\}$. 
Note that  in the definition of $S$ it suffices to consider
$Y=\mg(P)[i],\ P\in\spv,\ i\in\z$, since $[P]$ generate $DM_{pr}$.
Obviously, $S$ is the class of objects for a certain full
triangulated subcategory of $D(\ssc)$. We denote the localization (i.e.,
the Verdier quotient)
of $D$ by $S$ by $D_S$. By definition of $S$, the set $H=\{[P],\
P\in\spv\}$ weakly generates $D_S$. Since objects of $DM_{pr}$ are
compact, $S$ is closed with respect to arbitrary coproducts. It
follows that $D_S$ admits arbitrary coproducts. Note that
$DM_{pr}\subset D_S$ by Proposition III.2.10 of \cite{gelman};
hence we have a full embedding $\chowe\to D$.


By  Theorem \ref{madts}(I2) we obtain that $D_S$ supports
adjacent weight and $t$-structures which we will call Chow ones.
By  Theorem \ref{madts}(II) we have $\hrt_{Chow}=\chowe_*$.
Moreover, $\hw$ is the category $\chowe_{\oplus}$ of arbitrary
coproducts of effective Chow motives since $\chowe$ is idempotent
complete.

Note that the definition of $w_{Chow}$ is compatible with the
definition of the Chow weight structure on $DM_{pr}$. In particular,
this reasoning extends the weight complex functor to a functor $D\to
K_\w(\chowe_{\oplus})$.
 This would
give a (possibly, infinite) weight complex for any $X\in \obj
\dmge$. Recall that (by the results of \S8.3.2 of \cite{mymot})
$t(X)$ becomes (homotopy equivalent to) a finite complex after
tensoring the coefficients by $\q$. This weight complex functor
can be 'strengthened' (see  Remark \ref{rwndg}(3)) since
$D(\ssc)$ has a differential graded enhancement.


\subsection{Coniveau and truncated cohomology}\label{coniv}

Let $k$ be an arbitrary perfect field, $H\in \obj\dme$. We denote
$\tau_{\le i}H$ by $H'$ ($i$ is fixed, $\tau$ is the
$t$-truncation with respect to the homotopy $t$-structure).
 We denote by $H''$ the 'complement of $H'$ to $H$',
i.e., $H''=\co (H'\to H)$. Note that the homology of $H''$ is
concentrated in degrees $>i$. $j\in\z$  will be a fixed integral
number up to the end of the section.

Let
\begin{equation}\label{fgmot}
M= U_u\stackrel{d_u}{\to} U_{u-1}\stackrel{d_{u-1}}{\to}\dots
U_{1}\stackrel{d_{1}}\to U_0\in\obj \dmge 
\end{equation}
be a complex in $\smc$ ($U_l$ is in degree $-l$). We demand that  for
all $r$, any (closed) point  $u\in U_{r-1}$ the codimension of the
preimage (in the sense of Remark \ref{multva})
${\operatorname{codim}}_{U_{r}}d_r^{-1}(u)\ge
{\operatorname{codim}}_{U_{r-1}}u-1$; here we define the codimension
of a subvariety as the minimum of codimensions of its parts in the
corresponding  connected components.

We fix some $j\in\z$, $\dme$ will be denoted by $\cu$. In order to
write a formula for the $H'$ and $H''$-cohomology of $M$ and prove
it we will need some notation and certain orthogonality statements.

Let $(Y^0_l,Y^1_l)$ run through open subschemes of  $U_l$ such that
$U_l\setminus Y^k_l$ is everywhere of codimension $\ge
j-i-k+1-l$ in $U_l$ ($k=0,1$, $0\le l\le u$) and
 the images (in the sense of Remark \ref{multva})
$d_l(Y^k_l)\subset d_l(Y^k_{l-1})$ for all $k,l$.  We define  the
motives $L^k=Y^k_u\to \dots \to Y^k_{1}\to Y^k_0$
 for $k=0,1$ using
the corresponding restrictions of $d_l$.
  Note that  if $Y^1_l\subset Y^0_l$ for all $l$ then
we have natural morphisms $L^1\to L^0\to M$. We denote
$N^k=\co(L^k\to M)\in \obj \dmge$.

\begin{lem}\label{lger}
1. $\cu(N^k,H''[r])= 0$ for any $r\le j-k+1$ and any $Y^k_l$.

2. $\inli \cu(L^1,H'[j+1])=\inli \cu(L^1,H'[j])=
 \inli \cu(L^0,H'[j+1])=\ns$.

\end{lem}
\begin{proof}

1. It easily seen (by cutting $H[j]$ into its $t$-homology)
 that it suffices to prove a similar statement  for $H''$
  replaced by any homotopy invariant $S\in\ssc$ shifted by $v\le j-k$.

First let all $U_l$ except $U_t$ be  empty (and all $d_l=0$). Then
our (last) assertion can be easily deduced from Lemma 4.36 of
\cite{3} (and some cohomological comparison results of Voevodsky)
using the standard coniveau spectral sequence argument. We write
down a (short) proof here. The classical coniveau spectral sequence
is described (for example) in  \cite{suger} as follows: for any $U\in\sv$ by
formula (1.2) of ibid. there exists a spectral sequence
$$ E_1^{p,q}=\bigoplus_{x\in U^{(p)}}H^{p+q}_{x,Zar}(U,S)\implies
H_{Zar}^{p+q}(U,S);$$ here $U^{(p)}$ denotes the set of points of
$U$ of codimension $p$, $H^{p+q}_{x,Zar}(U,S)$ is the local Zariski
cohomology group (see \S4.6 of \cite{3} and Lemma 1.2.1 of
\cite{suger}). Note that this spectral sequence is functorial with
respect to open embeddings. Now, Lemma 4.36 of \cite{3} yields
cohomological purity in this case; in particular,
$H^{p+q}_{x,Zar}(U,S)=0$ for any $x\in U^{(p)}$ unless $q=0$. It
follows that the map $H_{Zar}^{v}(U_t,S)\to H_{Zar}^{v}(Y^k_t,S)$ is
bijective for $v<j-k-t$ and injective for $v=j-k-t$. Next, the
Zariski cohomology of $S$ coincides with its Nisnevich cohomology by
Theorem 5.3 of \cite{3}, whereas the latter equals $\cu(-,S[v])$ by
Proposition 3.2.3 of \cite{1}. Hence the long exact sequence of
relative cohomology for ($Y^k_t,Y_t$) yields our (last) claim
 in this case.

Now let $u=1$. We have an exact sequence
$$ 
 \ns= \cu(Y^k_{1}\to U_1,S[v-1])  \to
 \cu(N^k,S[v]) 
 \to \cu(Y^k_0\to U_0,S[v]) =\ns
$$ for $k=0,1$, $v\le j-k$; this yields the claim in this case. The
case  $u>1$ can be easily obtained from similar exact sequences
by induction.

2. The proof is similar to that of assertion 1. One should cut $H'$
into its $t$-pieces and apply the coniveau spectral sequence arguments.

To this end we recall that the inductive limit of (long) exact
sequences is exact, so we can pass to the limit in the coniveau
spectral sequence. Besides, the codimension condition (on $U_l$)
implies that for sets of $Y^k_l$ as in the assertion all (single)
$Y^k_l$ may be 'as small as possible'. This means that for any
open $Y\subset U_l$ such that $U_l\setminus Y$ is everywhere of
codimension $\ge j-i-k+1-l$ in $U_l$ ($k=0,1$) can be completed
to some set of $Y^k_l$; besides, we can intersect such sets
(componentwisely). It follows that the corresponding
$\inli_{Y_l^k} \bigoplus_{x\in
Y_l^k{}^{(p)}}H^{p}_{x,Zar}(U,S)=\ns$ since $\prli _{Y_l^k}
Y_l^k{}^{(p)}=\emptyset$ (for the  values of $p$ corresponding to
our situation).

\end{proof}

\begin{theo}\label{trcoh}  

I 
We have an isomorphism
\begin{equation}\label{ffbws} \cu(M,H''[j])\cong
\imm(\inli \cu(L^0,H[j])\to \inli \cu(L^1,H[j])).\end{equation} 

Here the connecting morphisms between the cohomology of $L^k$ for various sets
$(Y^k_l)$ are induced by open embeddings of varieties.

II For any $j$ we have
\begin{equation}\label{fbws} \cu(M,H'[j])\cong
 \imm(\inli \cu(N^0,H[j])\to
\inli \cu(N^1,H[j]));\end{equation}
the limit is defined as in assertion I.

III 1. The isomorphisms described above are functorial in the
obvious way with respect to 'nice' morphisms of complexes of
 correspondences $(f_l):M'\to M$. Here $M_l$ is a complex of
$U_l'$, $(f_l)$ is nice if for any   $l$, for any (closed) point
$u\in U_{l}$ we have ${\operatorname{codim}}_{U_{l}'}f_l^{-1}(u)\ge
{\operatorname{codim}}_{U_{l}}u$ (in the sense of Remark
\ref{multva}).

2. Furthermore, suppose that for some $(f_l)$ and fixed set of
$Y^0_l\subset U_l$ (satisfying the above conditions) we have
${\operatorname{codim}}_{U_{l}'}f_l^{-1}(U_l\setminus Y^0_l)\ge
j-i-l$. Then the morphism $f_{H''}^*:H''(M)\to H''(M')$ (resp.
$f_{H'}^*:H'(M)\to H'(M')$) is compatible with the natural morphism
$ \cu(L^0,H[j])\to \cu(L^1{'},H[j]))$ (resp. $ \cu(N^0,H[j])\to
\cu(N^1{'},H[j]))$) via the isomorphism of assertion I (resp.
assertion II).

\end{theo}
\begin{proof}

I Shifting $H$ we 
easily reduce the statement to the case $i=0$.

Lemma \ref{lger} allows us to argue similarly to the proof of
 Theorem \ref{sdt}(7); note that our assertion is an
analogue of part 8 of loc. cit.

Part 1 of the lemma yields  exact sequences
\begin{equation} \begin{aligned}\label{ger1} \ns=
\inli \cu(N^0,H''[j]) \to \cu(M,H''[j])\\
\to \inli \cu(L^0,H''[j]) \to \inli \cu(N^0,H''[j+1])=\ns
 \end{aligned}\end{equation}
and  \begin{equation} \label{ger2} \ns=\inli \cu(N^1,H''[j]) \to
\cu(M,H''[j])\to
\inli \cu(L^1,H''[j])  
\to\dots\end{equation}

By  Lemma \ref{lger}(2) we also have $\inli \cu(L^1,H[j])
\cong \inli \cu(L^1,H''[j])$.

Now we consider the commutative diagram 
\begin{equation}\label{badcd}\begin{CD}
@. \inli \cu(L^0,H[j])@>{}>>\inli  \cu(L^1,H[j])\\@.
@VV{g}V@VV{h}V \\
\cu(M,H''[j])@>{t}>> \inli \cu(L^0,H''[j]) @>{p}>>\inli
 \cu(L^1,H''[j])
\end{CD}\end{equation}

We have proved that $t$ and $h$ are bijective, $g$ is surjective,
and $p$ is injective. This immediately yields the assertion
required.

II 
Using Lemma \ref{lger}, we can argue exactly as in the proof of
assertion 7 of Theorem \ref{sdt} (and dually to the reasoning
above).

III1. We describe the functoriality in question for 'nice' $(f_l)$.

For  $r=0$ or $1$ let $Y^r_l$ be fixed for all $l$. Then we can take
$Y_l^r{'}=U_l'\setminus f_l^{-1}(U_l\setminus Y^r_l)$; $(f_i)$
induces morphisms $L{'}^r\to L^r$ and $N{'}^r\to N^r$. It remains to
note that the proofs of assertions I and II are compatible with these
morphisms.

III2. We check the compatibility desired for $H''$; the statement
for $H'$ can be proved similarly (and dually  in the categorical
sense).

We should verify the following. Let  $v\in \cu(M,H''[j])$ come from
some $w\in \cu(L^0,H[j])$ (for our fixed $Y_l^0$) via (\ref{ffbws}).
Denote by (\ref{ffbws}') the isomorphism (\ref{ffbws}) with $M$
replaced by $M'$, $L^r$ replaced by $L'^r$. Then $f^*_{H''}(v)$
should be mapped via (\ref{ffbws}') to the image of $f^*_{H}(w)$ in
$\inli \cu(M'^1,H[j])$. Here $f^*_{H}$ is the map
$\cu(L^0,H[j])\to \cu(L',H[j])$ induced by $(f_l)$, and $L'$ is the
complex of $U_l'\setminus f_l^{-1}(U_l\setminus Y^0_l)$.

The latter fact follows 
easily from the commutativity of the diagrams
$$
\begin{CD}
\cu(M,H''[j])@>{}>>  \cu(L^0,H''[j]) @>{}>> \inli \cu(L^0,H''[j])
\\@VV{}V
@VV{}V \\ \cu(L',H''[j])@>{}>>  \cu(L',H''[j]) @>{}>> \inli
\cu(L'{}^1,H''[j])
\end{CD}$$
and
$$
\begin{CD}
\cu(L^0,H[j])@>{}>>  \cu(L^0,H''[j])
\\@VV{}V
@VV{}V \\ \cu(L',H[j])@>{}>>  \cu(L',H''[j])
\end{CD}$$

\end{proof}

\begin{rema}\label{rws}

1. The main difference of this result from the usual comparison of
spectral sequences (as in \cite{paran}) is as follows: we calculate the
$D$-terms of the corresponding exact couple instead of  $E$-ones;
we compute cohomology of certain motives (instead of varieties as in
\cite{blog} and \cite{suger}).

2. Instead of applying assertion III2 to a single 'nice' set of
$Y^0_l$ one may consider a (directed) system of those.
 This is especially actual  if  the right hand sides of
(\ref{ffbws}) or (\ref{fbws}) can be calculated using such a
'nice' directed subset of the set of all possible $(Y^0_l)$ (which
is often the case). In this case assertion III2 allow 'calculating'
$f_{H''}^*$ (resp. $f_{H'}^*$) completely.

3. One can  generalize (\ref{fbws}) in the following way. Let
$r\ge i$; denote  $\tau_{\le r}H$ by $G$. Then for the corresponding
morphism of cohomology theories $ H'\to G$ we have
\begin{equation}\label{febws} \imm(G^j(M))\to H'[j](M))\cong
\inli \imm(H^j(N^0))\to H^j(N^{r-i+1})).\end{equation}

Here $N^{r-i+1}$ is defined similarly to $N^0,N^1$ in the theorem.
 This
statement can be easily obtained by calculating the $D$-terms
of the higher derived couples for the coniveau spectral sequence;
cf. Theorem \ref{sdt}.

4. Instead of considering limits of cohomology of motives we could
have considered the cohomology of the corresponding {\it pro-motives} as
it was done in \S4 of \cite{deg2}; this wouldn't have affected the
proof substantially. Unfortunately, the category of pro-motives is
not triangulated (if we define it in the obvious way).

A certain triangulated analogue of pro-motives
(a category of {\it comotives}) was
constructed in \cite{bger} (see \S1.5, \S3.1, and \S5 of ibid.). It
contains more information, yet is somewhat more difficult to deal
with (in particular, the author currently does not know how to
control the corresponding 'homotopy limits' of motives unless $k$
is countable). In this category the general weight structure
formalism can be applied directly; this allows one to get rid off the
codimension condition for $d_l$. Unfortunately, it seems difficult
to describe the corresponding weight decompositions of motives and their
morphisms in
the general case; one only knows that they can be constructed
from coniveau filtration of $U_r$ (as in Lemma \ref{twd}).


5. As was noted in (part 6 of) Remark \ref{rsdt}, all these
statements can be vastly generalized.

\end{rema}

\subsection{Expressing torsion motivic cohomology
(with compact support) in terms of \'etale one}\label{motet}

For  fixed $n>0$, $(l,p)=1$, $r\ge 0$, we  denote by
$\het(r)\in \obj \dmk$ some \'etale resolution of
$\mu_{l^n}^{\otimes r}$ by injective \'etale sheaves with transfers.
This object does not depend on the choice of a resolution (as an object
of the derived category) by obvious
reasons (it is the total derived image of $\mu_{l^n}^{\otimes r}$
with respect to the corresponding change of topologies functor).
$\het(r)$ is homotopy invariant, so we can  substitute it for $H$ in
the statements above (since for any fixed $M$ it suffices to
consider some homotopy $t$-truncation of
$\het(r)$, whereas the latter belong to $\obj\dme$).

\begin{pr}[The Beilinson-Lichtenbaum Conjecture]\label{pblc}\label{blc}

The  (well-known) cycle class map $\z/l^n\z(r)\to \het(r)$
identifies the former object with $\tau_{\le r}\het(r)$.

\end{pr}

 We recall that this statement is equivalent to the
Bloch-Kato Conjecture (see \cite{suvo} and \cite{gelev}). The
latter is well-known for $l=2$ (see \cite{vomi}), and was recently 
proved for an arbitrary $l$ (see \cite{vobe}).

For a motif $X$ we denote  $\dme(X,\z/l^n\z(s)[i])$ by
$\hml^{i}(X,\z/l^n\z(s))$;
$\hel^{i}(X,\z/l^n\z(s))=\dmk(X,\het(s)[i])$.

Now, Theorem \ref{trcoh} easily yields the following statement
(in the notation of loc. cit.). 

\begin{coro}\label{cmot}

1. For $M$ as in (\ref{fgmot}),  we have
\begin{equation}\label{eqmotcoh} H^{j}(M,\z/l^n\z(s))\cong
\imm(\inli \hel^{j}(N^0,\z/l^n\z(s))\to \inli
\hel^{j}(N^1,\z/l^n\z(s))).\end{equation} 

The corresponding functor $H_s''$ can be calculated as follows:
$$H_s''{}^{j}(M)\cong \imm(\inli \hel^{j}(L^0,\z/l^n\z(s))\to
\inli \hel^{j}(L^1,\z/l^n\z(s))).$$

These isomorphisms satisfy those functoriality properties that were
described in  Theorem \ref{trcoh}(III).

2. Let $Y_h,\ 1\le h\le u$, $u>0$, be smooth of the same dimension;
let $Y=\cup Y_h$  be a normal crossing scheme, i.e., all intersections of the components
(in some large basic scheme) are normal and smooth. Consider the
motif $M$ corresponding to the complex $(U_l)$; here
$U_l=\sqcup_{(i_j)}Y_{i_1}\cap Y_{i_2}\cap\dots \cap Y_{i_{l+1}}$
for all $1\le i_1\le\dots \le i_{r+1}\le u$, $d_l$ is the alternated
sum of $l+1$ natural maps $U_l\to U_{l-1}$.

Let  $(Y^0,Y^1)$ run through open subschemes of $Y$ such that
$Y\setminus Y_r$ is (everywhere) of codimension $\ge j-r-s+1$ in $Y$
($r=0,1$). Then we have
\begin{equation}\label{eqmotdef} H_s''{}^{j}(M)\cong
\imm(\inli \hel^{j}(Y^0,\z/l^n\z(s))\to \inli
\hel^{j}(Y^1,\z/l^n\z(s))).\end{equation} For $N^r=\mg(Y^r\to Y)$
($Y$ is in degree $0$) we have \begin{equation}\label{eqmotcohn} 
H^{j}(M,\z/l^n\z(s))  \cong \imm(\inli \hel^{j}(N^0,\z/l^n\z(s))\to
\inli (\hel^{j}(N^1,\z/l^n\z(s))) .\end{equation}

3. Let $Y'=\cup Y_i'$, $M'$ is defined similarly to $M$, let
$f:Y'\to Y$ be a morphism of schemes, suppose that for any (closed)
point $u\in Y$ we have ${\operatorname{codim}}_{Y'}f^{-1}(u)\ge
{\operatorname{codim}}_{Y}u-1$. Then the morphisms $f^*_{H'_s}$ and
$f^*_{H''_s}$ can be computed by the  way  described  in  Theorem \ref{trcoh}(III2) (see also its proof and  Remark
\ref{rws}(2)).

4. Suppose that  $U\in \sv$ equals $P\setminus Y=\cup Y_i, 1\le i\le m,$ where
$P \in \spv$, $Y$ is a smooth normal crossing divisor. Then for any
cohomological functor $G$ defined on $\dme$ we have a long exact
sequence \begin{equation}\label{exa}\dots G^j(\mg^c(U))\to G^j(X)\to
G^j(M) \to \dots\end{equation}
\end{coro}
\begin{proof}
1. This is immediate from  Theorem \ref{trcoh}(2) applied for
$H=\het(s)$, $i=s$.

2. We should prove that the formulas (\ref{eqmotdef}) and
(\ref{eqmotcohn}) compute the limits described in 
Theorem \ref{trcoh}(1,2).

First, we note that if $Y^r$ ($r=0,1$) satisfies the condition of
the assertion then 
\begin{equation}\label{inters} Y^r_l=\sqcup_{(i_j)}Y^r\cap Y_{i_1}\cap Y_{i_2}\cap\dots \cap
Y_{i_{l+1}}\end{equation} satisfies the conditions of Theorem \ref{trcoh}(1). Next, using proper descent we easily obtain
that the \'etale cohomology of $Y^r$ is isomorphic to that of
$L^r=(Y^r_l)$. It suffices to note that  any set of $Y^r_l$ as in
 Theorem \ref{trcoh}(1) can be shrunk to a one coming from
some $Y^r$ as in (\ref{inters}).

3. It suffices to note that the functoriality provided by  Theorem \ref{trcoh}(III2) is compatible with that of the
formulas (\ref{eqmotdef}) and (\ref{eqmotcohn}).

4. By definition,
$\hml_c^{j}(U,\z/l^n\z(s))=\hml^{j}(\mg^c(U),\z/l^n\z(s))$. Hence it
suffices to recall that $\mg^c(U)\cong \co(M\to \mg(X))$.  The
latter fact in the characteristic $0$ case  is Proposition 6.5.1 of
\cite{mymot}. In the characteristic $p$ case one can deduce the
statement from the results of \cite{de2} (and the Poincare duality
properties).

\end{proof}

\begin{rema}

1. Note that for $G=H'$ or $G=H''$ one can compute the map
$G^*(M)\to G^*(Y)$ in (\ref{exa}) using part 3 of the corollary.

Besides, one can write down the formula for the (motivic and
$H''$-) cohomology of $\mg^c(U)$ by substituting the 'complex' $Y\to
X$ for $M$ into part 1 of the corollary; here one should ignore the
fact that $Y$ may be singular.

2. Suppose that $K$ contains a primitive $l^n$-th root of unity. Then using it
one can identify all $\het(r)$, and so obtain certain maps
$\hml^{i}(-,\z/l^n\z(s)) \to \hml^{i+j}(-,\z/l^n\z(s+j))$ (induced
by the multiplication on the corresponding motivic Bott elements, as
in \cite{levbot}). Then (\ref{febws})  allows us to calculate the image
of these maps.

One can prove natural analogues of part 2 of the corollary and
part 1 of this remark.

3. It seems very interesting to replace \'etale cohomology in the
right hand side of (\ref{eqmotcohn}) by singular cohomology (in the
case when $k$ is the field of complex numbers). One can easily deduce
from Corollary 7.4 of
\cite{blog} that in the case $u=0$, $j=2s$, the formula would
calculate the group of algebraic cycles in $U_0$ of codimension
$s$ modulo algebraic equivalence. So it seems that the  homotopy
$t$-structure truncations of singular cohomology should be related
to a certain (non-existent yet) 'theory of mixed motives up to
algebraic equivalence'; see the end of \cite{voenilp}. 

4. Certainly, the cohomology of $M$ is a very natural candidate for
the cohomology of $Y$; note that $\cup Y_j\to Y$ is a $cdh$-covering
(see \cite{4}). Yet this does not automatically imply the
isomorphism of cohomology for all 'reasonable' cohomology
theories.

5. Recall that if $k$ admits resolution of singularities any smooth
quasi-projective $U$ can be presented as $X\setminus\cup Y_i$.

\end{rema}

\subsection{The cases when
$t_{\chow}$ can be easily calculated; relation to unramified
cohomology} \label{sunram}  

By Theorem \ref{sdt}, one can express the values of
$\dme(X,\tcho^{\le i}(Y)[j])$ (for example, in the case
$X\in\obj\dmge,\ Y\in \obj\dme$)
 in terms of the $Y$-cohomology of weight truncations of $X$.
It turns out that in some cases one obtains
very nice results this way; these are related to unramified
cohomology (see Remark \ref{runram}) below.

\begin{pr}\label{punram}

Let $U$ belong to $\sv$; let $P$ be any smooth compactification of $U$
(i.e., $U$ is open in $P$,
$P$ is smooth proper).

1. Suppose that $\cha k=0$. Then for any homotopy invariant
 $S\in \ssc$ we have $\cu (\mg(U),S^{\tcho=0})\cong S(P)$.

2. The same is true for any perfect $k$ if $S$ is a sheaf of
$\q$-vector spaces.


\end{pr}

\begin{proof}

For the proof of assertion 1 we take $\cu=\dme$; $\cu=\dmeq$ for assertion 2. We will use the same notation for motives of
varieties in these categories, and for their  weight and
$t$-structures.

We have morphisms $\mg(U)\stackrel{g}{\to} w_{\le 1}\mg(U)\stackrel
{f}{\to} w_{\le 0}\mg(U)$. By  Theorem \ref{sdt}(8), we have
$\cu(\mg(U), \tcho^{\le 0}(S))= \imm f^*(S)$.

Now, by Lemma \ref{lunram} below, we can assume that
$w_{\le 0}(\mg(U))=\mg(P)$.
 For this choice $(f\circ g)^*_S$ is injective by Corollary 4.19
of \cite{3}. Hence $f^*_S$ is injective also, and we get $S(P)\cong
\tcho^{\le 0}(S))$.

Lastly, we note that $S$ belongs to $\cu^{\tcho\le 0}$
(since smooth proper
varieties have no negative $S$-cohomology; note that the
construction of $\tcho$ uses Theorem \ref{madts}).
 Hence  $\tcho^{\le 0}(S)=S^{\tcho=0}$.

\end{proof}

\begin{rema}\label{runram}

1. The statement proved immediately yields that $S(P)$
is a birational invariant of $P$ (i.e., it depends only on the
function field of $P$); cf. Theorem 8.5.1 of \cite{suger}.

Now we relate the statement proved above to unramified cohomology (as
defined in \S4.1 of \cite{colbir}). Suppose that a cohomology
theory can be represented by $C\in \obj \cu$ (or by an object of  the unbounded
version of $\dme$); for example, this is (essentially) the
case for torsion
\'etale cohomology and de Rham cohomology. Then
for $U\in \sv$ the $i$-th unramified cohomology of $k(U)$ (the
function field) with coefficients in $C$ equals $S^i(P)$, where
$S^i=C^{t=i}$, $P$ is a smooth compactification of $U$. A similar
statement was verified in  Theorem 4.1.1 of \cite{colbir}; note
that  the results of  \cite{3}, \S4, yield all properties of
$\cu(-,C)$ that are necessary for the proof.

2. Now, let $C\in \cu^{t\ge 0}$ (with rational coefficients if
$\cha k>0$). Certainly, one has $\cu(\mg(U),\tcho^{\le
0}C)=\cu(\mg(P),C)$. Indeed, it suffices to apply the Proposition
above for $S=H^{t=0}(C)$.

More generally, by analyzing  weight
decompositions of $\mg(U)$ in more detail (and fixing their choice), one can check that for any
$i\ge 0$, $C\in \cu^{t\ge -2i}$ the map $\cu (w_{\le i}\mg(U),C)\to
\cu(w_{\le i+1}\mg(U),C)$ is injective. Therefore, we have $\cu(\mg(U),
C^{\tcho\ge -i})=\cu (w_{\le i}\mg(U),C)$ (for a 'nice' choice of
$\cu (w_{\le i}\mg(U)$).

Moreover, one can replace the functor $\cu (-,C)$ by any
cohomological  $H:\cu\to A$ ($A$ is an abelian category; in the case
$\cha k>0$ it should be $\q$-linear) such that
$H(\mg(U)[i])=0$ for any $U\in \sv$, $i<0$. Then the statements
above (of part 2 of this Remark)  will also be true if one defines
$H^{\tcho\ge -i}$ using virtual $t$-truncations (see Remark \ref{rtrfun}).

3. As mentioned in Remark \ref{rperf}, the author hopes to extend
the results of this section to the setting of
$\z[1/p]$-coefficients (in the case $\operatorname{char}k=p$).

\end{rema}

Now we prove the statement that was used in the proof of
Proposition \ref{punram}

\begin{lem}\label{lunram}

In the conditions of Proposition \ref{punram}  the map $\mg(U)\to
\mg(P)$ can be extended to a weight decomposition of $\mg(U)$.

\end{lem}

\begin{proof}

We should prove that $\co(\mg(U)\to\mg(P))\in \cu^{w\ge 1}$. In
the case $U=P\setminus Z$, $Z\subset P$ is smooth projective,
 this statement is immediate from the Gysin distinguished triangle
 $\mg(U)\to \mg(P)\to \mg(Z_i)(r_i)[2r_i]$,
 $Z_i$ are connected components of $Z$, $r_i$ are their
codimensions (see
 Proposition 3.5.4 of \cite{1} for $\cha k=0$ and Proposition
5.21 of \cite{deg2} for the general case).

In the general case the statement is also proved   by induction.
We choose a stratification of $P\setminus U=\cup Z_i$, 
where
$Z_i\setminus Z_{i+1}\in \sv$, $Z_m=\ns$ for some $m$. Then we can
apply the Gysin triangle for the pair ($P\setminus Z_{i+1},\
Z_i\setminus Z_{i+1})$. Part 1 of Remark \ref{mgc} yields that each
of $\mg(Z_i\setminus Z_{i+1})$ belongs to $\cu^{w\ge 0}$. Note
here that the latter statement is true for motives with rational
coefficients in any characteristic, since alterations yield that
the motif of  any smooth variety can be presented as a retract
of a motif of a complement of a smooth normal crossing
divisor; see Appendix B of \cite{hubka}.

\end{proof}


\section{Supplements}\label{ssupl}

We start the section by proving that weight structures have nice
'functorial' properties similar to those of $t$-structures (yet the
difference is substantial).

In \S\ref{swesloc} we show that a weight structure $w$ on $\cu$
which induces a weight structure on a triangulated $\du\subset\cu$
yields also a weight structure on the localization $\cu/\du$, and calculate $H(\cu/\du)$. 

In \S\ref{sglu} we prove the following converse to this statement:
 weight structures can be 'glued' in
a manner similar to that for $t$-structures. The author
 applied this fact for the construction of weight
structures
for relative motives 
 in \S2.3 of \cite{brelmot} and in \cite{bonivan}. 

 In \S\ref{scompt} we study the interaction of weight and
$t$-truncations. 
 In
\S\ref{sftriang} we prove (using an argument due to A. Beilinson)
that any $f$-category enhancement of $\cu$ yields a lift of $t$ to
a 'strong' weight complex functor $\cu\to K(\hw)$; cf. Remark
\ref{rctst}.

  In
\S\ref{svarwei} we discuss  other possible sources of conservative
'weight complex-like' functors and related spectral sequences.

We conclude the section by a discussion of three relevant types of
filtrations for triangulated categories; all of them 'should' be
actual (and closely related) for $\dmge\q$.

\subsection{Weight structures in localizations}\label{swesloc}

We call a category $\frac A B$ the {\it factor} of an additive
category $A$
by its full additive subcategory $B$ if $\obj \bl \frac A B\br=\obj
A$ and $(\frac A B)(X,Y)= A(X,Y)/(\sum_{Z\in \obj B} A(Z,Y) \circ
A(X,Z))$.

\begin{pr}\label{ploc}

1. Let $\du\subset \cu$ be a 
triangulated subcategory of
$\cu$; suppose
 that $w$ induces a weight structure on $\du$
(i.e., $\obj \du\cap \cu^{w\le
 0}$ and $\obj \du\cap \cu^{w\ge
 0}$ give a weight structure for $\du$).
We denote the heart of the latter weight structure  by $HD$.

 Then $w$ induces a weight structure on
 $\cu/\du$ (the localization, i.e., the Verdier quotient of $\cu$
by $\du$).
 This means
 that the Karoubi-closures of $\cu^{w\le
 0}$ and $\cu^{w\ge
 0}$ (in $\cu/\du$) give a weight structure for $\cu/\du$
(note that $\obj \cu=\obj \cu/\du$).

2. 
The heart $H(\cu/\du)$ of this weight structure 
is the Karoubi-closure  of $\frac { \hw} {HD}$ in
$\cu/\du$.

3. If $\cu,w$ is bounded (above, below, or both), then $\cu/\du$
also is.

\end{pr}
\begin{proof}
1. It clearly suffices to prove that for any $X\in (\cu/\du)^{w\ge
 0}$ and $Y\in (\cu/\du)^{w\le
 -1}$ we have $\cu/\du(X,Y)=\ns$; all other axioms of Definition
\ref{dwstr}
 are fulfilled automatically
 since $\cu/\du$ is a localization of $\cu$.

Recall now (see Lemma III.2.8 of \cite{gelman}) that any morphism in
$\cu/\du(X,Y)$ can be presented as $fs\ob$, where $f\in \cu(T,Y)$
for some $T\in\obj \cu$, $s\in \cu (T,X)$,  $\co (s)=Z\in\obj \du$.

By our assertion, there exists a choice of $Z^{w\ge 0}$ that belongs
to $\obj\du$. Since $X\perp w_{\le -1}Z $ we can factor the
morphism $X\to Z$ (corresponding to $s$) through $Z^{w\ge 0}$.

Hence (applying the octahedral axiom) we obtain that there exist
$T'\in \obj \cu$ and a morphism $d:T'\to T$ such that $\co d=w_{\le
-1}Z\in \obj \du$, whereas a cone of the composed morphism
$s':T'\to X$ equals $Z^{w\ge 0}$. It follows that $fs\ob
=(f\circ d)s'{}\ob$ in $\cu/\du$. Now note that $T'\in \cu^{w\ge 0}$ by
 Proposition \ref{bw}(3). Hence $T'\perp Y$, which yields
$f\circ d=0$. 


2. By construction, $\cu^{w=0}\subset (\cu/\du)^{w=0}$.

Now we prove that any object of $H(\cu/\du)$ is a retract of an
object of $\hw$ (in $\cu/\du$).

Let $Z$ belong to $ (\cu/\du)^{w=0}\subset \obj \cu$. We consider a distinguished triangle $w_{\ge 1}(Z)\to Z\to w_{\le 0}Z$ (a rotation of a weight decomposition of $Z$) in $\cu$.
In $\cu/\du$ we have $w_{\ge 1}(Z)\in (\cu/\du)^{w\ge 1}$, hence
$\cu/\du(w_{\ge 1}(Z), Z)=\ns$. Therefore $Z$ in $\cu/\du$ is a
retract of $Z^{w\le 0}$. Moreover,  $Z^{w\ge 1}\in (\cu/\du)^{w=0}$
since it is a retract of $Z^{w\le 0}\in (\cu/\du)^{w\le 0}$;
 therefore $Z^{w\le 0}\in (\cu/\du)^{w=0}$. Now applying the dual
 argument to $Z^{w\le 0}$ (see Remark \ref{rdd}), we obtain that
$Z$ in $\cu/\du$ is a retract of
 some $Z^0\in \obj \cu^{w=0}$. 

To conclude the proof it suffice to check that the natural functor
 $i:\hw/HD\to H(\cu/\du)$ is a full embedding. We
 consider the composition $c: \cu \stackrel{t}{\to} K_\w(\hw)\to
 K_\w(\hw/HD)$. This is a weakly exact functor that maps  all objects of $\du$ to $0$.
 Hence it maps all morphism whose cones belong to $\obj D$ into invertible ones; therefore it factors through the localization $\cu\to \cu/\du$. Since $c$ only kills those morphisms in $\hw$ that factor through $HD$, $i$ is injective on morphisms.


 It remains to prove that for $X,Y\in \cu^{w=0}$ any morphism $g:X\to Y$ in $\cu/\du$ comes
 from $\cu(X,Y)$. Applying 
 the same argument as in the proof of
 assertion 1, we obtain that $g$ can be presented as $fs\ob$
where $f\in
\cu(T,Y)$ for some $T\in\obj \cu$, $s\in \cu (T,X)$,  $\co s=Z\in
\du^{w\ge 0}$. Then $\cu(X,Y)$ surjects onto $\cu(T,Y)$. Now the
'calculus of fractions' yields the result. 

3. Since $\obj \cu/\du=\obj \cu$, we obtain the claim.

\end{proof}

\begin{coro}\label{loccor}
Let  $E\subset \hw$ be an additive subcategory. If  $X $ belongs to
the Karoubi-closure $ \obj \lan E\ra$, then  $t(X)$ is a retract of
some object of  $ K_\w^b(E)$ (here we assume that $ K_\w^b(E)\subset
K_\w(HC)$).

If $(\cu,w)$ is bounded then the converse implication also holds.

\end{coro}
\begin{proof}
We can assume that  $X\in \obj \lan E\ra$.  Then $X$ can be
obtained from objects of $E$ by repetitive consideration of cones of
morphisms. Since $t(\obj E)\subset \obj K_\w(E)$ and $t$ is a weakly
exact functor in the sense of Definition \ref{dkw},
 we obtain that
$t(X)\in \obj K_\w^b(E)$.

Conversely, let $t(X)$ be a retract of $Y\in \obj K_\w^b(E)\subset
\obj K_\w^b(\hw)$. By Proposition \ref{ploc} we obtain that $\cu/\lan
E\ra$ possesses a
 bounded weight structure whose heart contains $\frac{\hw}{E}$ as
a full subcategory.
  Hence, by Theorem \ref{wecomp}(V)
we obtain that $t_{\cu/\lan E\ra }$ is conservative. $Y\in \obj
K_\w^b(E)$ gives $ t_{\cu/\lan E\ra}(Y)=0$, hence $X $ and $Y$
belong to the Karoubi-closure of  $\lan E\ra$.

\end{proof}

\begin{rema}

1. Note that (in general) one cannot be sure that the
'factor weight structure' on $\cu/\du$ is non-degenerate.

2. Corollary \ref{loccor} is parallel to  Proposition 8.2.1(3)
of \cite{mymot}. In particular, it may be used for proving that the
motif of a smooth variety is mixed Tate if and only if its weight complex
(defined in \cite{gs}) is (this is Corollary 8.2.2 of \cite{mymot}).

3. Adding certain additional restrictions, one may also formulate
a criterion for $t(X) $ to belong to the Karoubi-closure of  $ \obj
K_\w(E)$ (instead of $\obj K_\w^b(E)$).

4. One can (easily) apply Proposition \ref{ploc} for the calculation
of $\homm(\tilde{M}(X), \tilde{M}(Y)[i])$ for $i\ge
0$; see Corollary 7.9 of \cite{kabir} (note that this
statement is not quite true for $i<0$). Here $X,Y\in \spv$,
$\tilde{M}(X), \tilde{M}(Y)$ are their birational
motives considered as objects of the triangulated
 category of birational motives (see \S5 of \cite{kabir}). 

5. 
Certainly, under the assumptions of Proposition \ref{ploc}
the inclusion $\du\to \cu$ and the projection $\cu\to \cu/\du$
are weight-exact (in the sense of Definition \ref{dadad}).

6. Note
that the description of $H(\cu/\du)$
is quite distinct from the 
description of the heart of the $t$-structure for a '$t$-exact localization'.
This is no surprise in view of  Theorem \ref{sdt}(3).

\end{rema}

\subsection{Gluing weight structures}\label{sglu}

Since weight structures are often 'dual' to $t$-structures (see
\S\ref{sadst}), this is no surprise that one can modify the 'gluing'
procedure of \S1.4 of \cite{BBD} (for $t$-structures) so that it can
be applied to weight structures (yet cf.  Remark \ref{rglu}(4)
below).

Our result is a certain converse to the result of the previous
subsection;
we study when a weight structure for $\cu$
can be recovered from  weight structures on a triangulated
$\du\subset \cu$ and on the
localization of $\cu$ by $\du$. As in the similar situation for
 $t$-structures, we need
{\it gluing data}, i.e., certain adjoint functors should exist.

We describe  gluing data for abstract triangulated
categories as it was done in \S1.4.3 of \cite{BBD} (we will only
change the notation for categories; see also the exercises at the
end of \S IV.4 of \cite{gelman}). 
Still recall that usually gluing
data sets come from certain derived categories of sheaves; this also
explains our notation for functors (yet note that we are actually
interested in the derived versions of the corresponding functors on
categories of sheaves). In particular, the category $\du$ (below)
usually comes from a closed subspace of the space corresponding to
$\cu$, whereas $\eu$ comes from its (open) complement.

It is well known (see Chapter 9 of \cite{neebook}) that
gluing data can be uniquely recovered from an inclusion
$\du\to \cu$ of triangulated categories that admits both a left
and a right adjoint functor. Then for $\eu$ being the Verdier
quotient
 $\cu/\du$ the projection $\cu\to \eu$ also admits both a
left and a right adjoint. Following \cite{BBD},
we summarize the properties of the functors obtained
(and introduce notation for them).

\begin{defi}\label{dglu}
The set $(\cu,\du,\eu, i_*, j^*, i^*, i^!, j_!, j_*)$  is called
{\it gluing data} if it satisfies the following conditions.

(i) $\cu,\du,\eu$ are triangulated categories; $i_*:\du\to\cu$,
$j^*:\cu\to\eu$, $i^*:\cu\to\du$, $i^!:\cu\to\du$, $j_*:\eu\to\cu$,
$j_!:\eu\to\cu$ are exact functors.

(ii) $i^*$ (resp. $i^!$) is left (resp. right) adjoint to $i_*$;
$j_!$ (resp. $j_*$) is left (resp. right) adjoint to $j^*$.

(iii)  $i_*$ is a full embedding; $j^*$ is isomorphic to the
localization (functor) of $\cu$ by
$i_*(\du)$.

(iv) For any $X\in \obj \cu$ the pairs of morphisms $j_!j^*X \to X
\to i_*i^*X$ and $i_*i^!X \to X \to j_*j^*X$ can be completed to
distinguished triangles (here the connecting
morphisms come from the adjunctions of (ii)).


(v) $i^*j_!=0$; $i^!j_*=0$.

(vi) All of the adjunction transformations $i^*i_*\to \id_{\du}\to
      i^!i_*$ and $j^*j_*\to \id_{\eu}\to j^*j_!$ are isomorphisms of
      functors.

\end{defi}

Now suppose that both $\du$ and $\eu$ are endowed with weight structures,
that we will (by an abuse of notation) denote by $w$.
We prove that there exists a (unique) weight structure for $\cu$ such that $i_*$ and $j^*$ are weight-exact (with respect
to weight structures mentioned; see Definition \ref{dadad}).
To this end
we consider
$\cu^{w\le 0}=\{X\in \obj
\cu:\ i^!X\in \du^{w\le 0} ,\ j^*X\in \eu^{w\le 0} \}$ and
$\cu^{w\ge 0}=\{X\in \obj \cu:\ i^*X\in \du^{w\ge 0} ,\ j^*X\in
\eu^{w\ge 0} \}$.

Before proving that these classes actually define a weight
structure, we will study how they behave with respect to the connecting functors.

\begin{lem}\label{lgluw}
1. $j_*$ maps $\eu^{w\le 0}$ to $\cu^{w\le 0}$.

2. $j_!$ maps $\eu^{w\ge 0}$ to $\cu^{w\ge 0}$.

3. $i_*$ maps $\du^{w\le 0}$ (resp. $\du^{w\ge 0}$, resp. $\du^{w=
0}$) to $\cu^{w\le 0}$ (resp. to $\cu^{w\ge 0}$, resp. to $\cu^{w=
0}$).

\end{lem}
\begin{proof}
1. As we know, for any $Z\in \obj \eu$ we have $j^*j_*(Z)\cong Z$,
$i^!j_*Z=0$. Hence for $Z\in \obj \eu^{\le 0}$ we have $j^*j_*(Z)\in
\eu^{w\le 0}$, $i^!j_*Z\in \du^{w\le 0}$.

2. Similarly to the proof of assertion 1, it suffices to note that for
any $Z\in \obj \eu$ we have $j^*j_!(Z)\cong Z$, $i^*j_*Z=0$.

3. For any $Y\in \obj \du$ we have $i^!i_*(Y)\cong i^*i_*Y\cong Y$,
$j^*i_*Y=0$. Hence we obtain the statement for $\du^{w\le 0}$ and
$\du^{w\ge 0}$ (similarly to assertion 1). Since $\du^{w= 0}=\du^{w\le
0}\cap  \du^{w\ge 0}$, the last part of the assertion follows from
the previous ones immediately.

\end{proof}

\begin{theo}\label{tglu}
The classes $\cu^{w\le 0}$ and $\cu^{w\ge 0}$ define a weight
structure for $\cu$.

\end{theo}
\begin{proof}
The definitions of $\cu^{w\le 0}$ and $\cu^{w\ge 0}$ immediately
yield that they are Karoubi-closed and semi-invariant with respect
to translations (in the sense of Definition \ref{dwstr}).
They are also extension-stable (see Definition
\ref{exstab}).

The proof of orthogonality is similar to that in Theorem 1.4.10 of
\cite{BBD}. Let $X\in \obj \cu^{w\ge 1}$, $X'\in \cu^{w\le 0}$. We
should prove that $X\perp X'$.

Part (iv) of Definition \ref{dglu} yields a (long) exact sequence
$\dots\to \cu (i_*i^*X,X')\to \cu (X,X')\to \cu (j_!j^*X,X')
\to\dots$. It remains to note that the adjunctions of functors (of
 Definition \ref{dglu}(ii)) and the definitions of
$\cu^{w\le 0}$ and $\cu^{w\ge 0}$ yield that
$$\cu (i_*i^*X,X')=\du(i^*X, i^!X')=\ns=\eu(j^*X, j^*X')=\cu
(j_!j^*X,X') .$$

The last axiom check is the existence of weight decomposition (for any $X\in \obj \cu$).  

First we note that for any object of $i_*(\du)$ a weight decomposition exists
(since  one can take  those weight decompositions that come from
$\du$ via $i_*$). Thus if $X'\in \obj \cu$ possesses a weight decomposition then any $X\in \obj \cu$ such that $j^*(X)\cong j^*(X')$ possesses a weight decomposition also; here we use the definition of the localization (of $\cu$ by $i_*(\du)$) and 
Remark \ref{rtwd} (note that $i_*(\du)$ is shift-stable!).

For an $X\in \obj \cu$ we consider some weight decomposition
$B[-1]\to j^*X\to A\stackrel{f}{\to} B$ of $j^*X$ in $\eu$.
Applying the definition of the localization of $\cu$ by $i_*(\du)$
again, we obtain the existence $C,D\in \obj \cu$ and
$f'\in\cu(C,D)$ such that $j^*(C\stackrel{f'}{\to} D)\cong
(A\stackrel{f}{\to} B)$. Hence $j^*(\co f')[-1]\cong j^*(X)$. Thus it suffices to verify the existence of weight
decompositions for $C$ and $D[-1]$. By Lemma \ref{lgluw}(1,2),
'trivial' weight decomposition exist for $j_*A$ and $j_!B[-1]$.
We have $j^*j_*A\cong j^*C$ and $j^*j_!B[-1]\cong j^*D[-1]$ (since $j_!$ and $j^*$ are one-sided inverses of $j^*$); this concludes
the proof. 

\end{proof}

\begin{rema}\label{rglu}

1. The argument used for the proof of existence of weight
decompositions also yields that $\cu^{w\le 0}$ equals the
Karoubi-closure of the smallest extension-stable subclass of $\obj
\cu$ containing $\obj j_*(\eu^{w\le 0})\cup\obj i_*(\du^{w\le 0})$
(cf. the proof of  Theorem \ref{recw}(II1)); $\cu^{w\ge 0}$
equals the (similarly defined) 'envelope' of $\obj j_!(\eu^{w\ge
0})\cup \obj i_*(\du^{w\ge 0})$.

2. Using adjunctions and   Proposition \ref{bw}(1,2) we obtain that
$X\in \cu^{w\le 0}$ (resp. $X\in \cu^{w\ge 0}$) if and only if $j^*X\in
\eu^{w\le 0}$  and for any $Y\in \du^{w\ge 1}$ we have
$i_*Y\perp X$ (resp. $j^*X\in \eu^{w\ge 0}$ and for any $Y\in
\du^{w\le -1}$ we have $X\perp i_*Y$).

It follows that in order to 'calculate' $w$ for $\cu$, it suffices
to know only $i_*$ and $j^*$. Note that these two functors could be
easier to describe than the remaining ones (in particular,  for the
 categories of relative motives; see below). 

3. In \S2.3 of \cite{brelmot}and in \cite{bonivan}  gluing of weight structures was used  for the construction of (one of the versions of) the Chow weight structure for the
(triangulated) category $DM(S)$ of relative motives, i.e., of motives
over an excellent finite-dimensional  base scheme $S$ (that is  not a field). 

Certainly, the existence of the Chow weight structure  immediately yields weights
for arbitrary cohomology of relative motives. 

4. To the surprise of the author, it seems that one cannot glue
adjacent weight and $t$-structures in a compatible way (in the
general case). Indeed, suppose that one has $t$-structures on $\du$
and $\eu$ that are 
left
adjacent to the corresponding weight
structures. Following Theorem 1.4.10 of \cite{BBD}, one should
define $t$ on $\cu$ as follows: $\cu^{t\le 0}=\{X\in \obj \cu:\
i^*X\in \du^{t\le 0} ,\ j^*X\in \eu^{t\le 0} \}$ and $\cu^{t\ge
0}=\{X\in \obj \cu:\ i^!X\in \du^{t\ge 0} ,\ j^*X\in \eu^{t\ge 0}
\}$. It follows that $j_!$ maps $\eu^{t\le 0}$ to $\cu^{t\le 0}$;
$j_*$ maps $\eu^{t\ge 0}$ to $\cu^{t\ge 0}$. Now, suppose that $w$
is left adjacent to $t$ on $\cu$. In order to formulate (one half
of) the condition for $X\in\obj \cu$ to belong to $\cu^{w\ge 0}$ one
should translate the condition $X\perp j_!(\obj
\eu^{t\le 0})$ into certain condition of the type $\eu(F(X),Y) =\ns,\
\forall Y\in \obj \eu^{t\le 0}$ for a certain functor $F:\cu\to\eu$.
The problem is that $F$ should be left adjoint to $j_!$, i.e., it is
not a part of our gluing data
(cf. also  Proposition \ref{padad}(3,4)).
One would also have a similar (actually, a dual)
problem
for $w$ that is right adjacent to $t$.
 
5. One may also try to 'glue' coniveau spectral sequences (this
would correspond to gluing of the corresponding Gersten
weight structures as defined in \cite{bger}). Yet this could be
 difficult; see the problem
with gluing of adjacent structures described above.

\end{rema}



\subsection{Multiple compositions of $t$- and weight
truncations}\label{scompt}

Now suppose that $\cu$ is endowed with some weight structures $w_i$,
$1\le i\le m$, such that there exist left adjacent $t$-structures
$t_i$.

Then applying  Theorem \ref{sdt}(7,8) one can easily and
naturally express
 the functors represented by all possible compositions
of $t_i$-truncations as certain images (as in loc. cit).

For example, applying part 7 of loc. cit. twice we obtain
\begin{equation}\label{compt}
\cu(X,\tau_{1,\le i}(\tau_{2,\le j}Y)) \cong \imm(\cu(w_{1,\le
i}(w_{2,\le j}X),Y)\to \cu(w_{1,\le i+1}(w_{2,\le j+1}X,Y));
\end{equation}
for all $i,j\in \z$; this isomorphism is functorial in both $X$ and
$Y$. We recall that morphisms of objects can be (non-uniquely)
extended to morphisms of their weight decompositions; one can also
apply Proposition \ref{trfun}. In these formulae one can also shift
$t$-truncations by $[l],\ l\in\z$, and compose truncations from
different sides. For example, (\ref{virtttr}) is essentially a
formula of this sort.

Now we somewhat extend these results. Note that a  duality (see  Remark \ref{rsdt}(5)) may be given by $\Phi(X,Y)=\cu(X,F(Y))$,
where $F:\du\to \cu$ is an exact functor; in this case $A=\ab$. The
case of adjacent structures then corresponds to $F=\id_{\cu}$.

So, suppose that  $\cu$, $\du$, and $\eu$ are triangulated
categories, let $F:\du\to \cu$, $G:\eu\to \du$ be exact functors;
let $w_1$, $w_2$ be weight structures for $\cu$, $t_1$ be a
$t$-structure for $\du$ and $t_2$ be a $t$-structure for $\eu$,
suppose that they satisfy the following orthogonality conditions:
$\cu(\cu^{w_1\le 0}), F(\du^{t_1\ge 1}))=\cu(\cu^{w_1\ge 1},
F(\du^{t_1\le 0}))=\cu(\cu^{w_2\le 0}, F\circ G(\eu^{t_2\ge 1}))
=\eu(G\circ F (\cu^{w_2\ge 1}), \eu^{t_2\le 0})=0.$ Then one can
express $\cu(X, F(\tau_{1,\le i}G(\tau_{2,\le j}Y))$ for $Y\in \obj
\eu$ as a certain image similar to (\ref{compt}).

Certainly, one may also consider compositions of more than two exact
functors.

 \subsection{A strong weight complex functor for triangulated
categories that admit
 $f$-triangulated enhancements}\label{sftriang}

Now we check that the strong weight complex functor $t$ exists if
there exists an $f$-category enhancement of our category (we will
define this notion very soon); see  Remark \ref{rctst}(3).
 The argument below was kindly communicated to the author by prof.
 A. Beilinson; it is described in somewhat more detail in \S7 of \cite{schnur}.
To make our notation compatible with that of \cite{be87} we will
denote our basic triangulated category (which is usually $\cu$) by
$D$. As usual, $D$ is endowed with a weight structure $w$.

The plan of the construction is the following one. Suppose that
there exists an $f$-category $DF$ over $D$. In particular, this
yields the existence of the 'forgetting of filtration' functor
$\omega:\ DF\to D$. We describe a class of objects $DF^s\subset
\obj DF$  that satisfies the following conditions:

(i) any object of $X\in \obj D$ can be 'lifted' to an element of
$X^*\in DF^s$;

(ii) For every $M,N \in DF^s$ the map
   $DF(N,M)\to D(\omega(N),\omega (M))$
is surjective;

(iii) There exists a functor $e:DF\to C^b(D)$  such that
$e(DF^s)\subset C^b(\hw)$ and for any $M,N \in DF^s$ the
functor $e$ maps
   $\ke DF(N,M)\to D(\omega(N),\omega (M))$ to morphisms that
are homotopic to $0$.

We will denote the induced functor $DF\to K^b(D)$ by $e'$.

 Then $X\to e(X^*)$ yields an additive
functor $T:D\to K$, where $K$ is a certain triangulated category
isomorphic to $K^b(\hw)$. Indeed, by (ii) any two choices of $X'$
are connected by (possibly, non-unique) morphisms. By (iii) these
morphisms become canonical isomorphisms after the application of
$e'$.
 Hence it suffices to take $K$ being  the  category obtained from
 $K^b(\hw)$ by
 factorizing $\obj K^b(\hw)$ by these isomorphisms. Indeed, this family respects
 coproducts
 since $\omega$ and $e$ do.

Now we recall the relevant definitions of the Appendix of
\cite{be87}.

\begin{defi}\label{beil}

I A triangulated category $DF$ will be called a filtered
triangulated one if it is endowed with strict triangulated
subcategories $DF (\le 0)$ and $DF (\ge 0)$; an exact
autoequivalence $s:DF\to DF$; and a morphism of functors
$\al:\id_{DF}\to s$, such that the following axioms hold (for $DF
(\le n)= s^n(DF (\le 0))$ and $DF (\ge n)= s^n(DF (\ge 0))$).

(i)  $DF (\ge 1)\subset DF (\ge 0)$; $DF (\le 1)\supset DF (\le
0)$; $\cup_{n\in\z}DF (\ge n)= \cup_{n\in\z}DF (\le n)=DF$.

(ii) For any $X\in\obj DF$ we have $\alpha_X=s(\alpha_{s\ob X})$.

(iii) For any $X\in\obj DF(\ge 1)$ and $Y\in\obj DF(\le 0)$ we have
$DF(X,Y)=\ns$; whereas $\alpha$ induces an isomorphism $DF(Y,s\ob
X)\cong DF(sY,X)\to DF (Y,X)$.

(iv) Any  $X\in\obj DF$ can be completed to a distinguished
triangle $A\to X\to B$ with $A\in\obj DF(\ge 1)$ and $B\in\obj
DF(\le 0)$.

II $DF$ is called an f-category over $D$ if $D\subset DF;\ \obj D=
\obj DF (\le 0)\cap \obj DF (\ge 0)$.

III We will denote by $\om$ (see Proposition A3 of \cite{be87})
the only exact functor $DF\to D$  that satisfies the following conditions:

(i) Its  restrictions are right adjoint to the inclusion $D\to  DF
(\le 0)$ and  left adjoint to the inclusion $D\to  DF (\ge 0)$
respectively.

(ii)  $\om (\al_X)$ is an isomorphism.

(iii) $DF(X,Y)=D(\om X,\om Y)$ for any $X\in \obj DF(\le 0)$,
$Y\in \obj DF(\ge 0)$.

\end{defi}

A simple example of this axiomatics is described in Example A2
loc. cit.

By Proposition A3 loc. cit. there also exist exact functors
 $\sigma_{\ge n}:\ DF\to DF(\ge n)$, and $\sigma_{\le n}:\ DF\to DF(\le
 n)$ that are respectively right and left adjoint to the
 corresponding inclusions.
We denote $\gr^{[a,b]}:=\sigma_{\le b}\sigma_{\ge a}$,
 $ \gr^a=\gr^{[a,a]}$. Note that 
 there exist canonical
 and functorial (in $X$) morphisms $d: \si_{\le 0}X\to \si_{\ge 1}X[1]$
 that can be completed to a distinguished triangle in I(iv)
of Definition \ref{beil}.

Now we define $e$. For $M \in \obj DF$ the complex $e(M)$ has
components equal to $s^{-a} \gr^a M [a]$ (this lies in $\obj D
\subset \obj DF$), the differential will be equal to $s^{-a-1}
(s(d')\circ \al_{\gr^a}) [a]$; here $d'$ is the boundary map of
the canonical triangle $\gr^{a+1} \to \gr^{[a,a+1]}\to
\gr^a\stackrel{d'}{\to} \gr^{a+1}[1]$. $M$ is a complex indeed by
the axiom I(ii). We have $e(s(X))\cong e(X)$.

Now for a weight structure $w$ on $D$ we define $DB^s=e\ob
(C^b(\hw))$, i.e., we demand $\gr^a(X) \in s^a D^{w=a}$.

We will use the following statement.

\begin{lem}\label{ii} For every $M,N \in DF^s$ the map
   $\al_*DF(N,M)\to DF(N,s(M))$ is surjective; all $DF(N,s^a (M))\to
   DF(N,s^{a+1}(M))$ for $a>0$
   are bijective.
\end{lem}
\begin{proof}

Set $P= \co (\al_M:M\to s(M))$. By the long exact sequence for
$DF(N,-)$, it suffices to show that  $DF(N,s^a (P)[b])=\ns$ for $a+b
\ge 0$.
   Since $s^a M[-a]\in DF^s$, it
suffices to show that $DF(N,P[b])=\ns$ for $b\ge 0$.

    By devissage, we can assume that  $\gr^a M$ and  $\gr^b N$
vanish for $a\neq m$,
$b\neq n$, $m,n\in\z$. In other words, $M=s^m (K)[-m]$,  $N = s^n
(L)[-n]$ for some
  $K,L \in \cu^{w=0}\subset \obj D\subset \obj DF$.

  One has $DF (N,P[b])= D(L[-n], \omega (\sigma_{\ge n} P)[b])$.
To see that this group vanishes, consider 3 cases.

 (a) Suppose
that $n>m+1$. Then $\sigma_{\ge n} P=0$.

(b) Suppose $n\le m$. Then $\sigma_{\ge n} P =P$, so
       $\omega (\sigma_{\ge n} P)=\omega (P)=0$.

  (c) Suppose $n=m+1$. Then  $\sigma_{\ge n} P =s(M)$, so
$\omega (\sigma_{\ge n} P)= K[-m]$ and $D(N,P[b])= D(L[-n],
K[-m+b])=D(L,K[b+1])=\ns$ since $w$ is a weight structure.
\end{proof}

Now  (ii) follows from Lemma \ref{ii} immediately since for any
$X,Y\in \obj DF$ we have $DF(X,s(Y))\cong D(\om(X),\om
(s^nY))\cong D(\om(X),\om (Y))$ for $n$ large enough 
 Definition \ref{beil}(III(ii--iii)).

(ii) easily yields (i). Indeed, we can prove the statement for
$X\in D^{[i,j]}$ by the induction on $j-i$. We have obvious
inclusions $D^{w=i}\to DF^s$ (that split $\om$).

To make the inductive step it suffices to consider $X\in
D^{[0,m]}$ for $m>0$. Then $X^{w\le 0}$ and $X^{\ge 1}$ can be
lifted to $DF^s$ by the inductive assumption. The map $X^{w\le
0}\to X^{\ge 1}$ lifts to $DF^s$ by Lemma \ref{ii}; its cone will
belong to $DF^s$ and so will be a lift of $X$.

Now we verify (iii). By  Lemma \ref{ii}, for $M,N\in DF^s$ we have
$$\ke (DF(N,M)\to D(\omega(N),\omega (M)))=\ke(\al_*DF(N,M)\to
DF(N,s(M))).$$ Since $\om (M)\to \om(s(M))$ is an isomorphism, we
obtain (iii). Hence $T$ is a well-defined functor.

Now we note that $T$ is an 'enhancement' of our 'weak' weight complex
functor $t$. Indeed, $T$ and $t$ coincide on $\hw$; both of them
respect weight decompositions of objects and morphisms in a
compatible way.

Lastly,  to check that $T$ is an exact functor one should apply the
method of Remark \ref{virtri} and of the proof of
 Theorem \ref{wecomp}. Thus it suffices to lift any
distinguished triangle $C{\to} X{\to} X'$ so that the sequence
$e(C^*)\to e(X^*)\to e(X'{}^*)$ splits termwisely (in $C^b(\hw)$).
Now, to find such lifting it suffices to choose the weight
decompositions of $X$ and $X'$ arbitrarily; choose a weight
decomposition of $C$ as in the proof of  Theorem
\ref{wecomp}(I); and lift to $DB^s$ the map $t(C)\to t(X)$ as in the
proof of (i). Then the map $a^*:C^*\to X^*$ will become split
surjective after the application of each $\gr^a$; hence we can
choose $\co a^*\in DF^s$ as a lift for $X'$. This yields the lift
desired.

Hence $T$ is a strong weight complex functor for $D,w$. This
argument is a certain weight structure counterpart of Proposition
A5 of \cite{be87}.

It also seems possible that an $f$-category enhancement of $D$ would
allow  defining certain {\it higher truncation} functors; see Conjecture
\ref{cortrf}.

\subsection{Possible variations of the weight complex functor;
reduction modulo $p$}\label{svarwei}

Now we try to tell whether the main results of this
paper can be generalized to a more general setting. We cannot prove any
if and only if conditions; however we try to clarify the picture.
Since we include this subsection only to explain our choice of
definitions, it is rather sketchy.

First we study the question where do exact conservative functors
come from.

Suppose that $f:\cu\to \cu'$ is an exact functor (here $\cu,\cu'$
are triangulated categories). We denote by $\ke f$ the class of
morphisms that are mapped to $0$ by $f_*$. $\ke f$ is a (two-sided)
ideal of $\mo \cu$ (see Definition \ref{didmo}).

Obviously, $f$ is conservative if and only if 
$\id_X\not\in \ke f$ for any $X\in \obj \cu,\ X\neq 0$.
 Note that in this case $\ke f$ may be called a radical ideal since
 for any $X\in\obj \cu$, $s\in \cu(X,X)\cap \ke f$,
 $id_X+s$ will be an automorphism.

Now we study the following inverse problem: which  ideals can
correspond to conservative exact functors. Unfortunately, it seems
that  there does not exist a nice way to kill morphisms in an
arbitrary $I$ unless $\cu$ has a differential graded enhancement.
So we suppose that $\cu=\trp(D)$ for a differential graded
category $D$; an ideal $I\nde\mo \cu$ comes from a differential
graded nilpotent (or formally nilpotent in an appropriate sense)
ideal $I'$ of $\prtp D$. Then one can form a category
$\cu'=\trp(D/I')$; using a certain spectral sequence argument for
representable functors $X_*$ for $X\in\obj \cu$ similar to that
described below (for realizations)
 one can verify that the natural exact functor
$\cu\to \cu'$
is conservative. However one cannot hope for a spectral sequence for
a realization $H$ unless $H(I)$ belongs to some nice radical ideal
(probably more conditions are needed). Note that this condition is obviously
fulfilled for representable functors.

We describe one of the cases when it makes sense to construct such a
theory (and which does not come from a weight structure). Let
$\cu,D$ be pro-$p$-categories (i.e., the morphism group is an
abelian pro-$p$-group for any pair of objects), $\cu=\trp D$,
and $I'=p\mo(\cu)$. Let $H=\trp(E)$ for a differential graded
functor $E:D\to
B(pro-p-\ab)$, where $B(pro-p-\ab)$ is the 'big' category of
complexes of abelian profinite $p$-groups (see \S\ref{twmot}). Then the complex that computes $H(X)$ for $X\in\obj
\cu$ has a natural filtration by subcomplexes given by $p^iE$. These
subcomplexes correspond to the functors $\trp(p^iE)$, and the factors
of the filtration are quasi-isomorphic to the complexes calculating the
functors $F_i=\trp(p^iE/p^{i+1}E)$. It remains to note that $F_i$
can be factored through the natural functor $\cu\to \trp(D/p)$.
Hence in this case the spectral sequence of a filtered complex has
properties similar to those of the spectral sequence $S$ in 
\S\ref{twmot}; $F_i$ 
are similar to truncated realizations (see \S\ref{twmot} above and \S7.3
of \cite{mymot}).

\subsection{Three types of filtrations for triangulated categories;
 the relation of $w_{Chow}$ to the weight filtration
for motives}\label{stfiltr}

The author believes that the following three  types
of filtrations for triangulated categories  are equally important: $t$-structures, weight
structures and {\it horizontal structures} (those correspond to weight filtrations of \S1.1 of \cite{btrans}). Here a (left) horizontal
structure denotes a filtration of $C$ by full triangulated
subcategories $C_i$ such that for any $i$ the inclusion $C_i\to C$
admits a (left) adjoint (we call a filtration of this type
horizontal since it is shift-invariant). The relations between 
these types of filtrations are described in more detail in \cite{btrans}, where several
examples are also given. Note that any filtration of any of
the three types of described defines canonical functorial spectral
sequences for any (co)homology of objects.

It may be interesting to study various 'configurations' of
the structures of these types. In particular, we know that $\dmge$
and $\dmge\q$ support the corresponding Chow weight structures (see \S\ref{hchow}).
Conjecturally,
$\dmge\q$ should also support the   motivic $t$-structure (cf.
\S\ref{smm}) and the (horizontal) {\it weight filtration}.
Here the latter
is given by $C_i$ that are generated (as triangulated categories)
by (the homology)
$H^{MM}_j(X),\ X\in\obj \chowe\q,\ j\ge -i$.
Mixed motivic homology 'should be' strictly compatible with
\'etale (co)homology;
for $P\in \spv$ one should have $H_j^{MM}(\mg(P))=0$ for $j>0$; hence
$C_{-1}=\ns$. We
obtain that the weight filtration and $t_{MM}$ induce
the same filtration on $\chowe\q$. Moreover, on $MM=\hrt_{MM}$ the weight filtration
should induce the same filtration as the Chow weight structure.
Hence  $X$ should belong to $\in \dmge\q^{w\le 0}$ if and only if for
 all $j\in \z$ we have $H_{-j}^{MM}(X)\in C_{j}$.

In \cite{wildat} the (conjectural) picture described above was
justified (in the
case when $k$ is a number field) for the triangulated category
$DAT\subset\dmge\q$ (of so-called  Artin-Tate motives). It was also
 shown that
the restriction of $w_{Chow}$ to $DAT$ can be completely characterized
 in terms of weights of singular homology. Actually, this corresponds
 to the fact that the Beilinson's derived category of graded polarizable mixed Hodge
complexes can be endowed with
a weight structure and also with a 'classical weight filtration';
these filtrations and the canonical $t$-structure for this category are
connected by the same relations as those that should connect
the corresponding filtrations of motives (see \cite{btrans}). It can be easily seen
the singular 
realization respects weight structures mentioned (i.e., it is weight-exact);
it should also 'strictly respect' them (and this was essentially
proved in \cite{wildat} for Artin-Tate motives).

It may also  be interesting to study the relations between the
homotopy $t$-structure and the {\it slice filtration} on $\dme$ (and on other motivic categories); see \cite{hubka}
and \S4.9 of \cite{bger}.

\end{document}